\documentclass[11pt]{article}
\usepackage{amssymb,amsmath,setspace,graphicx,multicol,wrapfig,amsfonts,amsthm,titling,url,array,parskip,enumitem,bbm,color,tikz,tikz-cd,caption}
\usepackage{hyperref}
\usepackage{mathtools}
\usepackage{placeins,comment}
\usepackage[capitalize,nameinlink,noabbrev,nosort]{cleveref}
\usepackage{aurical,pbsi}
\usepackage[T1]{fontenc}
\usepackage[hmargin=2.5cm,vmargin=2.5cm]{geometry}
\usepackage{array, makecell,calculator,young}
\usepackage{subcaption}
\usepackage{float}
\usepackage[title]{appendix}

\usepackage[makeroom]{cancel}
\usepackage[backend=bibtex,natbib=true,giveninits,doi=false,url=false,isbn=false,eprint=false,maxbibnames=99]{biblatex}
\renewbibmacro{in:}{}

\graphicspath{{Figures/}}

\theoremstyle{plain}
\newtheorem{theorem}{Theorem}[section]
\newtheorem{prop}[theorem]{Proposition}
\newtheorem{lemma}[theorem]{Lemma}
\newtheorem{corollary}[theorem]{Corollary}
\newtheorem{conjecture}[theorem]{Conjecture}
\theoremstyle{definition}
\newtheorem{definition}[theorem]{Definition}
\newtheorem{defn}[theorem]{Definition}
\newtheorem{example}[theorem]{Example}
\newtheorem{question}[theorem]{Question}
\newtheorem{remark}[theorem]{Remark}
\crefname{defn}{Definition}{Definitions}
\crefname{theorem}{Theorem}{Theorems}
\crefname{lemma}{Lemma}{Lemmas}

\renewcommand{\to}{\longrightarrow}

\newcommand{\tcr}[1]{\textcolor{red}{#1}}
\newcommand{\tcb}[1]{\textcolor{blue}{#1}}

\newcommand\marginal[1]{\marginpar[\raggedleft\tiny #1]{\raggedright\tiny #1}}

\newcommand{\JV}[1]{\marginal{\textcolor{blue}{JV: #1}}}
\newcommand{\OM}[1]{\marginal{\textcolor{magenta}{OM: #1}}}
\newcommand{\TODO}[2][To do: ]{{\textcolor{red}{\textbf{#1#2}}}}

\newcommand*\mycirc[1]{%
  \begin{tikzpicture}
    \node[draw,circle,inner sep=1pt] {#1};
  \end{tikzpicture}}
\DeclareMathOperator{\TASEP}{TASEP}

\DeclareMathOperator{\wt}{wt}
\DeclareMathOperator{\rot}{rot}

\DeclareMathOperator{\South}{South}
\DeclareMathOperator{\dg}{dg}
\DeclareMathOperator{\Par}{\pi}
\DeclareMathOperator{\MLQ}{MLQ}
\DeclareMathOperator{\M}{\mathcal{M}_{(2)}}
\DeclareMathOperator{\GMLQ}{GMLQ}

\DeclareMathOperator{\Pair}{\pi}
\DeclareMathOperator{\inv}{\texttt{inv}}
\DeclareMathOperator{\maj}{\texttt{maj}}

\DeclareMathOperator{\arm}{arm}

\DeclareMathOperator{\MLD}{bMLQ}
\DeclareMathOperator{\cocharge}{\texttt{cocharge}}
 
\DeclareMathOperator{\leg}{leg}
\DeclareMathOperator{\charge}{\texttt{charge}}
\DeclareMathOperator{\cw}{cw}
\DeclareMathOperator{\quinv}{\texttt{quinv}}

\DeclareMathOperator{\SSYT}{SSYT}
\DeclareMathOperator{\collapse}{\rho_N}
\DeclareMathOperator{\crw}{crw}
\DeclareMathOperator{\rrw}{rrw}
\DeclareMathOperator{\revcrw}{\overline{\crw}}
\DeclareMathOperator{\colw}{\overline{\crw}}
\DeclareMathOperator{\rw}{rw}

\DeclareMathOperator{\rev}{rev}
\DeclareMathOperator{\Icol}{I_{col}}
\DeclareMathOperator{\Irow}{I_{row}}

\DeclareMathOperator{\mlq}{mlq}
\DeclareMathOperator{\tab}{tab}
\DeclareMathOperator{\mRSK}{mRSK}
\DeclareMathOperator{\dRSK}{bRSK}

\DeclareMathOperator{\nona}{non-attacking}
\DeclareMathOperator{\perm}{perm}
\DeclareMathOperator{\quord}{quinv-sorted}

\DeclareMathOperator{\LS}{LS}
\DeclareMathOperator{\flip}{\collapse^{\uparrow}}

\DeclareMathOperator{\words}{\mathcal{W}}
\DeclareMathOperator{\dualRSK}{dualRSK}

\newlength\cellsize 
\savebox2{%
\begin{picture}(12,12)
\put(0,0){\line(1,0){12}}
\put(0,0){\line(0,1){12}}
\put(12,0){\line(0,1){12}}
\put(0,12){\line(1,0){12}}
\end{picture}}

\newcommand\cellify[1]{\def\thearg{#1}\def\nothing{}%
\ifx\thearg\nothing
\vrule width0pt height\cellsize depth0pt\else
\hbox to 0pt{\usebox2\hss}\fi%
\vbox to 12\unitlength{
\vss
\hbox to 12\unitlength{\hss$#1$\hss}
\vss}}

\newcommand\tableau[1]{\vtop{\let\\=\cr
\setlength\baselineskip{-16000pt}
\setlength\lineskiplimit{16000pt}
\setlength\lineskip{0pt}
\halign{&\cellify{##}\cr#1\crcr}}}

\savebox3{%
\begin{picture}(12,12)
\put(0,0){\line(1,0){12}}
\put(0,0){\line(0,1){12}}
\put(12,0){\line(0,1){12}}
\put(0,12){\line(1,0){12}}
\end{picture}}

\newcommand\expath[1]{%
\hbox to 0pt{\usebox3\hss}%
\vbox to 12\unitlength{
\vss
\hbox to 12\unitlength{\hss$#1$\hss}
\vss}}

\begin{document}

\title{Macdonald polynomials at $t=0$ through twisted multiline queues}

\author{Olya Mandelshtam and Jer\'onimo Valencia-Porras}

\maketitle

\abstract{
Multiline queues are versatile combinatorial objects that play a key role in understanding the remarkable connection between the asymmetric simple exclusion process (ASEP) on a circle and Macdonald polynomials. Specializing the results of Corteel--Mandelshtam--Williams (2018) to the $t=0$ case yields a formula for the $q$-Whittaker polynomials through the Ferrari--Martin (2007) algorithm with a major index ($\maj$) statistic. In this paper, we reinterpret the $\maj$ statistic as a $\charge$ statistic on reading words, thereby bypassing the Ferrari--Martin algorithm to obtain an elegant formula for the $q$-Whittaker polynomials. Our methods naturally extend to the case of \emph{bosonic multiline queues}, with which we obtain analogous results for the modified Hall--Littlewood polynomials using a $\cocharge$ statistic on reading words.

\emph{Twisted multiline queues (GMLQs)} are obtained from the action of the symmetric group on the rows of a multiline queue. The Ferrari--Martin algorithm was extended to GMLQs by Arita--Ayyer--Mallick--Prolhac (2011), and Aas--Grinberg--Scrimshaw (2020) showed it is preserved under this action. We extend these results by defining a $\maj$ statistic on GMLQs that is also preserved under this action. This yields a novel family of formulas, indexed by compositions, for the $q$-Whittaker polynomials.  Additionally, we define a procedure on both GMLQs and bosonic multiline queues that we call \emph{collapsing}, which can can be realized via the Kashiwara (crystal) operators on type-A Kirillov--Reshetikhin crystals. As an application, we naturally recover the Lascoux--Sch\"utzenberger $\charge$ formula for the $q$-Whittaker and modified Hall--Littlewood polynomials, and the classical and dual Cauchy identities for Schur functions.}

\setcounter{tocdepth}{1}
\tableofcontents

\section{Introduction}

Recently, intriguing connections have emerged between one-dimensional integrable particle systems and Macdonald polynomials. A notable example is the connection between the \emph{asymmetric simple exclusion process (ASEP) on a circle} and the symmetric Macdonald polynomials $P_{\lambda}(X;q,t)$ \cite{CGW-2015}. This connection was made explicit through \emph{multiline queues}, which interpolate between the stationary probabilities of the particle process and the Macdonald polynomials \cite{CMW18}. Similarly, an analogous link was established between the \emph{totally asymmetric zero range process (TAZRP) on a circle} and modified Macdonald polynomials $\widetilde{H}_{\lambda}(X;q,t)$ \cite{AMM20}. 

The Macdonald polynomials $P_{\lambda} (X;q,t)$ \cite{Mac88} are a family of symmetric functions indexed by partitions, in the variables $X = \{x_1, x_2,\dots \}$, 
 with coefficients in $\mathbb{Q}(q,t)$. They are characterized as the unique monic basis for the ring of symmetric functions that satisfies certain triangularity and orthogonality conditions. 
They are a unifying structure playing a central role in algebraic combinatorics, specializing to many other important classes of symmetric functions such as Hall--Littlewood, $q$-Whittaker, Jack, and Schur polynomials. 

Significant attention has been devoted to studying Macdonald polynomials combinatorially. The celebrated result of Haglund, Haiman, and Loehr \cite{HHL05} yields a tableaux formula in weights that are rational functions in $q$ and $t$. At $t=0$, their formula specializes to obtain a formula for $q$-Whittaker polynomials through \emph{semistandard key tabloids} of partition shape with a $\maj$ statistic. Multiline queues were famously first introduced by Ferrari and Martin in \cite{FM07} to compute the probabilities of the \emph{multispecies totally asymmetric simple exclusion process (TASEP)}. The construction was generalized to incorporate the parameter $t$ by Martin in \cite{martin-2020}. Corteel, Mandelshtam, and Williams later established the combinatorial connection with Macdonald polynomials by enhancing a variant of Martin's multiline queues with additional statistics \cite{CMW18}. This resulted in the formula 
\begin{equation}\label{eq:P}
 P_{\lambda}(X;q,t) = \sum_{M}\wt(M)x^M,
\end{equation}
where the sum is over enhanced multiline queues with weights in parameters $(X,q,t)$. At $t=0$, these enhanced multiline queues are precisely the original Ferrari--Martin multiline queues, and $\wt(M)=q^{\maj(M)}$, where $\maj(M)$ is a statistic we define in \cref{subsec:mlqs} and is related to the classical $\maj$ statistic on tableaux in \cite{HHL08}. In this way, \eqref{eq:P} specializes to a multiline queue formula for the $q$-Whittaker polynomial:
\begin{equation}\label{eq:Phl}
    P_{\lambda}(X;q,0) = \sum_{M}q^{\maj(M)}x^M,
\end{equation}
where the sum is over Ferrari--Martin multiline queues with weights in parameters $(X,q)$.

Our main focus in this article is to develop the combinatorics of multiline queues for the $t=0$ case to study $P_{\lambda}(X;q,0)$. We describe a procedure akin to Robinson--Schensted (RS) insertion on these multiline queues that we call \emph{collapsing}. The procedure can be described using Kashiwara crystal operators on binary matrices. Using this procedure we recover the classical Schur expansion of the $q$-Whittaker polynomials (see \cite{bergeron2020survey} for related expressions):
\begin{equation}
    P_{\lambda}(X;q,0) = \sum_{\mu} K_{\mu' \lambda'}(q)s_\mu(X).
\end{equation}
\cref{thm:charge,thm:wrapping-to-nonwrapping} yield an elementary bijective proof of this formula through multiline queues. Moreover, we obtain multiline queue formulas for the $(q,t)$-Kostka polynomials at $t=0$:

\begin{equation}\label{eq:kostka}
    K_{\lambda\mu}(q,0) = \sum_{\substack{M\in\MLQ(\mu',\lambda') \\ \collapse(M) = M(\lambda')}}q^{\maj(M)} =\sum_{N\in\MLQ_0(\lambda,\rev(\mu))} q^{\maj(\rot(N))}.
\end{equation}
(See \cref{cor:KF_mlq_version_1,cor:KF_mlq_version_2}.)



While multiline queues are in bijection with binary matrices with partition row content, \emph{twisted multiline queues} are in bijection with the set of all binary matrices. Proving a conjecture in \cite{AAMP}, Aas, Grinberg, and Scrimshaw showed in \cite{AGS20} that they also encode the stationary distribution of the ASEP. Building on this work, we define a {major index statistic} on twisted multiline queues that gives a family of formulas, indexed by compositions, for the $q$-Whittaker polynomials:
\begin{equation}
    P_{\lambda}(X;q,0)=\sum_{M\in\GMLQ(\alpha,n)} q^{\maj_G(M)}x^M,
\end{equation}
for any composition $\alpha$ whose sorted rearrangement is $\lambda'$.

Finally, it is natural to define \emph{bosonic multiline queues}, which are analogues of multiline queues in bijection with integer matrices with finite support and weakly decreasing row sums. These objects give formulas for the \emph{modified Hall-Littlewood} polynomials with (modified) $\maj$ and $\cocharge$ statistics. Using the same framework, we define a collapsing procedure on bosonic multiline queues that allows us to recover a bijective proof of the Cauchy identity, and leads to analogous formulas for the $t=0$ specialization of the \emph{modified} $(q,t)$-Kostka polynomials in terms of bosonic multiline queues:

\begin{equation}\label{eq:kostka modified}
    \widetilde{K}_{\lambda\mu}(q,0) = \sum_{\substack{D\in\MLD(\mu,\lambda') \\ \widetilde\collapse(D) = D(\lambda')}}q^{\widetilde\maj(D)} =\sum_{D\in\MLD_{0,R}(\lambda',\mu)} q^{\widetilde\maj(\rot(D))}.
\end{equation}
(See \cref{cor:KF_mld_version_1,cor:KF_mld_version_2}.)

This article is organized as follows. In \cref{sec:preliminaries}, we give the background on tableaux, crystal operators, and the $\charge$ and $\cocharge$ statistics. In \cref{sec:mlqs-gmlqs}, we introduce multiline queues with the major index statistic, twisted multiline queues with their corresponding generalized major index, and bosonic multiline queues with a bosonic major index. 
In \cref{sec:collapse}, we define \emph{collapsing} as an insertion procedure on multiline queues and bosonic multiline queues to recover bijective proofs of the Schur expansions of the $q$-Whittaker and modified Hall--Littlewood polynomials. In particular, collapsing on twisted multiline queues yields a new family of formulas for the $q$-Whittaker polynomials in terms of the twisted objects. In \cref{sec:mlq/d_RSK}, we present multiline queue versions of the RSK correspondence for multiline queues and bosonic multiline queues, recovering proofs of the dual and classical Cauchy identities and yielding expressions for Kostka--Foulkes and modified Kostka--Foulkes polynomials in terms multiline queues and bosonic multiline queues, respectively.

\section{Preliminaries on partitions, fillings, and statistics}\label{sec:preliminaries}

A \emph{partition} $\lambda$ of a positive integer $n$ is a weakly decreasing sequence of nonnegative integers $\lambda_1\geq \lambda_2\geq\ldots\geq \lambda_\ell>0$ that sum to $n$. The numbers $\lambda_i$ are called the \emph{parts} of $\lambda$, the \emph{length}, denoted $\ell(\lambda)$, is the number of nonzero parts, and the \emph{size} is $|\lambda|=n$. The \emph{Young diagram} associated to a partition $\lambda$ consists of $\lambda_i$ left-justified boxes in the $i$-th row for $1\leq i\leq \ell(\lambda)$. To simplify notation, we write $\lambda$ to mean both the partition and its diagram. We use the convention that rows of the diagram are labeled from bottom to top, in accordance with French notation for Young diagrams. The \emph{conjugate} $\lambda'$ of a partition $\lambda$ is the partition obtained by reflecting $\lambda$ across the line $y=x$, with parts $\lambda_i'=|\{j\colon \lambda_j\geq i\}|$.  

We write $\mu\subseteq\lambda$ if $\mu_i\leq\lambda_i$ for all $i$. The \emph{dominance order} on partitions, denoted $\leq$, is defined by $\lambda\leq \mu$ if $\lambda_1+\cdots+\lambda_k\leq \mu_1+\cdots\mu_k$ for all $k$. 

A \emph{(weak) composition} $\alpha$ of a positive integer $n$ is an ordered tuple of nonnegative integers $(\alpha_1,\ldots,\alpha_\ell)$ such that the entries add up to $n$. We denote it by $\alpha\models n$. Also define $\alpha^+$ to be the partition obtained by rearranging the parts of $\alpha$ in weakly decreasing order. We call $\alpha$ a \emph{strong composition} if its parts are all positive.

For a partition $\lambda$, a \emph{filling} or \emph{tableau} of shape $\lambda$ is an assignment of positive integers to every box in the diagram of $\lambda$. A \emph{semistandard Young tableau} is a filling in which entries are strictly increasing from bottom to top and weakly increasing from left to right. Denote by $\SSYT(\lambda)$ the set of semistandard Young tableaux of shape $\lambda$. For $T\in\SSYT(\lambda)$, the \emph{content} of $T$ is the composition $c(T) = (c_1,c_2,\ldots)$ where $c_i$ is the number of occurrences of $i$ in $T$. We represent the content of a filling $T$ by the monomial $x^T \coloneqq \prod_{i\geq 1}x_i^{c_i}$. For a composition $\alpha$, we define $\SSYT(\lambda,\alpha)$ to be the set of semistandard Young tableaux of shape $\lambda$ with content $\alpha$. Finally, for a positive integer $n$, we define $\SSYT(\lambda,n)$ to be the set of semistandard Young tableaux of shape $\lambda$ whose entries are at most $n$.

Define $\words$ to be the set of words in the alphabet $\mathbb{N}$. 

\begin{defn}\label{def:readingorder}
    For a semistandard tableau $T$, define its \emph{row reading word} to be the word obtained from scanning the rows of $T$ from top to bottom and recording the entries from left to right within each row. Similarly, define its \emph{column reading word} to be the word obtained from scanning the columns of $T$ from left to right and recording the entries from top to bottom within each row. Denote the row and column reading words of $T$ by $\rrw(T)$ and $\crw(T)$, respectively.
\end{defn}

\begin{example}
    Let $\lambda = (6,4,3,2)$. We show the diagram of shape $\lambda$ with its cells labeled with respect to row reading order (left) and column reading order (right).
    \[{\tableau{1&2 \\ 3&4&5\\6&7&8&9\\10&11&12&13&14&15}} \qquad\qquad\qquad {\tableau{1&5 \\ 2&6&9\\3&7&10&12\\4&8&11&13&14&15}}
    \]
    We show a particular semistandard filling with content $\alpha=(3,2,2,0,3,0,4,1,0,0\ldots)$, i.e. with 
    $x^T=x_1^3x_2^2x_3^2x_5^3x_7^4x_8$, for the previous diagram of $\lambda$: 
    \[
    T \;\; = \quad \raisebox{15pt}{\tableau{7&8 \\ 5&7&7\\2&3&5&7\\1&1&1&2&3&5}}
    \]
    The row and column reading words, with $|$ separating the rows and columns respectively, are:
    \[\rrw(T)=7\,8\,|\,5\,7\,7\,|\,2\,3\,5\,7\,|\,1\,1\,1\,2\,3\,5\,\qquad\quad
    \crw(T)=7\,5\,2\,1\,|\,8\,7\,3\,1\,|\,7\,5\,1\,|\,7\,2\,|\,3\,|\,5.\]  
\end{example}

\subsection{Classical and cylindrical parentheses matching on words}\label{sec:operators}
We will make use of two types of related operations on words, described below. 

\begin{defn}[Bracketing rule]\label{def:bracketing}
Let $n$ be a positive integer and let $w$ be a word in the alphabet $\{1,\ldots,n\}$. For $1\leq i<n$, define $\Par_i(w)$ to be a word in open and closed parentheses $\{\,(\,\,,\,)\,\}$ that is obtained by reading $w$ from left to right and recording a "$($" for each $i+1$ and a "$)$" for each $i$. The \emph{bracketing rule} is the procedure of iteratively matching pairs of open and closed parentheses whenever they are adjacent or whenever there are only matched parentheses in between. Then $\Par_i(w)$ contains the data of which instances of $i$ and $i+1$ in $w$ are matched or unmatched following the bracketing rule applied to $\Par_i(w)$.
\end{defn}

In \cite{BumpSchilling17}, the bracketing rule is referred as the \emph{signature rule}. In the context of representation theory, it is a standard combinatorial tool to describe the action of raising and lowering operators on tensor products of crystals. We describe the action of these operators on words.

\begin{defn}[Raising and lowering operators]
Define the operator $E_i$ as follows. If $\Par_i(w)$ has no unmatched $i+1$'s, $E_i(w)=w$. Otherwise, $E_i(w)$ is $w$ with the leftmost unmatched $i+1$ changed to an $i$. Define the operator $F_i$ as follows. If $\Par_i(w)$ has no unmatched $i$'s, $F_i(w)=w$. Otherwise, $F_i(w)$ is $w$ with the rightmost unmatched $i$ changed to an $i+1$. Define $E_i^\star(w)$ to be the word $w$ with all unmatched $i+1$'s changed to $i$'s. 
\end{defn}

\begin{defn}[Cylindrical bracketing rule]\label{def:cylindrical matching}
Let $\Par^c_i(w)$ represent the word $\Par_i(w)$ on a circle, so that open and closed parentheses may match by wrapping around the word. Then the \emph{cylindrically unmatched} $i+1$'s and $i$'s in $w$ correspond respectively to the (cylindrically) unmatched open and closed parentheses in $\Par^c_i(w)$, according to the signature rule executed on a circle. The \emph{wrapping} $i+1$'s and $i$'s in $w$ correspond respectively to the cylindrically matched open and closed parentheses in $\Par^c_i(w)$ that are unmatched in $\Par_i(w)$.
\end{defn}


\begin{defn}[Reflections/Lascoux--Sch\"utzenberger involutions]
Define $S_i:\words\rightarrow\words$ as follows. If $\Par_i(w)$ has $a$ unmatched $i$'s and $b$ unmatched $i+1$'s, then $S_i(w)$ is the word obtained from $w$ by replacing the subword $i^a(i+1)^b$ with $i^b(i+1)^a$. 
\end{defn}

\begin{remark}\label{rem:reflections-content}
    If $w$ is a word with partition content $\alpha$, then $S_i(w)$ is a word with partition content $s_i\cdot\alpha$, where $s_i$ is the transposition swapping $i$ and $i+1$ and $\cdot$ represents the action of the permutation on indices. Explicitly, $s_i\cdot\alpha = (\alpha_1,\ldots,\alpha_{i-1},\alpha_{i+1},\alpha_i,\alpha_{i+2},\ldots).$
\end{remark}

\begin{remark}\label{rem:S_i-and-cylindrical-bracketing}
    The operators $S_i$ can be computed using the cylindrical bracketing rule: $S_i(w)$ is the word $w$ with all cylindrically unmatched $i$'s in $\Par^c_i$ converted to $i+1$'s, and vice versa. 
\end{remark} 

\begin{example}
    Consider the word $w=3\,1\,2\,2\,1\,4\,3\,4\,2\,1\,3\,1\,2\,3\,2$. The bracketing rule yields the following information, where the unmatched parentheses are show in red/bold and the $\_$ represent positions of the word that are ignored by this process: $$\Par_1(w)=\_\ \tcr{\boldsymbol{)}}\ (\ (\ )\ \_\ \_\ \_\ (\ )\ \_\ )\ \tcr{\boldsymbol{(}}\ \_\ \tcr{\boldsymbol{(}}$$ The unmatched 1's and 2's from $w$ are underlined: $3\,\underline{1}\,2\,2\,1\,4\,3\,4\,2\,1\,3\,1\,\underline{2}\,3\,\underline{2}$. Therefore, the action of the operators $E_1$, $F_1$ and $S_1$ in the word is 
    $$E_1(w) = 3\,\underline{1}\,2\,2\,1\,4\,3\,4\,2\,1\,3\,1\,\underline{\tcr{1}}\,3\,\underline{2} \qquad  F_1(w) = 3\,\underline{\tcr{2}}\,2\,2\,1\,4\,3\,4\,2\,1\,3\,1\,\underline{2}\,3\,\underline{2}$$ For the cylindrical bracketing rule we have: $$\Par_1^c(w)=\_\ )\ (\ (\ )\ \_\ \_\ \_\ (\ )\ \_\ )\ \tcr{\boldsymbol{(}}\ \_\ ($$ 
    Writing $w$ with a hat over the cylindrically unmatched elements (a subset of the underlined classically unmatched elements): $w=3\,\underline{1}\,2\,2\,1\,4\,3\,4\,2\,1\,3\,1\,\underline{\hat{2}}\,3\,\underline{2}$. The operator $S_1$ acts by switching the cylindrically unmatched elements: 
    \[S_1(w) = 3\,\underline{1}\,2\,2\,1\,4\,3\,4\,2\,1\,3\,1\,\underline{\hat{\tcr{1}}}\,3\,\underline{2}.
    \]
\end{example}

\subsection{$\charge$ and generalized $\charge$}
The $(q,t)$-Kostka polynomials are the polynomials originally appearing as the coefficients in the expansion
\begin{equation}\label{eq:Pschur} 
J_{\mu}(X;q,t) = \sum_{\lambda} K_{\lambda\mu}(q,t) s_{\lambda}[X(1-t)],
\end{equation}%
which can be equivalently written as
\begin{equation}\label{eq:Hschur} 
\widetilde{H}_{\mu}(X;q,t) = \sum_{\lambda} \widetilde K_{\lambda\mu}(q,t) s_{\lambda}(X).
\end{equation}
where $\widetilde K_{\mu\lambda}(q,t) = t^{n(\lambda)} K_{\mu\lambda}(q,t^{-1}).$
Although it is known due to \cite{Hai99} that $K_{\lambda\mu}(q,t)\in\mathbb{N}[q,t]$, combinatorial formulas for these coefficients are only known for general $\lambda,\mu$ in the $q=0$ case. This also yields a formula in the $t=0$ case since $K_{\lambda\mu}(q,t) = K_{\lambda'\mu'}(t,q)$. The first such result was given by Lascoux and Sch\"utzenberger \cite{LS78} in a seminal paper where they introduced the $\charge$ statistic.

\begin{defn}[Charge and cocharge]\label{charge}
The \emph{charge} of a permutation $\tau\in S_n$ is defined as
\[
\charge(\tau) 
=\sum_{\substack{i\in[n-1]\\\tau^{-1}_i<\tau^{-1}_{i+1}}} (n-i).
\]
which can be described as cyclically scanning the one-line notation of $\tau$ from left to right to read the entries in order from largest to smallest, and adding $n-i$ to  $\charge$ whenever one wraps around to reach $i$ from $i+1$. The $\cocharge$ is equal to $\cocharge(\tau)=\binom{n}{2}- \charge(\tau)$.
\end{defn}

The definition of $\charge$ generalizes to words with partition content by splitting the word into \emph{$\charge$ subwords}.

\begin{defn}\label{def:chargeT}
 Let $w$ be a word with partition content $\mu$. Extract the first subword $w^{(1)}$ by scanning $w$ from left to right and finding the first occurrence of its largest letter $k\coloneqq \mu'_1$, then $k-1,\ldots,2,1$, looping back around the word to the beginning whenever needed. This subword $w^{(1)}$ is then extracted from $w$, and the remaining $\charge$ subwords are obtained recursively from the remaining letters, which now have content $(\mu_1-1,\mu_2-1,\ldots,\mu_k-1)$. The subword $w^{(i)}$ can now be treated as a permutation in one-line notation. 
The $\charge$ of a word $w$ with partition content $\mu$ is the sum over its $\charge$ subwords $w^{(1)},\ldots,w^{(\ell)}$ where $\ell=\mu_1$:
\[
\charge(w)=\charge(w^{(1)})+\cdots+\charge(w^{(\ell)}).
\]
\end{defn}

The $\charge$ of a word with partition content can be equivalently described in terms of the classical and cylindrical matching operators on words, where the cylindrical matching operators pick out the sets of $\charge$ subwords of equal length, and the classical operators determine the contribution to $\charge$ from each set of $\charge$ subwords. 
\begin{lemma}\label{lem:charge with operators}
    Let $w$ be a word with content $\lambda$, with $L=\lambda_1$. Define $w_{(L)},w_{(L-1)},\ldots,w_{(1)}$ to be the decomposition of $w$ into (possibly empty) subwords, obtained as follows. Define $\overline{w}_{(L)}=w$, and for $r=L-1,\ldots,1$, let $\overline{w}_{(r)}=w\backslash\{w_{(L)}\cup\cdots\cup w_{(r+1)}\}$. Sequentially, for $r=L,L-1,\ldots,1$, set $w_{(r)}\coloneqq w_{(r,r)}\cup w_{(r,r-1)}\cup \cdots\cup w_{(r,1)}$ where $w_{(r,r)}$ is the subword consisting of the letters $r$ in $\overline{w}_{(r)}$, and for $k=r-1,\ldots,1$ $w_{(r,k)}$, let $\overline{w}_{(r,k+1)}$ be the subword of $\overline{w}_{(r)}$ that excludes all letters greater than $k$ that are not in $w_{(r,r)}\cup\cdots\cup w_{(r,k+1)}$. Then $w_{(r,k)}$ consists of the letters $k$ in $\overline{w}_{(r,k+1)}$ that are cylindrically matched in $\Par_k^c(\overline{w}_{(r,k+1)})$. Let $a_{r,k}$ be the number of letters $k$ that are cylindrically matched in $\Par_k^c(\overline{w}_{(r,k+1)})$, but are not classically matched in $\Par_k(\overline{w}_{(r,k+1)})$. Then the contribution to $\charge$ from the subword $w_{(r,k)}$ is $a_{r,k}(r-k)$, and the total $\charge$ of $w$ is
    \[
    \charge(w)=\sum_{2\leq r\leq L}\sum_{1\leq k\leq r-1} a_{r,k}(r-k).
    \]
\end{lemma}

\begin{proof}
    The equivalence to \cref{def:chargeT} is a straightforward check that we outline below, illustrated in \cref{ex:charge operators}
    \begin{itemize}
        \item For $r=L-1,\ldots,1$, $\overline{w}_{(r)}$ is the subword of $w$ after the $\charge$ subwords of length greater than $r$ have been extracted.
        \item For $2\leq r\leq L$, the letters ``$r$'' in the $\charge$ subwords of length $r$ are all the ``$r$'''s in $\overline{w}_{(r)}$, which by construction corresponds to the subword $w_{(r,r)}$. 
        \item For $k=r-1,\ldots,1$, the letters ``$k$'' in the $\charge$ subwords of length $r$ are the ``$k$'''s in $\overline{w}_{(r)}$ that are cylindrically matched to the $k+1$'s in $w_{(r,k+1)}$, which by construction corresponds to the subword $w_{(r,k)}$. 
        \item Moreover, the subset of $w_{r,k}$ that is not classically matched to the $k+1$'s in $w_{r,k+1}$ is precisely the set of letters $k$ in the length-$r$ subwords of $w$ that contribute to $\charge(w)$, and that contribution is $r-k$. 
    \end{itemize}
\end{proof}

\begin{example}\label{ex:charge operators}
Consider $w=3\,3\,4\,2\,2\,3\,2\,2\,1\,1\,1\,1\,1\,2\,3\,4$. We show the entries of the subwords $w_{(r,k)}$ circled with the entries \emph{not} in $\bar{w}_{(r)}$ replaced by ``$\cdot$''.
\begin{center}
\begin{tabular}{c|c|c|c|c|c|c}
$k$&$w_{(4,k)}$&$a_{4,k}$&$w_{(3,k)}$&$a_{3,k}$&$w_{(2,k)}$&$a_{2,k}$\\
4&$33\mycirc{4}223221111123\mycirc{4}$&&&&&\\
3&$\mycirc{3}3{4}22\mycirc{3}221111123{4}$&1&$\cdot\mycirc{3}\cdot\cdot2\cdot\cdot2\cdot\cdot1112\mycirc{3}\cdot$&&\\
2&${3}3{4}\mycirc{2}2{3}\mycirc{2}21111123{4}$&0&$\cdot{3}\cdot\cdot\mycirc{2}\cdot\cdot\mycirc{2}\cdot\cdot1112{3}\cdot$&1&$\cdot\cdot\cdot\cdot\cdot\cdot\cdot\cdot\cdot\cdot\cdot\cdot1\mycirc{2}\cdot\cdot$\\
1&${3}3{4}{2}2{3}{2}2\mycirc{1}\mycirc{1}11123{4}$&0&$\cdot{3}\cdot\cdot{2}\cdot\cdot{2}\cdot\cdot\mycirc{1}\mycirc{1}12{3}\cdot$&0&$\cdot\cdot\cdot\cdot\cdot\cdot\cdot\cdot\cdot\cdot\cdot\cdot\mycirc{1}{2}\cdot\cdot$&1
\end{tabular}
\end{center}
The contribution to $\charge(w)$ is thus $1(4-3)+0(4-2)+0(4-1)$ from $w_{(4)}$, $1(3-2)+0(3-1)$ from $w_{(3)}$, $1(2-1)$ from $w_{(2)}$, and 0 from $w_{(1)}$, which is empty. Thus $\charge(w)=3$.
\end{example}

Next we turn our attention to semistandard tableaux. 
For a partition $\mu$, define $n(\mu)=\sum_{i} \binom{\mu'_i}{2}.$

\begin{defn}
    Let $T$ be a semistandard tableau with partition content $\mu$. Then the \emph{charge of $T$} is defined as
    \[
    \charge(T)=\charge(\rrw(T))
    \]
    and the \emph{cocharge of $T$} is given by
    \[
    \cocharge(T)=n(\mu)-\charge(T).
    \]
\end{defn}

In fact, any reading order that is a linear extension of the northwest-to-southeast partial order (including the column reading order) can be used to compute the $\charge$ of a semistandard tableau. One way to see this is that any such reading word will produce the same Robinson-Schensted-Knuth (RSK) insertion tableau \cite{Fomin-Greene-1993}, implying Knuth equivalence, which in turn preserves $\charge$ \cite{butler1994subgroup}. We record this fact as a lemma to refer to it later on.
\begin{lemma}\label{lem:rw colw}
Let $T$ be a semistandard tableau with partition content. Then 
\[\charge(T)=\charge(\crw(T)).\]
\end{lemma}

\begin{example}
    For the semistandard tableau in $\SSYT(\lambda)$ for $\lambda=(7,6,3)$ below, the $\charge$ is computed as follows:
    \[
    T=\raisebox{7pt}{\tableau{3&3&5\\2&2&2&4&4&5\\1&1&1&1&2&3&4}}
    \]
The row reading word is $w=\rrw(T)=3\,3\,5\,|\,2\,2\,2\,4\,4\,5\,|\,1\,1\,1\,1\,2\,3\,4$. The $\charge$ words are $w^{(1)}=5\,2\,4\,1\,3$, $w^{(2)}=3\,2\,5\,1\,4$, $w^{(3)}=3\,2\,1\,4$, and $w^{(4)}=1\,2$. Then, the total $\charge$ is
\[
\charge(T)=\charge(5\,2\,4\,1\,3)+\charge(3\,2\,5\,1\,4)+\charge(3\,2\,1\,4)+\charge(1\,2)=3+2+1+1=7,
\]
and the cocharge is
\[
\cocharge(T)=n(\lambda)-\charge(T)=27-7=20.
\]
Alternatively, the column reading word is $v=\crw(T)=3\,2\,1\,|\,3\,2\,1\,|\,5\,2\,1\,|\,4\,1\,|\,4\,2\,|\,5\,3\,|\,4$, which has $\charge$ words $v^{(1)}=2\,1\,5\,4\,3$, $v^{(2)}=3\,2\,1\,5\,4$, $v^{(3)}=3\,2\,1\,4$, $v^{(4)}=1\,2$, having total $\charge$ also equal to $3+2+1+1=7$.
\end{example}

Using the $\charge$ statistic, the $q=0$ specialization of $K_{\lambda\mu}(q,t)$ is the generating polynomial for the $\charge$'s of the set of semistandard Young tableaux of shape $\lambda$ and with partition content $\mu$, namely $\SSYT(\lambda,\mu)$. As mentioned before, since $K_{\lambda\mu}(q,t) = K_{\lambda'\mu'}(t,q)$ (see \cite[Page 32]{MR2371044}), this yields the following Schur expansion for the $q$-Whittaker polynomials:
\begin{equation}\label{eq:P-LS}
P_{\mu}(X;q,0) = \sum_{\lambda} \left (\sum_{T\in \SSYT(\lambda',\mu')} q^{\charge(T)} \right ) s_{\lambda}(X).
\end{equation} 
The \emph{modified Kostka--Foulkes polynomials} are well-known to be related to $K_{\lambda\mu}(0,t)$ as 
$$\widetilde{K}_{\lambda\mu}(q,0)=\widetilde{K}_{\lambda\mu'}(0,q)=q^{n(\mu')}K_{\lambda\mu'}(0,q^{-1}) = q^{n(\mu')}K_{\lambda' \mu}(q^{-1},0),$$ and are thus given in terms of $\cocharge$:
\begin{equation}\label{eq:modified K-F}
\widetilde{K}_{\lambda\mu}(q,0) = \sum_{T\in \SSYT(\lambda,\mu')} q^{\cocharge(T)}. 
\end{equation}
The previous formulas give the Schur expansion for modified Hall-Littlewood polynomials:
\begin{equation}\label{eq:H-LS}
\widetilde{H}_{\mu}(X;q,0) = \sum_{\lambda} \left ( \sum_{T\in \SSYT(\mu,\lambda')} q^{\cocharge(T)} \right ) s_{\lambda}(X).
\end{equation}

Now we consider words with general composition content. The definition of $\charge$ can be extended to this case, by using reflection operators to \emph{straighten} the word in view of \cref{rem:reflections-content}, and then compute $\charge$ in the usual way. We will use this statistic in \cref{sec:GMLQ maj} to give a formula for $K_{\lambda\mu}(q,0)$ in terms of \emph{twisted multiline queues.}


\begin{defn}
    \label{def:chargeG}
    Let $\alpha$ be a (weak) composition with $n$ parts and let $w$ be a word with content $\alpha$. Let $\tau\in S_n$ be such that $\tau\cdot\alpha = \alpha^+$ and suppose it can be written as a reduced word $\tau = s_{i_k}\cdot s_{i_{k-1}} \cdots s_{i_1}$. The \emph{generalized $\charge$} of $w$ is defined as $$\charge_G(w) = \charge(S_{i_k}\circ S_{i_{k-1}} \circ\ldots\circ S_{i_1}(w)).$$
\end{defn}

\begin{remark}
    The permutation $\tau$ in the previous definition is not unique. However, if $\alpha_i = \alpha_{i+1}$, then $S_i$ acts trivially on a word with content $\alpha$. 
\end{remark}

\begin{example}
    Consider the word $w = 1\,4\,3\,3\,2\,1\,2\,4\,2\,4\,2$ with content $\alpha = (2,4,2,3)$. We show the steps in the straightening of $w$, underlining the letters that change in each step:
    \begin{align*}
        w = 1\,4\,3\,3\,2\,1\,2\,4\,2\,4\,2 \qquad &\longleftarrow \qquad \alpha = (2,4,2,3) \\
        S_1(w) = 1\,4\,3\,3\,2\,1\,\underline{1}\,4\,\underline{1}\,4\,2 \qquad &\longleftarrow \qquad s_1\cdot\alpha = (4,2,2,3) \\
        S_3\circ S_1(w) = 1\,4\,3\,3\,2\,1\,1\,\underline{3}\,1\,4\,2 \qquad &\longleftarrow \qquad s_3\cdot s_1\cdot\alpha = (4,2,3,2) \\
        S_2 \circ S_3\circ S_1(w) = 1\,4\,\underline{2}\,3\,2\,1\,1\,3\,1\,4\,2 \qquad &\longleftarrow \qquad s_2\cdot s_3\cdot s_1\cdot\alpha = (4,3,2,2) \\
    \end{align*}
    To compute the generalized $\charge$ of $w$ we use the straightened word and its corresponding $\charge$ subwords as follows:
    \begin{align*}
        \charge_G(w) &= \charge(1\,4\,2\,3\,2\,1\,1\,3\,1\,4\,2) \\ &= \charge(4\,3\,2\,1)+\charge(1\,3\,4\,2)+\charge(2\,1)+\charge(1) = 4.
    \end{align*}
\end{example}

\section{Multiline queues and twisted multiline queues}\label{sec:mlqs-gmlqs}

\subsection{Multiline queues}\label{subsec:mlqs}
Multiline queues were originally introduced by
Ferrari and Martin \cite{FM07} to compute the stationary probabilities for the \emph{multispecies totally asymmetric simple exclusion process (TASEP)} on a circle, building upon earlier work of Angel \cite{Angel06}. The TASEP is a one-dimensional particle process on a circular lattice describing the dynamics of interacting particles of different species (or priorities), in which each site of the lattice can be occupied by at most one particle, and particles can hop to an adjacent vacant site or swap places according to some Markovian process. The \emph{Ferrari--Martin (FM) algorithm} associates each multiline queue to a state of the TASEP via a \emph{projection map}; then, the stationary probability of each state is proportional to the number of multiline queues projecting to that state.

\begin{defn}\label{def:MLQ}
Let $\lambda$ be a partition, $L=\lambda_1$, and $n\geq \ell(\lambda)$ be an integer. A \emph{multiline queue of shape $(\lambda,n)$} is an arrangement of balls on an $L\times n$ array with rows numbered $1$ through $L$ from bottom to top, such that row $j$ contains $\lambda'_j$ balls. Columns are numbered $1$ through $n$ from left to right periodically modulo $n$, so that $j$ and $j+n$ correspond to the same column number. The site $(r,j)$ of $M$ refers to the cell in column $j$ of row $r$ of $M$. A multiline queue can be represented as a tuple $M=(B_1,\ldots,B_L)$ of $L$ subsets of $\{1,\ldots,n\}$ where the $j$-th subset has size $\lambda_j'$ and corresponds to the sites containing balls in row $j$. We denote the set of multiline queues of size $\lambda,n$ by $\MLQ(\lambda,n)$. At times we will omit specifying $n$, and just write $\MLQ(\lambda)$. Formally,
\[
\MLQ(\lambda,n) = \left\{ (B_1,\ldots,B_L):B_j\subseteq \{1,\ldots,n\},|B_j|=\lambda'_j\ \mbox{for}\ 1 \leq j \leq L
\right\}.
\]
\end{defn}

\begin{remark}
It is natural to regard a multiline queue as the tensor product of Kirillov--Reshetikhin crystals of single‐column shape.  This viewpoint is used in the presentation of multiline queues in \cite{KMO-multiline-2016} and will be mentioned frequently throughout this article. In particular, from this perspective, many of our definitions (cylindrical pairing in \cref{def:cylindrical matching,def:matched above below}, labeling rules in \cref{def:FM,def:mld labeling}, major index statistic in \cref{def:maj_mlqs,def:maj_mlds}) and results (\cref{lem:charge cw}, \cref{cor:dualCauchy,cor:Cauchy}, \cref{lem:localcrystalscommute_mlq,lem:localcrystalscommute_mld}) can be viewed as special cases or combinatorial proofs of classical statements about Kashiwara crystals.  We nonetheless provide a self‐contained setup with elementary proofs focused on the multiline queue setting.
\end{remark}

A multiline queue can be encoded  either by its row word or its column word.

\begin{defn}\label{def:words-mlq}
    For a multiline queue $M$, its \emph{row word} $\rw(M)$ is obtained by recording the column number of each ball in $M$ by scanning the rows from \emph{bottom to top} and from \emph{left to right} within each row. 
    Similarly, its \emph{column word} $\cw(M)$ is obtained by recording the row number of each ball by scanning the columns from \emph{left to right} and from \emph{top to bottom} within each column.
\end{defn}

\newcommand{\fr}{fr}

\begin{remark}\label{rem:biword}
    The words $\rw(M)$ and $\cw(M)$ are related through the interpretation of $M$ as \emph{biwords} of the associated binary matrix, where rows are labeled from bottom to top, and columns from left to right. Define the \emph{row biword} $B_r(M)$ to consist of entries $\binom{i}{j}$ for each pair $(i,j)$ such that $M_{n-i+1,j}=1$ (equivalently, the ball $j$ is in row $n-i+1$ of $M$), sorted such that $\binom{i}{j}$ is left of $\binom{i'}{j'}$ if and only if $i<i'$ or $i=i'$ and $j<j'$ (i.e. in lexicographic order). With this definition, $\rw(M)$ is the bottom word of $B_r(M)$. The \emph{column biword} $\overline{B_c}(M)$ consists of entries $\binom{j}{i}$for each pair $(i,j)$ such that $M_{n-i+1,j}=1$, sorted such that $\binom{j}{i}$ is left of $\binom{j'}{i'}$ if and only if $j<j'$ or $j=j'$ and $i>i'$ (i.e. in \emph{antilexicographic} order). Then $\cw(M)$ is the bottom row of $\overline{B_c}(M)$.
    
    Note that if $B_r(M) = \binom{w_1}{w_2}$ then $\overline{B_c}(M) = \binom{w_2'}{w_1'}$ where $w_1'$ and $w_2'$ are reorderings of $w_1$ and $w_2$ to match the antilexicographic order of the column biword. 
\end{remark}

\begin{example}
    For the multiline queue $M=(\{1,2,3,4\},\{1,3,5,6\},\{2,3\},\{3,5\})$ from \cref{ex:label} the words described in \cref{def:words-mlq} are $$\rw(M) = 1\,2\,3\,4\,|\,1\,3\,5\,6\,|\,2\,3\,|\,3\,5 \quad \text{and} \quad \cw(M) = 2\,1\,|\,3\,1\,|\,4\,3\,2\,1\,|\,1\,|\,4\,2\,|\,2$$
    As with tableaux, the bars ``$|$'' serve as delimiters for the rows and columns, respectively. In terms of biwords, we have
    \[B_r(M)=\binom{1\ 1\ 1\ 1\ 2\ 2\ 2\ 2\ 3\ 3\ 4\ 4}{1\ 2\ 3\ 4\ 1\ 3\ 5\ 6\ 2\ 3\ 3\ 5}\quad,\quad \overline{B_c}(M)=\binom{1\ 1\ 2\ 2\ 3\ 3\ 3\ 3\ 4\ 5\ 5\ 6}{2\ 1\ 3\ 1\ 4\ 3\ 2\ 1\ 1\ 4\ 2\ 2}.\]
\end{example}

The FM algorithm is a labeling procedure that deterministically assigns a label to each ball in a multiline queue $M$ to obtain a labeled multiline queue $L(M)$. 
We shall first give a description of the labeling procedure following \cite{KMO-multiline-2016}, through an iterative application of the cylindrical matching (NY) rule of \cref{def:cylindrical matching}.

\begin{defn}\label{def:matched above below}
Let $M=(B_1,\ldots,B_L)$ be a multiline queue. We call a ball in $B_r$ \emph{matched below} if it is paired in $\Par_r(\cw(M))$ for any $1\leq r<L$, and we call it \emph{matched above} if it is paired in $\Par_{r-1}(\cw(M))$ for any $1<r\leq L$. Otherwise, we call it \emph{unmatched below} (resp.~\emph{above}). Analogously, we say a ball is \emph{cylindrically matched/unmatched} by referring to $\Par_r^c(\cw(M))$.
\end{defn}

\begin{example}
    For the multiline queue $M=(\{1,2,3,4\},\{1,3,5,6\},\{2,3\},\{3,5\})$ in \cref{ex:label}, the balls matched below are at sites $(1,1), (1,2), (1,3),(1,4),(2,3),(2,5)$ and $(3,3)$, and the balls matched above are at sites $(2,1),(2,3),(3,2),(3,3)$, and $(4,3)$.
\end{example}

\begin{defn}[Multiline queue labeling process]\label{def:FM} Let $M=(B_1,\ldots,B_L)$ be a multiline queue of shape $(\lambda,n)$. Define the \emph{labeled multiline queue} $L(M)$ by replicating $M$ and sequentially labeling the balls, as follows. For each row $r$ for $r=L,\ldots,2$, each unlabeled ball in $B_r$ is labeled $r$. Next, for $\ell=L,\ldots,r$, let $\cw(M)^{(\ell,r)}$ be the restriction of $\cw(M)$ to the balls labeled ``$\ell\,$'' in $B_r$ and the unlabeled balls in $B_{r-1}$. The balls in row $r-1$ that are cylindrically matched in $\Par^c_{r-1}(\cw(M)^{(\ell,r)})$ acquire the label ``$\ell\,$''. 
To complete the process, all unpaired balls in row 1 are labeled ``$1$''. Such a labeling is shown in \cref{ex:label,ex:collapse}.
\end{defn}

\begin{remark}
    To distinguish between the multiline queue $M$ and its labeled version $L(M)$, we refer to the elements of $M$ as `balls' and the ones in $L(M)$ as `particles'. 
\end{remark}

\begin{defn}[Major index of a multiline queue]\label{def:maj}
Let $M\in\MLQ(\lambda,n)$ with labeling $L(M)$, and let $m_{\ell,r}$ be the number of particles labeled ``$\ell\,$'' in row $r$ of $L(M)$ that wrap when paired to the particles labeled ``$\ell\,$'' in row $r-1$. In other words, $m_{\ell,r}$ is the number of letters $r$ in $\cw(M)^{(\ell,r)}$, 
that are not matched above in $\Par_{r-1}(\cw(M)^{(\ell,r)})$. The \emph{major index} ($\maj$) of $M$ is:
\[
\maj(M)=\sum_{2\leq \ell\leq L}\ \sum_{2\leq r\leq \ell}m_{\ell,r}\,(\ell-r+1).
\]
\end{defn}

\begin{remark}
    We call this statistic $\maj$ in reference to the classical major index of a tableau, which is a sum over the legs of its descents. When a multiline queue is represented as a tableau (see \cref{sec:tableaux}), the wrapping pairings precisely correspond to the descents, making the major indices of the multiline queue and the tableau coincide.
\end{remark}

The FM algorithm \cite{FM07} is commonly described as a queueing process, in which for each row $i$, balls are paired between row $i$ and row $i-1$ one at a time 
with respect to a certain priority order. This form of the procedure produces the labeling $L(M)$ in addition to a \emph{set of pairings} between balls with the same label in adjacent rows. 

\begin{defn}[FM algorithm]\label{def:FM2}
Let $M=(B_1,\ldots,B_L)$ be a multiline queue of size $(\lambda,n)$. For each row $r$ for $r=L,\ldots,2$:
\begin{itemize}
\item Every unlabeled ball is labeled ``$r$''.
\item Once all balls in row $r$ are labeled, each of them is sequentially paired to the first unlabeled ball weakly to its right in row $r-1$, \emph{wrapping around from column $n$ to column 1} if necessary. The order in which balls are paired from row $r$ to row $r-1$ is from the largest label ``$L$'' to the smallest label ``$r$'', and (by convention) from left to right among balls with the same label. There is a unique choice for every such pairing. The resulting strands may be referred to as ``bully paths'' in the literature. 
\end{itemize}
To complete the process, all unpaired balls in row $1$ are labeled ``$1$''. We associate to $M$ the multiset $\Pair(M)=\{(r(p),\ell(p),\delta(p))\colon \, p\ \text{is a pairing in}\ M\}$ which records the following data for each pairing $p$: $r(p)$ is the row from which the pairing originates, $\ell(p)$ is the label corresponding to the pairing, and $\delta(p)$ is equal to $1$ if the pairing wraps and $0$ otherwise. See \cref{ex:label}. 
\end{defn}

The FM algorithm was originally introduced to compute the stationary distribution of the TASEP. Let $\TASEP(\lambda,n)$ be the set of TASEP states with particles of type $\lambda$ on $n$ sites. Namely, $\TASEP(\lambda,n)=S_n(\lambda_1,\ldots,\lambda_k,0^{n-k})$. Now define the \emph{FM projection map} $\mathfrak{p}:\MLQ(\lambda,n)\rightarrow\TASEP(\lambda,n)$ to be the word obtained by reading the labels of the bottom row of $L(M)$ from left to right. It was shown in \cite{FM07} that the stationary probability of a state of the TASEP is proportional to the number of multiline queues projecting to it.
\begin{theorem}[{\cite{FM07}}]
Let $\lambda$ be a partition and $n$ an integer. The stationary probability of a state $\mu\in\TASEP(\lambda,n)$ is equal to
\begin{equation}\label{eq:FM result}
\frac{1}{|\MLQ(\lambda,n)|}\Big| \{M\in\MLQ(\lambda,n)\colon\mathfrak{p}(M)=\mu\} \Big|.
\end{equation}
\end{theorem}

\begin{remark}\label{rem:order}
    The left-to-right order of pairing for balls of the same label is simply a convention: in fact, based on \cref{def:FM}, any pairing order among balls of the same type will yield the same labeling in the multiline queue (see, e.g. \cite[Lemma 2.2]{AGS20}). In particular, the multiset $\Pair(M)$ is invariant of the order of pairing of balls of the same label (see \cref{ex:label}). Of course, if one wishes to keep track of the \emph{strands} linking paired balls, those indeed depend on the pairing order, which becomes an important technical point when mapping multiline queues to tableaux; see, for instance, \cite[Section 5]{CMW18}.
\end{remark}

For a multiline queue $M=(B_1,\ldots,B_L)$, define its \emph{content} $x^M$ to be to be the monomial in $x_1,x_2,\ldots$ recording the numbers of balls contained in each column of $M$:
    \[x^M=\prod_{j=1}^L \prod_{b\in B_j} x_b.
    \]

    With $\Pair, \delta,\ell,r$ as given in \cref{def:FM2}, we can alternatively define $\maj(M)$ by 
\begin{equation}\label{def:maj_mlqs}
\maj(M)=\sum_{(r(p),\ell(p),\delta(p))\in\Pair(M)}\delta(p)(\ell(p)-r(p)+1).
\end{equation}

\begin{defn}\label{def:MLQ-weight}
    The \emph{weight} of a multiline queue is defined as $\wt(M)=x^M q^{\maj(M)}$. 
\end{defn}

The definition of $\wt(M)$ matches the $t=0$ restriction of the definition of the multiline queue weight in \cite{CMW18}.
    
\begin{example}\label{ex:label}
In \cref{fig:mlq_pairing}, we show a multiline queue $M=(\{1,2,3,4\},\{1,3,5,6\},\{2,3\},\{3,5\})$ of shape $(\lambda,n)$ with $\lambda=(4,4,2,2)$ and $n=6$. We use two different pairing orders to get the same labeled multiline queue $L(M)$. There are three wrapping pairings: one of label ``$4$'' from row 4, one of label ``$4$'' from row 2, and one of label ``$2$'' from row 2. The set of pairings of $M$ (which is independent from the pairing order) is given by 
\[
\Pair(M)=\{(4,4,0),(4,4,1),(3,4,0),(3,4,0),(2,4,0),(2,4,1),(2,2,0),(2,2,1)\}.
\]
Thus $\maj(M)=(4-4+1)+(4-2+1)+(2-2+1)=5$, so the weight of the multiline queue is $$\wt(M) = x_1^2x_2^2x_3^4x_4x_5^2x_6\,q^5.$$

\begin{figure}[!ht]
    \centering
    \captionsetup{width=.9\linewidth}
    \resizebox{.8\linewidth}{!}{
    \begin{tikzpicture}[scale=0.7]
    \def \w{1};
    \def \h{1};
    \def \r{0.25};
    
    \foreach \i in {0,...,4}
    {
    \draw[gray!50] (0,\i*\h)--(\w*6,\i*\h);
    }
    \foreach \i in {0,...,6}
    {
    \draw[gray!50] (\w*\i,0)--(\w*\i,4*\h);
    }
    \foreach \xx\yy in {2/3,4/3,1/2,2/2,0/1,2/1,4/1,5/1,0/0,1/0,2/0,3/0}
    {
    \fill (\w*.5+\w*\xx,\h*.5+\h*\yy) circle (\r cm);
    }
    
    
    \begin{scope}[xshift=7cm]
    \foreach \i in {0,...,4}
    {
    \draw[gray!50] (0,\i*\h)--(\w*6,\i*\h);
    }
    \foreach \i in {0,...,6}
    {
    \draw[gray!50] (\w*\i,0)--(\w*\i,4*\h);
    }
    \foreach \xx\yy\i in {2/3/4,4/3/4,1/2/4,2/2/4,0/1/2,2/1/4,4/1/4,5/1/2,0/0/4,1/0/2,2/0/4,3/0/2}
    {
    \draw (\w*.5+\w*\xx,\h*.5+\h*\yy) circle (\r cm);
    \node at (\w*.5+\w*\xx,\h*.5+\h*\yy) {\i};
    }

    \draw[black!50!green] (\w*2.5,\h*3.5-\r)--(\w*2.5,\h*2.5+\r);%
    \draw[blue,-stealth] (\w*4.5,\h*3.5-\r)--(\w*4.5,\h*2.95)--(\w*6.3,\h*2.95);%
    \draw[blue] (-.2,\h*2.95)--(\w*1.5,\h*2.95)--(\w*1.5,\h*2.5+\r);%
    
    \draw[black!50!green] (\w*2.5,\h*2.5-\r)--(\w*2.5,\h*2.05)--(\w*4.5,\h*2.05)--(\w*4.5,\h*1.5+\r);
    \draw[blue] (\w*1.5,\h*2.5-\r)--(\w*1.5,\h*1.87)--(\w*2.5,\h*1.87)--(\w*2.5,\h*1.5+\r);
    
    \draw[blue] (\w*2.5,\h*1.5-\r)--(\w*2.5,\h*.5+\r);
    \draw[black!50!green,-stealth] (\w*4.5,\h*1.5-\r)--(\w*4.5,\h*1.0)--(\w*6.3,\h*1.0);
    \draw[black!50!green] (-.2,\h*1.0)--(\w*0.5,\h*1.0)--(\w*0.5,\h*.5+\r);
    \draw[red,-stealth] (\w*5.5,\h*1.5-\r)--(\w*5.5,\h*.87)--(\w*6.3,\h*.87);
    \draw[red] (-.2,\h*.87)--(\w*3.5,\h*.87)--(\w*3.5,\h*.5+\r);
    \draw[orange] (\w*.5,\h*1.5-\r)--(\w*.5,\h*1.1)--(\w*1.5,\h*1.1)--(\w*1.5,\h*.5+\r);
    
    \end{scope}

    \begin{scope}[xshift=14cm]
    \foreach \i in {0,...,4}
    {
    \draw[gray!50] (0,\i*\h)--(\w*6,\i*\h);
    }
    \foreach \i in {0,...,6}
    {
    \draw[gray!50] (\w*\i,0)--(\w*\i,4*\h);
    }
    \foreach \xx\yy\i in {2/3/4,4/3/4,1/2/4,2/2/4,0/1/2,2/1/4,4/1/4,5/1/2,0/0/4,1/0/2,2/0/4,3/0/2}
    {
    \draw (\w*.5+\w*\xx,\h*.5+\h*\yy) circle (\r cm);
    \node at (\w*.5+\w*\xx,\h*.5+\h*\yy) {\i};
    }

    \draw[black!50!green] (\w*2.5,\h*3.5-\r)--(\w*2.5,\h*2.5+\r);%
    \draw[blue,-stealth] (\w*4.5,\h*3.5-\r)--(\w*4.5,\h*2.95)--(\w*6.3,\h*2.95);%
    \draw[blue] (-.2,\h*2.95)--(\w*1.5,\h*2.95)--(\w*1.5,\h*2.5+\r);%
    
    \draw[black!50!green] (\w*2.5,\h*2.5-\r)--(\w*2.5,\h*1.5+\r);
    \draw[blue] (\w*1.5,\h*2.5-\r)--(\w*1.5,\h*2.05)--(\w*4.5,\h*2.05)--(\w*4.5,\h*1.5+\r); 
    
    \draw[black!50!green] (\w*2.5,\h*1.5-\r)--(\w*2.5,\h*.5+\r);
    \draw[blue,-stealth] (\w*4.5,\h*1.5-\r)--(\w*4.5,\h*1.0)--(\w*6.3,\h*1.0);
    \draw[blue] (-.2,\h*1.0)--(\w*0.5,\h*1.0)--(\w*0.5,\h*.5+\r);
    \draw[red,-stealth] (\w*5.5,\h*1.5-\r)--(\w*5.5,\h*.87)--(\w*6.3,\h*.87);
    \draw[red] (-.2,\h*.87)--(\w*1.5,\h*.87)--(\w*1.5,\h*.5+\r);
    \draw[orange] (\w*.5,\h*1.5-\r)--(\w*.5,\h*1.1)--(\w*3.5,\h*1.1)--(\w*3.5,\h*.5+\r);
    
    \end{scope}
    \end{tikzpicture}
    }

    \caption{We show two different pairing orders giving the same labeling for the multiline queue on the left, where the first is the canonical left-to-right pairing order, and the second is reversed.}
    \label{fig:mlq_pairing}
\end{figure}
\end{example}

\begin{remark}
The equivalence of \cref{def:FM,def:FM2} is due to the order independence of pairing mentioned in \cref{rem:order}. Indeed, one can choose an order of pairing in the procedure from \cref{def:FM2} such that the pairings themselves correspond to the matching parentheses in $\Par_{r-1}^c(\cw(M)^{(\ell,r)})$ for each row $r$ and each label $r\leq \ell$ of balls in that row. Moreover, the coefficients $m_{\ell,r}$ in \cref{def:FM} correspond to the number of pairings $p\in\Pair(M)$ such that $(\ell(p),r(p),\delta(p))=(\ell,r,1)$, which implies that the two definitions of $\maj(M)$ are equivalent.
\end{remark}


%
Since both $L(M)$ and $\charge(M)$ can be defined through cylindrical matching on $\cw(M)$, we obtain the following characterization of $\maj(M)$ that bypasses the FM algorithm. Note that the result can be derived from the connection between the NY rule and the energy statistic on the corresponding KR crystals established in \cite{NY95}, but we give a concise, self-contained proof to highlight the role of the column reading word in a multiline queue. See \cref{rem:H} for further details. 

\begin{lemma}\label{lem:charge cw}
Let $M\in\MLQ(\lambda)$ be a multiline queue. Then
\[
\maj(M)=\charge(\cw(M)).
\]
\end{lemma}

\begin{proof}
Let $L$ be the height of $M$, let $w=\cw(M)$, and let $w_{(L)},w_{(L-1)},\ldots,w_{(1)}$ be the decomposition into subwords of $w$ according to the construction in \cref{lem:charge with operators}, so that $w_{(\ell)}$ contains all the length-$\ell$ $\charge$ subwords of $\cw(w)$, for $1\leq \ell\leq L$. 
Since the labels $L,L-1,\ldots,1$ in $L(M)$ and the subwords $w_{(L)},w_{(L-1)},\ldots,w_{(1)}$ in $\cw(M)$ are obtained sequentially using the same cylindrical matching rule, for each $\ell=L,L-1,\ldots,1$ 
the set of ``$\ell\,$''-labeled balls in 
$L(M)$ precisely corresponds to the subword $w_{(\ell)}$ 
in $\cw(M)$. 

Moreover, the contribution to $\charge(\cw(M))$ from $w_{(\ell)}$ is given by the difference between entries matched in $\Par_r^c(w_{(\ell)})$ and $\Par_r(w_{(\ell)})$ for $1\leq r<\ell$. Similarly, the contribution to $\maj(M)$ from the ``$\ell\,$''-labeled balls in $L(M)$ is given by the difference between balls matched in $\Par_r^c(\cw(M)^{(\ell,r)})$ and $\Par_r(\cw(M)^{(\ell,r)})$ for $1\leq r<\ell$. Every wrapping element that is cylindrically matched (from row $r+1$ to row $r$) contributes $\ell-r$ for a given label ``$\ell\,$'' and row $r$ for both $\charge(\cw(M))$ and $\maj(M)$, and since there is an equal number of such elements contributing in both statistics, we get that $\charge(\cw(M))=\maj(M)$ as desired.
\end{proof}

See \cref{ex:charge} for an example of \cref{lem:charge cw}. Notably, this theorem eliminates the need for the FM algorithm to determine $\maj(M)$. Thus we obtain the following formula for $P_{\lambda}(X;q,0)$.

\begin{corollary}\label{thm:P MLQ}
    Let $\lambda$ be a partition. The $q$-Whittaker polynomial is given by
    \begin{equation}\label{eq:P MLQ}
        P_{\lambda}(x_1,\ldots,x_n;q,0) = \sum_{M\in\MLQ(\lambda,n)} q^{\maj(M)}x^M = \sum_{M\in\MLQ(\lambda,n)} q^{\charge(\cw(M))}x^M,
        \end{equation}
        where the first equality is due to \cite{CMW18}.
\end{corollary}

\begin{example}\label{ex:charge}
    For the multiline queue in \cref{fig:mlq_pairing}, $\cw(M)=2\,1\,|\,3\,1\,|\,4\,3\,2\,1\,|\,1\,|\,4\,2\,|\,2$. The $\charge$ subwords of $\cw(M)$ are $w^{(1)}=4\,3\,2\,1$, $w^{(2)}=1\,3\,4\,2$, $w^{(3)}=2\,1$, $w^{(4)}=1\,2$. We indicate the entries of the $\charge$ subwords by the subscripts $1,2,3,4$:
    \[\cw(M)=2_31_23_21_34_13_12_11_11_44_22_22_4\]
    One sees that the entries with subscripts 1 and 2 correspond to the balls with label ``$4$'' in $M$, and those with subscripts 3 and 4 correspond to the balls with label ``$2$'' (as $|w^{(1)}|=|w^{(2)}|=4$ and $|w^{(3)}|=|w^{(4)}|=2$). The contribution to $\maj(M)$ from balls with label ``$4$'' is $1+3=4$ and the contribution to $\maj(M)$ from balls with label ``$2$'' is $1$. This matches the contribution to $\charge(\cw(M))$ from the length-4 subwords: $\charge(w^{(1)})+\charge(w^{(2)})=0+(1+3)=4$, and the length-2 subwords: $\charge(w^{(3)})+\charge(w^{(4)})=0+1=1$, respectively.
\end{example}

An important subset of multiline queues is the one where $\maj$ equal to zero, as the generating functions of those sets for fixed shapes correspond to Schur functions.

\begin{defn}
If $M$ satisfies $\maj(M)=0$, we call it \emph{nonwrapping}. Denote the set of nonwrapping multiline queues of shape $\lambda$ and size $n$ by $\MLQ_0(\lambda,n)$.
\end{defn}

Since nonwrapping multiline queues correspond to the $q=t=0$ restriction of \eqref{eq:P}, we have
\begin{equation}\label{eq:schur-mlqs}
s_{\lambda}(x_1,\ldots,x_n)=\sum_{M\in\MLQ_0(\lambda,n)} x^M. 
\end{equation}
There are two ways to see this bijectively. 
An immediate bijection from $\MLQ_0(\lambda,n)$ to $\SSYT(\lambda,n)$ is to build a semistandard tableau from $M\in\MLQ_0(\lambda,n)$ by sending a ball in $M$ at site $(r,j)$ to the content $n-j+1$ in row $r$ of the tableau. This bijection is not weight-preserving in the content monomials corresponding to $M$ and the tableau respectively, so one must rely on the symmetry of $s_{\lambda}$ to complete the argument. 
The second bijection, which is indeed weight-preserving, is given by (row or column) RS insertion of the row word of $M$, and is stated in \cref{sec:bijection_mlq0_ssyt}.

\subsection{A tableaux formula in bijection with multiline queues}\label{sec:tableaux}

The link between Macdonald polynomials and the asymmetric simple exclusion process (ASEP) is due to Cantini, De Gier, and Wheeler \cite{CGW-2015}, who found that the \emph{partition function} of the ASEP of type $\lambda$ on $n$ sites is a specialization  of the Macdonald polynomial $P_{\lambda}(x_1,\ldots,x_n;q,t)$ at $q=x_1=\cdots=x_n=1$. In \cite{martin-2020}, Martin introduced enhanced multiline queues to compute the stationary distribution for the ASEP with a parameter $t$; this parameter describes the relative rate of particles hopping left vs.~right (in the TASEP, $t=0$, which means particles only hop in one direction).  Building upon \cite{martin-2020,CGW-2015}, the multiline queue formula for $P_{\lambda}$ was given in \cite{CMW18}. At $t=0$, this formula coincides with \eqref{eq:P MLQ}. In particular, the $q$-Whittaker polynomial $P_{\lambda}(x_1,\ldots,x_n;q,0)$ analogously specializes to the partition function of the TASEP of type $\lambda$ on $n$ sites. 

The formula of \cite{CMW18} inspired the discovery of new a tableau formula for the modified Macdonald polynomials $\widetilde{H}_{\lambda}(X;q,t)$ in terms of a \emph{queue inversion statistic} ($\quinv$) \cite{AMM20, AMM22}, which is a statistic on tableaux that encodes multiline queue dynamics. In particular, this statistic established the connection between $\widetilde{H}_{\lambda}$ and the \emph{totally asymmetric zero range process}, whose relation to the ASEP is captured by the plethystic relationship between $\widetilde{H}_{\lambda}$ and $P_{\lambda}$. Using this relationship, the first author found a new tableau formula for $P_{\lambda}(X;q,t)$ in terms of the $\quinv$ statistic on \emph{$\quinv$-sorted, $\quinv$-non-attacking} fillings \cite{Man23,Man24}. 

It should be emphasized that although these new $\quinv$ formulas look very similar to the well-known Haglund--Haiman--Loehr formulas with the analogous statistic $\inv$ \cite{HHL05,HHL08}, they are fundamentally different in that the $\quinv$ statistic encodes properties of the pairings in the enhanced multiline queues with the parameter $t$. In the $t=0$ case, this formula yields a tableau representation of the FM pairing algorithm on multiline queues. To stay within the scope of this article, we only provide the formula in the $t=0$ case; for a full treatment, see \cite{Man24}.


For a partition $\lambda$, define $\dg(\lambda)$ to be the Young diagram of shape $\lambda'$, consisting of bottom-justified columns corresponding to the parts of $\lambda$. Let $x=(r,i), y=(r-1,i), z=(r-1,j)$ with $i<j$ be a triple of cells (if $\lambda_i=r-1$, then $x=\emptyset$ and the triple is called degenerate). For a filling $\tau:\dg(\lambda)\rightarrow \mathbb{Z}^+$, we call $(x,y,z)$ a \emph{$\quinv$ triple} if the entries are cyclically increasing (up to standardization) when read in counterclockwise order. That is, $\tau(x)<\tau(y)< \tau(z)$ or $\tau(y)< \tau(z)< \tau(x)$ or $\tau(z)<\tau(x)< \tau(y)$. The statistic $\quinv(\tau)$ counts the number of $\quinv$ triples in the filling $\tau$.

We also define the $\maj$ of a filling. For a cell $u=(r,c)\in\dg(\lambda)$, define $\leg(u)=\lambda_c-r$. Let $\South(u)$ be the cell directly below $u$ in $\dg(\lambda)$ (if it exists). Then define 
\begin{equation}\label{eq:maj tab}
\maj(\tau)=\sum_{\substack{u\in\dg(\lambda)\\\tau(u)>\tau(\South(u))}} \leg(u)+1.
\end{equation}
Then from \cite[Theorem 1.1]{Man24}, we obtain $P_\lambda(X;q,0)$ as a sum over fillings with maximal $\quinv$ equal to $n(\lambda)$:
\begin{equation}\label{eq:P 0}
    P_{\lambda}(x_1,\ldots,x_n;q,0)=\sum_{\substack{\tau:\dg(\lambda)\rightarrow [n]\\
    \quinv(\tau)=n(\lambda)}} q^{\maj(\tau)}x^{\tau}.
\end{equation}
\begin{remark}
The statement of \cite[Theorem 1.1]{Man24} includes the additional conditions that $\tau$ must be non-attacking and $\quinv$-sorted in the sum on the right hand side. However, it immediately follows from the definition of $\quinv$ that any filling $\tau$ with $\quinv(\tau)=n(\lambda)$ is necessarily both non-attacking and $\quinv$-sorted. Thus we may omit these two conditions.
\end{remark}

The right hand side of \eqref{eq:P 0} is a sum over all \emph{max-$\quinv$} fillings of $\dg(\lambda,n)$, which in particular implies there is a unique filling appearing in the sum for every given set of row contents. Moreover, in a max-$\quinv$ filling, the content of each row is a subset of $[n]$. Identifying such a filling with the (unique) multiline queue that has the same row contents, one gets the following correspondence, which can be deduced from \cite[Section 5]{Man24}.

\begin{theorem}\label{thm:queue to tableau}
    For a fixed $n$, there is a weight-preserving bijection between multiline queues in $\MLQ(\lambda,n)$ and max-$\quinv$ fillings $\tau:\dg(\lambda)\rightarrow [n]$, establishing a direct correspondence between the terms in the formulas \eqref{eq:P 0} and \eqref{eq:P MLQ}.
\end{theorem}

\begin{proof}
A content-preserving bijection between $\MLQ(\lambda,n)$ and max-$\quinv$ fillings of $\dg(\lambda)$ with entries in $[n]$ is given in \cite[Section 5]{Man24}, specialized at $t=0$. We will show this bijection preserves the $\maj$ statistic. In particular, for each row $r$ of a multiline queue $M\in\MLQ(\lambda,n)$ and each label $\ell$ of particles in that row, this bijection sends the particles with labels ``$\ell\,$'' to the cells in row $r$ in columns of height $\ell$ of a filling $\tau$ of $\dg(\lambda)$. The max-$\quinv$ condition on the filling captures the cylindrical pairing rule of the multiline queue. In particular, the number of descents in row $r+1$ of $\tau$ in columns of height $\ell$ (i.e. the cells whose content is greater than the content of the cell directly below) is precisely equal to the number of wrapping pairings of label ``$\ell\,$'' from row $r+1$ to row $r$ in $M$. The contribution from each of those cells in the set $\{u:\tau(u)>\tau(\South(u))\}$ is $\leg(u)+1=\ell-r+1$, which is equal to the contribution of the corresponding wrapping pairings to $\maj(M)$. Thus the total contribution to $\maj(\tau)$ in \eqref{eq:maj tab} is equal to $\maj(M)$.
\end{proof} 

\subsection{Twisted multiline queues}\label{sec:Gmlqs}
Relaxing the restriction on the row content of a multiline queue, we obtain \emph{twisted multiline queues}, which are in bijection with binary matrices with finite support. These were introduced in \cite{AAMP} and treated as operators on words of fixed length in a reformulation of a generalized FM algorithm. Our presentation of twisted multiline queues will follow the work of Aas, Grinberg, and Scrimshaw in \cite{AGS20}, where they consider multiline queues as a tensor product of Kirillov--Reshetikhin crystals with a so-called ``spectral weight''. In our context, this spectral weight is the content weight $x^M$ as defined in the previous section.  

\begin{defn}
Let $\lambda$ be a partition, $\alpha$ a composition such that $\alpha^+=\lambda'$, and $n\geq \ell(\lambda)$ a positive integer. 
A \emph{twisted multiline queue} of type $(\alpha,n)$ is a tuple of subsets $(B_1,\ldots,B_L)$ such that $B_j\subseteq[n]$ and $|B_j|=\alpha_j$ for $1 \leq j \leq L$. Denote the set of twisted multiline queues corresponding to a composition $\alpha$ by $\GMLQ(\alpha,n)$. 
\end{defn}

\begin{remark}
    According to the previous definition we have that $\MLQ(\lambda,n)=\GMLQ(\lambda',n)$. 
\end{remark}

A labeling procedure for $\GMLQ(\alpha,n)$ generalizing \cref{def:FM2} was introduced in \cite{AAMP} and reformulated in \cite[Section 2]{AGS20}, treating the components of the multiline queue as operators on words. In the procedure, vacancies in the multiline queue are considered to be \emph{anti-particles}, which are paired \emph{weakly to the left}. Labels are assigned sequentially to both the particles and the anti-particles by pairing sites between adjacent rows from top to bottom in a certain priority order such that particles (\emph{resp.}~anti-particles) are paired weakly to the right (\emph{resp.}~left), and propagating the labels upon pairing. When $B$ is a multiline queue, the labeling of the particles coincides with the FM algorithm.

\begin{defn}[Generalized multiline queue pairing process]\label{def:gmlq_particlewise} Let $\alpha=(\alpha_1,\ldots,\alpha_L)$ be a (weak) composition with $\lambda=\alpha^+$, let $n\geq \ell(\lambda')$, and let $B=(B_1,\ldots,B_L)\in\GMLQ(\alpha,n)$. The algorithm will produce a labeling of each site of the $n\times L$ lattice such that each particle and each anti-particle in $B$ has an associated label, and we will denote this labeled configuration by $L_G(B)$. 
\begin{itemize}
    \item Begin by labeling each particle in the topmost row by ``$L$'' and each anti-particle in the topmost row by ``$L-1$''. 
    \item Sequentially, for $r=L-1,L-2,\ldots,1$, do the following. Assuming row $r+1$ has been labeled in the previous step, let $w=(w_1,\ldots,w_n)$ be the set of labels read off row $r+1$ from left to right. The labeling process of row $r$ occurs in two independent phases. Let $(i_1,\ldots,i_n)\in S_n$ be the shortest permutation (with respect to Bruhat order) that fixes a weakly decreasing order $w_{i_1}\geq w_{i_2}\geq \cdots \geq w_{i_n}$ on the elements of $w$, namely if $w_{i_k}=w_{i_{k+1}}$, then $i_k<i_{k+1}$. Let $s=|B_r|$ be the number of particles in row $r$.  
\begin{enumerate}
\item {\bf Particle Phase.}  
For $k=1,2,\ldots,s$, find the first unlabeled particle in row $r$ weakly to the right (cyclically) of site $i_k$ and label it ``$w_{i_k}$''. 
\item {\bf Anti-particle Phase.}  
For $k=n,n-1,\ldots,s+1$, find the first unlabeled anti-particle weakly to the left (cyclically) of site $i_k$ and label it ``$w_{i_k}-1$''. 
\end{enumerate}
\end{itemize}
\end{defn}

We note that during the pairing/labeling process of each row, the outcome is independent of the order of pairing chosen among sites with the same label. However, the condition that the pairing order permutation is the shortest with respect to Bruhat order implies that sites with the same label are ordered from left to right, as in \cref{ex:GMLQ pairing}.

\begin{example}\label{ex:GMLQ pairing}
Let $B$ be a twisted multiline queue on $6$ columns such that $w = 2\,5\,2\,3\,4\,2$ is the labeling of row $i+1$ of $L_G(B)$, and let $B_i=\{1,5\}$ be the queue at row $i$. We show the labeling of row $i$ in $L_G(B)$ according to \cref{def:gmlq_particlewise}.

\begin{figure}[h!]
\begin{center}
\captionsetup{width=.8\linewidth}
\resizebox{0.8\linewidth}{!}{
    \begin{tikzpicture}[scale=0.7]
    \def \w{1};
    \def \h{1};  
    \def \r{0.3};
    \def \off{0.09};

    \begin{scope}[xshift=0cm]
    
    \foreach \xx\yy\i in {0/1,4/1}
    {
    \draw[blue] (\w*.5+\w*\xx,\h*.5+\h*\yy) circle (\r cm);
    }
    \foreach \xx\yy\i in {1/1,2/1,3/1,5/1}
    {
    \draw[red] (\w*.5+\w*\xx-\w*\r,\h*.5+\h*\yy+\h*\r) -- (\w*.5+\w*\xx+\w*\r,\h*.5+\h*\yy+\h*\r) -- (\w*.5+\w*\xx+\w*\r,\h*.5+\h*\yy-\h*\r)--(\w*.5+\w*\xx-\w*\r,\h*.5+\h*\yy-\h*\r)--(\w*.5+\w*\xx-\w*\r,\h*.5+\h*\yy+\h*\r);
    }
    \node[black] at (0.5,2.5) {2};
    \node[black] at (1.5,2.5) {5};
    \node[black] at (2.5,2.5) {4};
    \node[black] at (3.5,2.5) {2};
    \node[black] at (4.5,2.5) {4};
    \node[black] at (5.5,2.5) {2};
    \node[gray] at (0.5,3.1) {\tiny 4};
    \node[gray] at (1.5,3.1) {\tiny 1};
    \node[gray] at (2.5,3.1) {\tiny 2};
    \node[gray] at (3.5,3.1) {\tiny 5};
    \node[gray] at (4.5,3.1) {\tiny 3};
    \node[gray] at (5.5,3.1) {\tiny 6};

    \node[blue] at (0.5,1.5) {4};
    \node[red] at (1.5,1.5) {3};
    \node[red] at (2.5,1.5) {1};
    \node[red] at (3.5,1.5) {1};
    \node[blue] at (4.5,1.5) {5};
    \node[red] at (5.5,1.5) {1};

        \draw[blue] (\w*1.5-\off,\h*2.5-\r)--(\w*1.5-\off,\h*1.95)--(\w*4.5-\off,\h*1.95)--(\w*4.5-\off,\h*1.5+\r);
    
        \draw[blue,-stealth] (\w*2.5-\off,\h*2.5-\r)--(\w*2.5-\off,\h*2.1)--(\w*6.2,\h*2.1);
        \draw[blue] (-.2,\h*2.1)--(\w*0.5-\off,\h*2.1)--(\w*0.5-\off,\h*1.5+\r);
       
        \draw[red,-stealth] (\w*0.5+\off,\h*2.5-\r)--(\w*0.5+\off,\h*2.0)--(-0.2,\h*2.0);
        \draw[red] (6.2,\h*2.0)--(\w*5.5+\off,\h*2.0)--(\w*5.5+\off,\h*1.5+\r);

        \draw[red] (\w*5.5-\off,\h*2.5-\r)--(\w*5.5-\off,\h*2.025)--(\w*2.5+\off,\h*2.025)--(\w*2.5+\off,\h*1.5+\r);

        \draw[red] (\w*3.5,\h*2.5-\r)--(\w*3.5,\h*1.5+\r);

        \draw[red] (\w*4.5+\off,\h*2.5-\r)--(\w*4.5+\off,\h*1.875)--(\w*1.5+\off,\h*1.875)--(\w*1.5+\off,\h*1.5+\r);
        
    \end{scope}


    \begin{scope}[xshift=9cm]
    
    \foreach \xx\yy\i in {0/1,4/1}
    {
    \draw[blue] (\w*.5+\w*\xx,\h*.5+\h*\yy) circle (\r cm);
    }
    \foreach \xx\yy\i in {1/1,2/1,3/1,5/1}
    {
    \draw[red] (\w*.5+\w*\xx-\w*\r,\h*.5+\h*\yy+\h*\r) -- (\w*.5+\w*\xx+\w*\r,\h*.5+\h*\yy+\h*\r) -- (\w*.5+\w*\xx+\w*\r,\h*.5+\h*\yy-\h*\r)--(\w*.5+\w*\xx-\w*\r,\h*.5+\h*\yy-\h*\r)--(\w*.5+\w*\xx-\w*\r,\h*.5+\h*\yy+\h*\r);
    }
    \node[black] at (0.5,2.5) {2};
    \node[black] at (1.5,2.5) {5};
    \node[black] at (2.5,2.5) {4};
    \node[black] at (3.5,2.5) {2};
    \node[black] at (4.5,2.5) {4};
    \node[black] at (5.5,2.5) {2};
    \node[gray] at (0.5,3.1) {\tiny 6};
    \node[gray] at (1.5,3.1) {\tiny 1};
    \node[gray] at (2.5,3.1) {\tiny 3};
    \node[gray] at (3.5,3.1) {\tiny 4};
    \node[gray] at (4.5,3.1) {\tiny 2};
    \node[gray] at (5.5,3.1) {\tiny 5};

    \node[blue] at (0.5,1.5) {4};
    \node[red] at (1.5,1.5) {3};
    \node[red] at (2.5,1.5) {1};
    \node[red] at (3.5,1.5) {1};
    \node[blue] at (4.5,1.5) {5};
    \node[red] at (5.5,1.5) {1};

        \draw[blue] (\w*1.5-\off,\h*2.5-\r)--(\w*1.5-\off,\h*1.9)--(\w*4.5-\off,\h*1.9)--(\w*4.5-\off,\h*1.5+\r);
    
        \draw[blue,-stealth] (\w*4.5-\off,\h*2.5-\r)--(\w*4.5-\off,\h*2.1)--(\w*6.2,\h*2.1);
        \draw[blue] (-.2,\h*2.1)--(\w*0.5-\off,\h*2.1)--(\w*0.5-\off,\h*1.5+\r);
       
        \draw[red,-stealth] (\w*0.5+\off,\h*2.5-\r)--(\w*0.5+\off,\h*2.0)--(-0.2,\h*2.0);
        \draw[red] (6.2,\h*2.0)--(\w*5.5+\off,\h*2.0)--(\w*5.5+\off,\h*1.5+\r);

        \draw[red] (\w*5.5-\off,\h*2.5-\r)--(\w*5.5-\off,\h*2.0)--(\w*3.5+\off,\h*2.0)--(\w*3.5+\off,\h*1.5+\r);

        \draw[red] (\w*3.5-\off,\h*2.5-\r)--(\w*3.5-\off,\h*2.0)--(\w*2.5+\off,\h*2.0)--(\w*2.5+\off,\h*1.5+\r);

        \draw[red] (\w*2.5-\off,\h*2.5-\r)--(\w*2.5-\off,\h*2.0)--(\w*1.5+\off,\h*2.0)--(\w*1.5+\off,\h*1.5+\r);  
    \end{scope}
    \end{tikzpicture}
    }

\caption{We show the pairings according to \cref{def:gmlq_particlewise}, with two different pairing orders (producing the same labeling). The small grey numbers above the word $w$ show the permutation giving the order of pairing priority among the letters of $w$.} 
\end{center}
\end{figure}
\end{example}

\begin{remark}\label{rem:only-word-is-needed-GMLQs}
    Note that in order to assign labels in row $i$ from a labeled row $i+1$, we don't need to know the set $B_{i+1}$. Indeed, from the description of the pairing algorithm from \cref{def:gmlq_particlewise}, only the word $w$ is needed, as shown in the previous example. 
\end{remark}

Let $\alpha$ be a composition and let $\lambda=\alpha^+$. We may extend the definition of the FM projection map to the generalized projection map $\mathfrak{p}:\GMLQ(\alpha,n)\rightarrow\TASEP(\lambda',n)$ to be the word obtained by reading the labels of the bottom row of $L_G(M)$ from left to right.

\begin{remark}
    The pairing process can alternatively be described in terms of a \emph{bicolored reading word} that keeps track of both the particles and the anti-particles, through the functions $\Par^c_i$ sequentially applied to this word. We omit this description, as it is a straightforward generalization of \cref{def:FM}.
\end{remark}

\begin{lemma}\label{lem:antiparticles have smaller label}
Let $B=(B_1,\ldots,B_L)\in\GMLQ(\alpha,n)$ be a twisted multiline queue with labeling $L_G(B)$. 
\begin{itemize}
\item[(i)] Within each row of $L_G(B)$, the largest label of an anti-particle is strictly smaller than the smallest label of a particle.
\item[(ii)] If $\alpha=\alpha^+$, $L_G(B)$ restricted to the particles coincides with $L(B)$ from \cref{def:FM2}, and for each $1 \leq r \leq L$, the entries in the $r$'th row of $L_G(B)$ corresponding to anti-particles are labeled $r-1$.
\end{itemize}
\end{lemma}

\begin{proof}
    The proof for (i) is by induction on the rows $r$ of $B$. It holds trivially for the base case $r=L$. Assuming the claim holds for row $r+1$ for $1\leq r< L$, let ``$k$'' be the smallest label of a site that pairs during the particle phase. Then every particle in row $r$ will have a label that is at least $k$, and every anti-particle will have a label that is at most $k-1$ since labels are decremented by 1 during the anti-particle phase. Thus the claim holds for row $r$ as well.
 
    We again prove (ii) by induction on the rows $r$ of $M$. Since particles are labeled ``$L$'' in both $L_G(B)$ and $L(B)$ and anti-particles are labeled ``$L-1$'' in $L_G(B)$, the base case $r=L$ holds. Suppose the claim holds for row $r+1$ for $1\leq r<L$. Since $|B_r|\geq |B_{r+1}|$, every particle in row $r+1$ pairs during the particle phase, which coincides with the ordinary FM particle pairing process. To complete the particle phase, the leftmost $|B_r|-|B_{r+1}|$ anti-particles in row $r+1$ (which are labeled $r$ by the induction hypothesis) pair to the remaining unpaired particles in row $r$, giving each of them the label ``$r$'', the same label they would have received in $L(B)$. In the anti-particle phase, the remaining anti-particles in row $r+1$ pair to the anti-particles in row $r$, giving each of them the label ``$r-1$''. Thus the statement holds for row $r$ as well.
\end{proof}

There is a natural symmetric group action on the rows of twisted multiline queues, generated by row-swapping involutions $\sigma_i$ that establish an isomorphism $\GMLQ(\alpha)\rightarrow\GMLQ(s_i\cdot\alpha)$.

\begin{defn}
    For $B=(B_1,\ldots,B_L)\in\GMLQ(\alpha)$, for $1\leq i\leq L-1$, define the involution $\sigma_i:\GMLQ(\alpha)\rightarrow\GMLQ(s_i\cdot\alpha)$ to swap the numbers of particles between rows $i+1$ and $i$ 
    as follows. If $|B_{i+1}|=|B_i|$, $\sigma_i(B)=B$. Otherwise, $\sigma_i(B)$ is obtained by exchanging cylindrically unmatched particles in $\Par^c_i(\cw(B))$ between $B_{i+1}$ and $B_i$.  
    \end{defn}

\begin{example}\label{ex:sigma}
    For $\alpha=(2,2,3)$ and $B=(\{2,3\},\{1,4\},\{2,3,4\})\in\GMLQ(\alpha,4)$, we show $\sigma_2(B)=(\{2,3\},\{1,2,4\},\{3,4\})\in\GMLQ(s_2\cdot\alpha,4)$ and $\sigma_1(\sigma_2(B))=(\{2,3,4\},\{1,2\},\{3,4\})\in\GMLQ(s_1s_2\cdot\alpha,4)$.
    \begin{center}
    \resizebox{0.95\linewidth}{!}{
    \begin{tikzpicture}[scale=0.7]
    \def \w{1};
    \def \h{1};
    \def \r{0.3};
    \def \off{0.09};
    \begin{scope}[xshift=0cm]
    \foreach \i in {0,...,3}
    {
    \draw[gray!50] (0,\i*\h)--(\w*4,\i*\h);
    }
    \foreach \i in {0,...,4}
    {
    \draw[gray!50] (\w*\i,0)--(\w*\i,3*\h);
    }
    \foreach \xx\yy\i in {1/0,2/0,0/1,3/1,1/2,2/2,3/2}
    {
    \draw[blue] (\w*.5+\w*\xx,\h*.5+\h*\yy) circle (\r cm);
    }
    \foreach \xx\yy\i in {0/0,3/0,1/1,2/1,0/2}
    {
    \draw[red] (\w*.5+\w*\xx-\w*\r,\h*.5+\h*\yy+\h*\r) -- (\w*.5+\w*\xx+\w*\r,\h*.5+\h*\yy+\h*\r) -- (\w*.5+\w*\xx+\w*\r,\h*.5+\h*\yy-\h*\r)--(\w*.5+\w*\xx-\w*\r,\h*.5+\h*\yy-\h*\r)--(\w*.5+\w*\xx-\w*\r,\h*.5+\h*\yy+\h*\r);
    }

    \draw[black,<->] (5,1.5) -- (7,1.5);
    \node[black] at (6,2) {$\sigma_2$};
    \end{scope}

    \begin{scope}[xshift=8cm]
    \foreach \i in {0,...,3}
    {
    \draw[gray!50] (0,\i*\h)--(\w*4,\i*\h);
    }
    \foreach \i in {0,...,4}
    {
    \draw[gray!50] (\w*\i,0)--(\w*\i,3*\h);
    }
    \foreach \xx\yy\i in {1/0,2/0,0/1,1/1,3/1,2/2,3/2}
    {
    \draw[blue] (\w*.5+\w*\xx,\h*.5+\h*\yy) circle (\r cm);
    }
    \foreach \xx\yy\i in {0/0,3/0,2/1,1/2,0/2}
    {
    \draw[red] (\w*.5+\w*\xx-\w*\r,\h*.5+\h*\yy+\h*\r) -- (\w*.5+\w*\xx+\w*\r,\h*.5+\h*\yy+\h*\r) -- (\w*.5+\w*\xx+\w*\r,\h*.5+\h*\yy-\h*\r)--(\w*.5+\w*\xx-\w*\r,\h*.5+\h*\yy-\h*\r)--(\w*.5+\w*\xx-\w*\r,\h*.5+\h*\yy+\h*\r);
    }
    \draw[black,<->] (5,1.5) -- (7,1.5);
    \node[black] at (6,2) {$\sigma_1$};
    
    \end{scope}

    \begin{scope}[xshift= 16cm]
    \foreach \i in {0,...,3}
    {
    \draw[gray!50] (0,\i*\h)--(\w*4,\i*\h);
    }
    \foreach \i in {0,...,4}
    {
    \draw[gray!50] (\w*\i,0)--(\w*\i,3*\h);
    }
    \foreach \xx\yy\i in {1/0,2/0,3/0,0/1,1/1,2/2,3/2}
    {
    \draw[blue] (\w*.5+\w*\xx,\h*.5+\h*\yy) circle (\r cm);
    }
    \foreach \xx\yy\i in {0/0,2/1,3/1,1/2,0/2}
    {
    \draw[red] (\w*.5+\w*\xx-\w*\r,\h*.5+\h*\yy+\h*\r) -- (\w*.5+\w*\xx+\w*\r,\h*.5+\h*\yy+\h*\r) -- (\w*.5+\w*\xx+\w*\r,\h*.5+\h*\yy-\h*\r)--(\w*.5+\w*\xx-\w*\r,\h*.5+\h*\yy-\h*\r)--(\w*.5+\w*\xx-\w*\r,\h*.5+\h*\yy+\h*\r);
    }

    \end{scope}
    
    \end{tikzpicture}
}
\end{center}
\end{example}

In fact, if one views the multiline queue as a tensor product of Kirillov--Reshetikhin (KR) crystals, $\sigma_i$ precisely corresponds to the Nakayashiki-Yamada (NY) rule \cite{NY95}, describing the action of the combinatorial $R$ matrix on these crystals.  With the perspective of $\sigma_i$ as a combinatorial R matrix, one immediately obtains the following properties (see also \cite[Proposition 6.3]{AGS20} or \cite[Lemma 2.3]{vanleeuwen2006double} for a different approach in which $\sigma_i$ is built from co-plactic operators defined in \cite[Section 5.5]{Lot02} together with a cyclic shift operator).
\begin{lemma}\label{lem:coxeter}
        The $\sigma_i$'s satisfy the relations: 
        \begin{enumerate}[label=(\roman*)]
            \item $\sigma_i^2=id$, 
            \item $\sigma_i\sigma_j=\sigma_j\sigma_i$ if $|i-j|\geq 2$, 
            \item $\sigma_i\sigma_{i+1}\sigma_i=\sigma_{i+1}\sigma_i\sigma_{i+1}$.
        \end{enumerate}
\end{lemma}

In \cite{AGS20}, it is shown that  $\sigma_i$ preserves the generalized FM projection map $\mathfrak{p}$: 
\[
\mathfrak{p}(M)=\mathfrak{p}(\sigma_i(M)),
\]
thus obtaining a version of \eqref{eq:FM result} for twisted multiline queues. In particular, this implies that the stationary distribution of the TASEP can be computed as the cardinality of the sets of twisted multiline queues with fixed bottom-row labelings. 
More generally, they showed the following, with an example given in \cref{ex:sigma theorem}.

\begin{theorem}[{\cite[Theorem 3.1]{AGS20}}]\label{thm:AGS}
Let $B\in\GMLQ(\alpha,n)$, and let $1\leq i\leq L$ where $L=\ell(\alpha)$. Then the labeled arrays $L_G(B)$ and $L_G(\sigma_i(B))$ coincide on all rows except for row $i+1$. 
\end{theorem}

In \cref{sec:GMLQ maj}, we will strengthen the above result by defining a major index statistic $\maj_G$ on twisted multiline queues, which is also preserved by the involution $\sigma_i$.

\begin{example}\label{ex:sigma theorem}
    For $\alpha=(4,3,2,5,1)$ and $n=6$, we show $B\in\GMLQ(\alpha,n)$ on the left and $\sigma_2(B)\in\GMLQ(s_2\cdot\alpha,n)$ on the right along with their labelings $L(B)$ and $L(\sigma_2(B))$. Notice that the labels of the sites coincide on all rows except for row 3. 
    
    \begin{center}
        \resizebox{!}{1.6in}{\begin{tikzpicture}
            \node[] at (0,0) {\includegraphics[scale = 0.4]{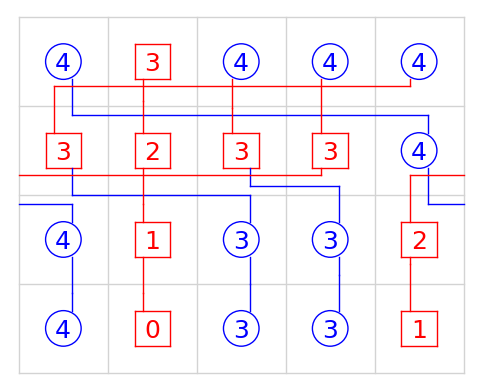}};

            \node[] at (7,0) {\includegraphics[scale = 0.4]{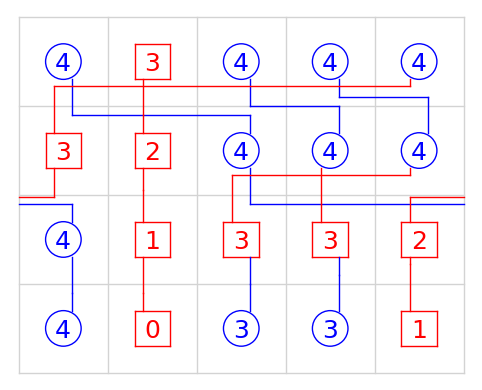}};

            \draw[<->] (2.75,0.0) -- (4,0.0) {};
            \node[black] at (3.4,0.25) {$\sigma_2$};

            \node[gray] at (-2.5,-1.35) {\tiny 1};
            \node[gray] at (-2.5,-0.45) {\tiny 2};
            \node[gray] at (-2.5,0.45) {\tiny 3};
            \node[gray] at (-2.5,1.35) {\tiny 4};

            \node[gray] at (4.5,-1.35) {\tiny 1};
            \node[gray] at (4.5,-0.45) {\tiny 2};
            \node[gray] at (4.5,0.45) {\tiny 3};
            \node[gray] at (4.5,1.35) {\tiny 4};
            
        \end{tikzpicture}
        }
    \end{center}
\end{example}

\subsection{Bosonic multiline queues}\label{sec:bmlqs}
By allowing each cell of the multiline queue to contain any number of particles (rather than only 0 or 1), we obtain \emph{bosonic multiline queues}, which first appeared in the context of a family of statistical mechanics models called the \emph{totally asymmetric zero range process} (see, for instance, \cite{KMO-multiline-2016,AMM22}). These objects are in bijection with $\inv$-free Haglund--Haiman--Loehr tableaux \cite{HHL08} and in (weight preserving) bijection with $\quinv$-free tableaux \cite{AMM20}, yielding combinatorial formulas for the modified Hall--Littlewood polynomials. 

 From now on, we may refer to classical multiline queues as \emph{fermionic multiline queues} when the context is not clear.

\begin{definition}
    Let $\lambda$ be a partition and $n$ a positive integer. A \emph{bosonic multiline queue} of type $(\lambda,n)$ is a configuration of particles on a $\lambda_1\times n$ grid, such that each site can contain any number of particles, and each row $j$ contains a total of $\lambda_j'$ particles. Like with multiline queues, rows are labeled from bottom to top. We can represent a bosonic multiline queue by the tuple $D=(D_1,\ldots,D_{\lambda_1})$, where each $D_i$ is a \emph{multiset} of $[n]$ of size $\lambda_i'$. Denote the set of bosonic multiline queues of type $(\lambda,n)$ by $\MLD(\lambda,n)$.
\end{definition}

\begin{defn}\label{def:cw_bmlqs}
Define the bosonic multiline queue column reading word
\[
\widetilde{\cw}(D):=\rev(\cw(D)),
\]
where $\cw(D)$ is given by the standard multiline queue reading order. See \cref{ex:mld cw}.
\end{defn}

One can alternatively obtain $\widetilde{\cw}(D)$ from the integer matrix representation of a bosonic multiline queue as follows. Scanning the columns from left to right and the cells from top to bottom within each column, construct the reading word $\widetilde{\cw}(D)$ by recording $k$ copies of ``$r$'' for a cell $(r,j)$ (in row $r$, column $j$) that contains a nonzero entry $k$. 

\begin{example}\label{ex:mld cw}

    Consider the bosonic multiline queue $D=(\{1,2,3,3,4,4\},\{2,3,3\},\{1,4\})$. Below we show two equivalent representations of $D$: as a ball diagram on the left, and as a positive-integer matrix on the right. In the matrix, the number of balls in the given cell is represented by an integer, with a vacancy corresponding to 0.

    \begin{center}
        \includegraphics[scale = 0.25]{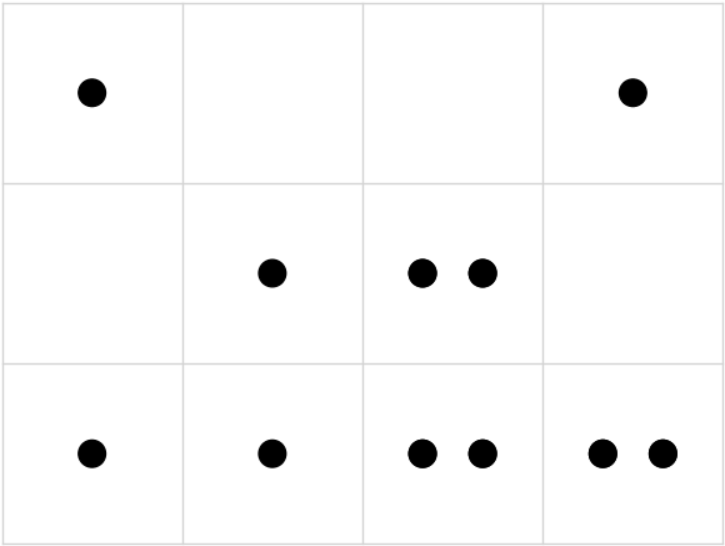}
        \hspace{1cm}
        \includegraphics[scale=0.25]{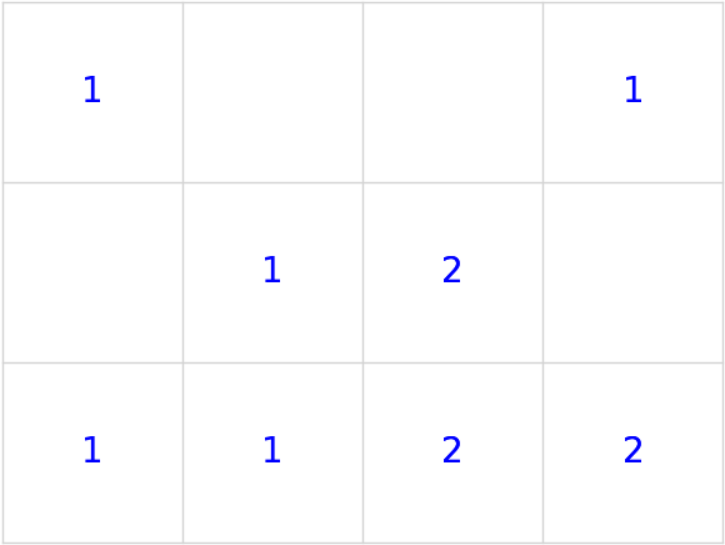}
    \end{center}
    Using the multiline queue reading order we obtain $\cw(D) = 3\,1\,|\,2\,1\,| 2\, 2\, 1\, 1\, | \, 3\,1\,1$, and thus $\widetilde\cw(D) = 1\,1\,3\,|\,1\,1\,2\,2\,|\,1\,2\,|\,1\,3$.

\end{example}

\begin{defn}\label{def:matched above below MLD word}
Let $D=(D_1,\ldots,D_L)$ be a bosonic multiline queue. We call a ball in $D_i$ \emph{matched below} if it is paired in $\Par_i(\widetilde{\cw}(D))$ for any $1\leq i<L$, and we call it \emph{matched above} if it is paired in $\Par_{i-1}(\widetilde{\cw}(D))$ for any $1<i\leq L$. Otherwise, we call it \emph{unmatched below} (resp.~\emph{above}).
\end{defn}

There is a pairing process on bosonic multiline queues that is analogous to the FM algorithm that produces a labeling of its particles and a major index statistic. To make the connection with the modified Hall--Littlewood polynomials, the particles are paired \emph{strictly to the left} rather than weakly to the right as for multiline queues.\footnote{Sending $t\mapsto t^{-1}$ and then specializing $t^{-1}=0$ in the plethystic definition of $\widetilde H_\lambda(X;q,t)$ recovers the inverted pairing rule for the corresponding multiline queues (see~\cite[Section 8]{AMM22} for details).} Aside from this rule change, the pairing process is identical to the FM algorithm given in \cref{def:FM}, which has an equivalent formulation in terms of cylindrical matching operators with respect to the modified column reading word $\widetilde{\cw}$.

\begin{defn}[Bosonic multiline queue labeling process]\label{def:mld labeling} Let $D=(D_1,\ldots,D_L)$ be a bosonic multiline queue of shape $(\lambda,n)$. Define the \emph{labeled multiline queue} $L(D)$ by replicating $D$ and sequentially labeling the balls, as follows. For each row $r$ for $r=L,\ldots,2$, each unlabeled ball in $B_r$ is labeled $r$. Next, for $\ell=L,\ldots,r$, let $\widetilde{\cw}(D)^{(\ell,r)}$ be the restriction of $\widetilde{\cw}(D)$ to the balls labeled $\ell$ in $D_r$ and the unlabeled balls in $D_{r-1}$. The balls in row $r-1$ that are cylindrically matched in $\Par^c_{r-1}(\widetilde{\cw}(D)^{(\ell,r)})$ acquire the label $\ell$. 
To complete the process, all unpaired balls in row 1 are labeled "$1$". Such a labeling is shown in \cref{ex:labeled bMLQ}.
\end{defn}

The labeling process above can be reformulated analogously to that of multiline queues, to produce a labeled diagram with balls paired between adjacent rows. 

\begin{defn}[Bosonic multiline queue pairing process]\label{def:MLDpairing}
Let $D=(D_1,\ldots,D_L)$ be a bosonic multiline queue of type $(\lambda,n)$. For each row $r$ for $r=L,\ldots,2$:
\begin{itemize}
\item Every unlabeled ball is labeled "$r$".
\item Once all the balls in row $r$ are labeled, each of them is sequentially paired to the first unlabeled ball \emph{strictly to its left} in row $r-1$, wrapping around from column $1$ to column $n$ if necessary. The order in which balls are paired from row $r$ to row $r-1$ is from the largest label $L$ to the smallest label $r$, and (by convention) from right to left among balls with the same label.
\end{itemize}
To complete the process, all unpaired balls in row $1$ are labeled "$1$". We associate to $D$ the multiset $\Pair(D)=\{(r(p),\ell(p),\delta(p)):\ p\ \text{is a pairing in}\ D\}$ where for each pairing $p$, $r(p)$ is the row from which the pairing originates, $\ell(p)$ is the label corresponding to the pairing, and $\delta(p)$ is equal to $1$ if the pairing wraps and $0$ otherwise.
\end{defn} 

Recall that a pairing from a particle at site $(r,a)$ to a particle at site $(r-1,b)$ is wrapping in a multiline queue if $a>b$, whereas in a bosonic multiline queue the pairing is wrapping if $a\leq b$. This switch is reflected in the definition of $\widetilde{\cw}$ (\cref{def:cw_bmlqs}) as the reverse of $\cw$. 

\begin{center}
\resizebox{0.7\linewidth}{!}{\begin{tikzpicture}[scale=0.7]
    \def \w{1};
    \def \h{1};
    \def \r{0.25};
    \def \a{6};
    \def \b{-1};
    
    \foreach \i in {0,1,2}
    {
    \draw[gray!50] (0+\b,\i*\h)--(\w*4+\b,\i*\h);
    \draw[gray!50] (\a+0,\i*\h)--(\a+\w*4,\i*\h);
    \draw[gray!50] (2*\a+0,\i*\h)--(2*\a+\w*4,\i*\h);
    }
    \foreach \i in {0,...,4}
    {
    \draw[gray!50] (\w*\i+\b,0)--(\w*\i+\b,2*\h);
    \draw[gray!50] (\a+\w*\i,0)--(\a+\w*\i,2*\h);
    \draw[gray!50] (2*\a+\w*\i,0)--(2*\a+\w*\i,2*\h);
    }
    
    \foreach \xx\yy\i in {3+\b/1/$a$,1+\b/0/$b$,1+\a/1/$a$,3+\a/0/$b$,1+2*\a/1/$a$,1+2*\a/0/$\,$}
    {
    \fill (\w*.5+\w*\xx,\h*.5+\h*\yy) circle (\r cm);
    \node at (\w*.5+\w*\xx,\h*2.3) {\i};
    }

    \node at (2+\b,-.5) {multiline queue};
    \node at (1.5*\a+2,-.5) {bosonic multiline queue};
    \node at (1.5*\a+2,1) {or};
    \draw[blue,-stealth] (\w*3.5+\b,\h*1.5-\r)--(\w*3.5+\b,\h*.95)--(\w*4.3+\b,\h*.95);%
    \draw[blue] (-.2+\b,\h*.95)--(\w*1.5+\b,\h*.95)--(\w*1.5+\b,\h*.5+\r);
   \draw[blue,-stealth] (\w*1.5+\a,\h*1.5-\r)--(\w*1.5+\a,\h*.95)--(\w*-.3+\a,\h*.95);%
    \draw[blue] (4.3+\a,\h*.95)--(\w*3.5+\a,\h*.95)--(\w*3.5+\a,\h*.5+\r);%
    \draw[blue,-stealth] (\w*1.5+2*\a,\h*1.5-\r)--(\w*1.5+2*\a,\h*1.1)--(\w*-.3+2*\a,\h*1.1);%
    \draw[blue] (4.3+2*\a,\h*.9)--(\w*1.5+2*\a,\h*.9)--(\w*1.5+2*\a,\h*.5+\r);%
\end{tikzpicture}
}
\end{center}

\begin{example}\label{ex:labeled bMLQ}
    Consider the bosonic multiline queue $D$ from \cref{ex:mld cw}. In \cref{fig:mld-example} we show the corresponding pairings (left) and labels on the individual balls in $L(D)$ (right).  

    \begin{figure}[b]
        \centering
        \resizebox{1.5in}{!}{\includegraphics[height = 4.5cm]{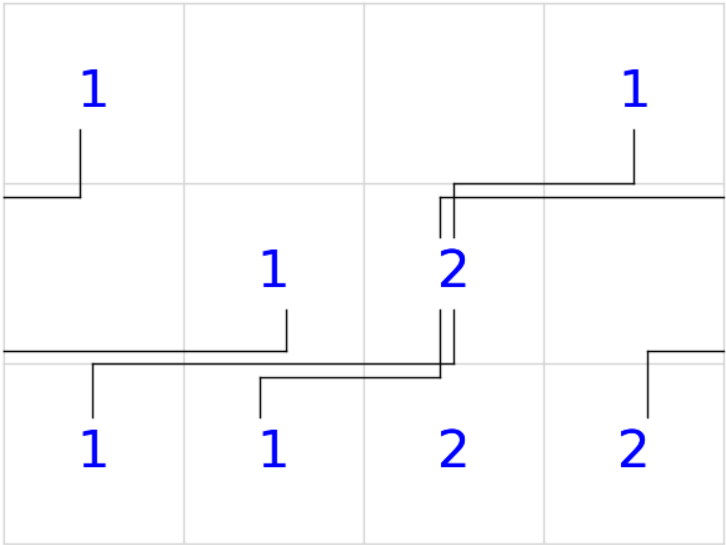}  }    
        \hspace{1cm}
         \resizebox{1.5in}{!}{\includegraphics[height = 4.5cm]{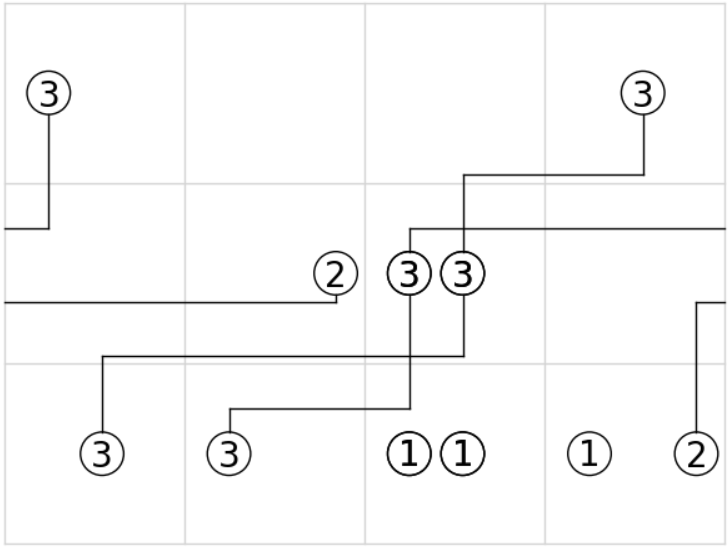}}
        \caption{Representations of pairings of the bosonic multiline queue $D=(\{1,2,3,3,4,4\},\{2,3,3\},\{1,4\})$.}
        \label{fig:mld-example}
    \end{figure}
\end{example}

As with multiline queues, the content of a bosonic multiline queue $D = (D_1,\ldots,D_L)$ is the monomial recording the numbers of balls in each column $$x^D = \prod_{j=1}^{L}\prod_{b\in D_j} x_b.$$ The bosonic major index, on the other hand, is determined by the \emph{nonwrapping pairings}. 

\begin{defn}[Major index for bosonic multiline queues]\label{def:maj_mlds}
The \emph{major index} of a bosonic multiline queue $D$, denoted by $\widetilde{\maj}(D)$, computed from $\Pair(D)$ as 
\[
\widetilde{\maj}(D)=\sum_{(r,\ell,\delta)\in\Pair(D)}(1-\delta)(r-\ell+1).
\]
\end{defn} 

\begin{defn}\label{def:MLD-weight}
    The \emph{weight} of a bosonic multiline queue is defined as $\wt(D)=x^D q^{\widetilde{\maj}(D)}$. 
\end{defn}

\begin{example}
    Let $D=(\{1,2,3,3,4,4\},\{2,3,3\},\{1,4\})$ be the bosonic multiline queue from \cref{fig:mld-example}. $\widetilde\maj(D)=5$, coming from the nonwrapping pairings between the cells $(3,4)$ and $(2,3)$, the cells $(2,3)$ and $(1,2)$, and the cells $(2,3)$ and $(1,1)$. 
\end{example}

\begin{remark}
In the bosonic analogue of \cref{sec:tableaux}, one obtains bosonic multiline queues by reinterpreting the tableau formula for modified Hall--Littlewood polynomials via $\quinv$-free fillings of a Young diagram. Just as in the fermionic case, for any fixed row content there is a unique $\quinv$-free filling, so we can identify its row content with the corresponding bosonic multiline queue. Under this identification, the statistic $\maj$ on fillings agrees exactly with the $\widetilde\maj$ defined on bosonic multiline queues. 
\end{remark}
   
The fact that $\widetilde\maj(D)=\cocharge(\widetilde{\cw}(D))$ follows from the same argument as the proof of \cref{lem:charge cw} based on the bosonic multiline queue labeling procedure in \cref{def:mld labeling}.
   
Thus we have the following formula for the modified Hall--Littlewood polynomials. 
\begin{theorem}\label{thm:H MLD}
\begin{equation}\label{modifiedHL}
\widetilde{H}_{\lambda}(x_1,\ldots,x_n;q)=\sum_{D\in\MLD(\lambda,n)} x^Dq^{\widetilde\maj(D)}=\sum_{D\in\MLD(\lambda,n)} x^D q^{\cocharge(\widetilde{\cw}(D))}.
\end{equation}
\end{theorem}

\begin{example}\label{ex:cw D}
    For the bosonic multiline queue $D=(\{1,2,3,3,4,4\},\{2,3,3\},\{1,4\})$ in $\MLD(\lambda,4)$ with $\lambda=(3,3,2,1,1,1)$ in \cref{fig:mld-example}, the unique $\quinv$-free filling of $\dg(\lambda)$ with the same row content as $D$ is
    \[\sigma=\raisebox{9pt}{\tableau{4&1\\3&3&2\\2&1&4&4&3&3}}
    \]
    We have that $\widetilde{\cw}(D)=\rev(31.21.2211.311)=1\,1\,3\,|\,1\,1\,2\,2\,|\,1\,2\,|\,1\,3$ so \[\cocharge(\widetilde{\cw}(D)) = \cocharge(321)+\cocharge(213)+\cocharge(12)+3\cocharge(1)=5=\maj(\sigma).\] 
\end{example}


\begin{definition}
    A bosonic multiline queue $D$ is called \emph{nonwrapping} if for each pairing from a particle at site $(r,a)$ to a particle at site $(r-1,b)$, it satisfies $a>b$. Denote the set of nonwrapping bosonic multiline queues by $\MLD_0(\lambda,n)$.
\end{definition}

Note that the technical condition for a bosonic multiline queue $D=(D_1,\ldots,D_L)$ to be nonwrapping is that for each $2\leq i\leq L$, and for each $k\in[n]$, the number of elements smaller than $k$ in $D_{i-1}$ is not less than the number of elements smaller than or equal to $k$ in $D_i$. An immediate consequence is that the set $\MLD_0(\lambda,n)$ is in direct correspondence with the set of semistandard tableaux $\SSYT(\lambda',n)$ via a natural map that sends the row contents of a bosonic multiline queue to the row contents of a semistandard tableau.

\begin{lemma}\label{lem:schur-mld}
\[
s_{\lambda}(x_1,\ldots,x_n)=\sum_{D\in\MLD_0(\lambda',n)}x^D
\]
\end{lemma}

\begin{proof}
Let $D=(D_1,\ldots,D_L)\in\MLD_0(\lambda',n)$, where $L=\lambda'_1$. Let $F$ be the filling of the Young diagram of shape $\lambda$ such that row $i$ contains the entries in $D_i$ in weakly increasing order. Then $F\in\SSYT(\lambda)$. Conversely, for any filling $F\in \SSYT(\lambda,n)$, the bosonic multiline queue $D=(D_1,\ldots,D_L)$ obtained by letting $D_i$ be the row content of row $i$ of $F$, is nonwrapping. By construction, this map preserves the content, and so $x^D=x^F$, and the lemma follows. 
\end{proof}

A nonwrapping bosonic multiline queue has maximal $\maj$, equal to $n(\lambda')$. This is equivalent to specializing $\widetilde{H}_{\lambda}(X;q,0)$ at $q=\infty$ and normalizing, to obtain $s_{\lambda'}(x)$ and recover \cref{lem:schur-mld}.

\section{Collapsing procedure on multiline queues}\label{sec:collapse}
This section introduces a \emph{collapsing procedure} for binary matrices using crystal operators, which we will use as a tool to study the $\maj$ statistic on multiline queues. Considering the matrices as twisted multiline queues, this procedure produces a bijection from the latter to pairs consisting of a nonwrapping multiline queue and a semistandard Young tableau, such that the $\maj$ of the multiline queue is transferred to the $\charge$ statistic of the semistandard Young tableau. In \cref{sec:mlq/d_RSK}, we further use bi-directional collapsing to define a form of Robinson--Schensted--Knuth (RSK) correspondence for multiline queues (where the recording object is a nonwrapping multiline queue). 

\subsection{Collapsing on (twisted) multiline queues via row operators}\label{sec:raising-lowering}

Let $\M$ be the set of binary matrices with finite support, and let $\M(L,n)$ be the set of such matrices on $L$ rows and $n$ columns. We can consider a matrix $B\in\M(L,n)$ as a twisted multiline queue (see \cref{sec:Gmlqs}) where balls and vacancies are the 1's and 0's respectively.

Throughout this section, unless explicitly specified, $B\in\M(L,n)$ is a binary matrix given by $B = (B_1,B_2,\ldots,B_L)$, where $B_j \subseteq [n]$ is the set of column labels of the balls (1's) of row $j$ of $B$ for $1\leq j \leq L$. 

\begin{definition}\label{def:dropping}
The \emph{dropping operator} $e_i$ acts on $M$ by moving the smallest entry that is unmatched above in $\Par(B_{i+1},B_i)\coloneqq \Par_i(\cw(M))$ from $B_{i+1}$ to $B_{i}$. In $B$, this corresponds to the \emph{leftmost} ball that is unmatched above in row $i+1$ dropping to row $i$. Additionally, define $e_i^{\star}$ as the operator that moves entries unmatched above in $\Par(B_{i+1},B_i)$ from row $i+1$ to row $i$. Then $e_i^{\star}(M)=e_i^{\phi_i(M)}(M)$, where $\phi_i(M)$ is the total number of entries unmatched above in $\Par(B_{i+1},B_i)$.
\end{definition}

\begin{definition}\label{def:lifting}
The \emph{lifting operator} $f_i$ acts on $M$ by moving the largest entry that is unmatched below in $\Par(B_{i+1},B_i)=\Par_i(\cw(M))$ from $B_{i}$ to $B_{i+1}$. In $B$, this corresponds to the \emph{rightmost} ball that is unmatched below in row $i$ being lifted to row $i+1$.
\end{definition}

In fact, the dropping and lifting operators on twisted multiline queues are simply the classical crystal operators on its column reading word. In particular, the collapsing operators $e_i^\star$ satisfy the defining relations of the $0$‐Hecke algebra of type A (equivalently the Demazure operators at $q=0$) \cite{Norton79}, so that 
$e_i^\star$ gives a well-defined 0-Hecke action on multiline queues.

\begin{lemma}\label{prop:rais/low-and-drop/lift}
The dropping and lifting operators satisfy
\begin{itemize}
    \item $\cw(e_i(B)) = E_i(\cw(B)) $
    \item $\cw(f_i(B)) = F_i(\cw(B)) $
\end{itemize}
where $E_i$ and $F_i$ are the standard raising and lowering \emph{crystal operators in type A on words}.
\end{lemma}
See \cref{sec:preliminaries} and \cite{BumpSchilling17} for further information about the crystal operators in the context of crystal bases. Along with the fact that the classical operators $E_i$ and $F_i$ are inverses when they act non-trivially, the previous lemma implies the following.
\begin{lemma}\label{lem:inverses}
For $1\leq i\leq n-1$, when $e_i(B)\neq B$, then $f_i(e_i(B))=B$. When $f_i(B)\neq B$, then $e_i(f_i(B))=B$. 
\end{lemma}

We borrow the terminology from tableaux insertion for the following definition. 

\begin{defn}
   Let $N=(N_1,\ldots,N_L)$ be a nonwrapping multiline queue. We say that $x$ \emph{bumps} $a\in B_r$ at row $r$ when $a$ is unmatched above in $\Par(N_r\cup\{x\},N_{r-1})$. In this case, necessarily $a<x$.
\end{defn}

Throughout this article, our convention is that sequences of operators act from right to left. For ease of reading, we will use multiplicative notation on the operators $e_i^\star$ and omit the composition symbol $\circ$. 

The following result can be found, for instance, in \cite{vanleeuwen2006double,Lot02}. 

\begin{theorem}\label{thm:dropping-operators-algebra}
    The collapsing operators $e_i^\star$ satisfy the following algebraic relations:
    \begin{itemize}
        \item[(i)] $(e_i^\star)^2 = e_i^\star$,
        \item[(ii)] $e_i^\star e_j^\star = e_j^\star e_i^\star$ whenever $|i-j|\geq 2$,
        \item[(iii)] $e_i^\star e_{i+1}^\star e_i^\star = e_{i+1}^\star e_i^\star e_{i+1}^\star.$
    \end{itemize}
\end{theorem}

Define $e^{\star}_{[a,b]}\coloneqq e_{a}^\star e_{a+1}^\star\cdots e_{b}^\star$ for $1\leq a\leq b$ as a sequential application of operators (read from right to left). As an operator on twisted multiline queues, $e_{[a,b]}^\star$ sequentially drops all unmatched above balls from row $b+1$ to $b$, then from $b$ to $b-1$, and so on, down to row $a$.  

\begin{defn}[Collapsing]\label{def:collapsing_mlqs}
Let $L$ and $n$ be positive integers. Define \emph{collapsing} on binary matrices as the map 
\begin{align}
    \rho \; \colon\; \M(L,n) & \;\to\; \bigcup_\mu \MLQ_0(\mu,n)\times \SSYT(\mu',L)\\
    B &\;\longmapsto\; (\rho_N(B),\rho_Q(B)),
\end{align}
where $\rho_N(B)$ is given by
\begin{equation}\label{eq:dropping} 
\rho_N(B)\coloneqq e_{[1,L-1]}^\star e_{[1,L-2]}^\star\cdots e_{[1,2]}^\star e_{[1,1]}^\star(B),
\end{equation}
and $\rho_Q(B)$ is the semistandard tableau with content $(|B_1|,\ldots,|B_L|)$, whose entries $i$ record the difference in row content between $e^\star_{[1,i-1]}\cdots e^\star_{[1,1]}(B)$ and $e^\star_{[1,i-2]}\cdots e^\star_{[1,1]}(B)$ for $1\leq i\leq L$.
\end{defn}

\begin{remark}\label{rem:adapted-strings}
    The data of the number of times each dropping operator $e_i$ acts in a given multiline queue $B$ is the \emph{adapted string} of the element in the sense of Littelmann \cite{Littelmann1998}. Indeed, the \emph{nice decomposition} introduced therein for the longest permutation
    $$w_0 = s_1(s_2s_1)(s_3s_2s_1)\ldots(s_{L-1}s_{L-2}\ldots s_1).$$
    reflects the order of application of the dropping operators in \eqref{eq:dropping}. 
\end{remark}

\begin{remark}\label{rem:collapse in diagram}
It is convenient to visualize collapsing as a procedure occurring directly on the diagram of a binary matrix, in which sequentially, row by row from bottom to top, balls that are unmatched above are dropped to the row below, until a nonwrapping multiline queue is reached.

Let $B=(B_1,\ldots,B_L)\in\M(n,L)$ be a binary matrix. We build the nonwrapping multiline queue $\rho_N(B)$ and the recording tableau $\rho_Q(B)$ recursively, as follows. Initiate $N_1$ to be a copy of row 1 of $B$, and set $Q_1$ to be a single row with $|B_1|$ boxes containing the entry 1. Sequentially, for each row $r=2,3,\ldots,L-1$, let $N_r$ be the nonwrapping multiline queue obtained from collapsing the bottom $r$ rows of $B$. Let $Q_r$ be the partially built recording tableau whose shape is the conjugate of the shape of $N_r$, and whose content is $1^{|B_1|}2^{|B_2|}\ldots r^{|B_r|}$. Place the balls from row $r+1$ of $B$ in row $r+1$ of $N_r$ to obtain $N'_{r+1}$. Set $u=r+1$.

\begin{itemize}
        \item If there are no balls unmatched above in row $u$ of $N'_{r+1}$, the collapsing for $r+1$ is complete. 
        \item Otherwise, drop all balls that are unmatched above in row $u$ of $N'_{r+1}$ to row $u-1$, update $N'_{r+1}$, and repeat the step for $u=u-1$. 
\end{itemize}
Once the collapsing for row $r+1$ is complete, set $N_{r+1}$ to be the fully collapsed $N'_{r+1}$ (which is by construction nonwrapping). For each row $\ell=1,\ldots,r+1$, take the difference between the number of entries in row $\ell$ of $N_{r+1}$ and row $\ell$ of $N_r$, and add that many boxes filled with the entry ``$r+1$'' to row $\ell$ of the recording tableau $Q_r$ to obtain $Q_{r+1}$. It is then immediate that the shape of $Q_{r+1}$ is conjugate to the shape of $N_{r+1}$, and its content matches the row sizes of the first $r+1$ rows of $B$. Once row $L$ of $B$ is collapsed, we set $\rho_N(B)=N_L$ and $\rho_Q(B)=Q_L$. We illustrate this procedure in \cref{ex:collapse}.
\end{remark}

Since the RSK correspondence commutes with the Kashiwara crystal operators (see, e.g.~\cite[Corollary 9.2]{BumpSchilling17}), one can deduce that \eqref{eq:dropping} is a bijection. 
We provide a self-contained elementary proof. 

\begin{lemma}\label{prop:rho}
Let $B=(B_1,\ldots,B_L)\in\M(L,n)$ be a binary matrix. Then $\rho_N(B)\in\MLQ_0(\mu,n)$ for some partition $\mu$, and $\rho_Q(B)$ is a semistandard Young tableau of shape $\mu'$ with content $(|B_1|,\ldots,|B_L|)$.
\end{lemma}

\begin{proof}
    The proof of both statements is by induction on the number of rows of $B$. Let $2\leq i\leq L$ and denote the collapsing operator on the bottom $k$ rows by $$\rho_N^{(k)}=e_{[1,k-1]}^\star\cdots e_{[1,2]}^\star e_{[1,1]}^\star$$ so that $\rho_N(B)=\rho_N^{(L)}(B)$. 
    We will show that the bottom $k$ rows of $\rho_N^{(k)}(B)$ form a nonwrapping multiline queue and the partial recording object $Q^{(k)}$ is a semistandard tableau of shape $\mu^{(k)}$ where, for $1\leq j \leq i$, $\mu^{(k)}_j$ is the number of balls in row $j$ of $\rho_N^{(k)}(B)$, and with content $(|B_1|,|B_2|,\ldots,|B_k|).$ Note that to show that the bottom $j$ rows of a multiline queue are nonwrapping, it is enough to show that the operators $e_1^\star,\ldots,e_{j-1}^\star$ act trivially.
    
    For the base case, $\rho_N^{(1)}(B)=B$, and by definition the first row of $B$ forms a nonwrapping multiline queue and the partial recording object $Q^{(1)}$ is a tableau with one row of size $|B_1|$, filled with $1$'s. Now assume the statement holds for some $r\geq 1$: the bottom $r$ rows of $\rho_N^{(r)}(B)$ are nonwrapping. We claim that $e_1^\star,\ldots,e_r^\star$ act trivially on $\rho_N^{(r+1)}(B)=e_{[1,r]}^\star\circ \rho_N^{(r)}(B)$. Indeed, $e_1^\star\circ\rho_N^{(r+1)}(B)=\rho_N^{(r+1)}(B)$, and the relations in \cref{thm:dropping-operators-algebra} imply that $e_{j}^\star e_{[1,r]}^\star= e_{[1,r]}^\star e_{j-1}^\star$ for any $1<j\leq r$; and so for $1< j\leq r$, we get
    \begin{align*}
e_j^\star\circ \rho_N^{(r+1)}(B)=e_j^\star e_{[1,r]}^\star \circ \rho_N^{(r)}(B)=e_{[1,r]}^\star e_{j-1}^\star\circ \rho_N^{(r)}(B)=e_{[1,r]}^\star\circ \rho_N^{(r)}(B)
=\rho_N^{(r+1)}(B)
    \end{align*}
where the third equality is due to the fact that $e_{j-1}^\star$ acts trivially on $\rho_N^{(r)}(B)$ for $1<j\leq r$.

    To construct $Q^{(r+1)}$ from $Q^{(r)}$, the new balls appearing in each row of $\rho_N^{(r+1)}(B)$ after applying each $e_{[1,i]}^\star$ to $\rho_N^{(r)}(B)$ are recorded in the corresponding row of $Q^{(r)}$ as the entry ``$r+1$''. By the above, $\mu_{j}^{(r)}\leq \mu_{j}^{(r+1)}\leq \mu_{j-1}^{(r)}$ for each $2\leq j\leq r+1$. Thus row $j$ of $Q^{(r+1)}$ contains $\mu_{j}^{(r+1)}-\mu_{j}^{(r)}\geq 0$ entries ``$r+1$'' for $1\leq j\leq r+1$, and the shape of $Q^{(r+1)}/Q^{(r)}$ is a horizontal strip of size $|B_{r+1}|$. Therefore, the shape of $Q^{(r+1)}$ is the partition $\mu^{(r+1)}$ and its content is $(|B_1|,\ldots,|B_{r+1}|)$. \\
Thus we have $\rho_N(B)=\rho_N^{(L)}(B)\in\MLQ_0(\mu,n)$, $\rho_Q(B)=Q^{(L)}\in\SSYT(\mu',\alpha)$, where $\mu'=\mu^{(L)}$ and $\alpha=(|B_1|,\ldots,|B_L|)$.
\end{proof}

The following result is a powerful property of collapsing on multiline queues: the $\maj$ of a multiline queue $M$ becomes the $\charge$ of the recording tableau $\rho_Q(M)$.  This fact is well-known to experts (see, for instance, \cite[Section 4]{LenartSchilling11} and the references therein), following from the fact that energy is preserved under the action of non-zero crystal operators. However, we shall give a self-contained multiline queue proof in \cref{sec:mlq/d_RSK}.
 
\begin{theorem}\label{thm:charge}
    Let $M\in\MLQ(\lambda,n)$ be a multiline queue. 
    Then
    \[\maj(M)=\charge(\rho_Q(M)).\]
\end{theorem}

Since the $\charge$ of a tableau can be computed from sequential applications of the (classical and cylindrical) matching functions $\Par_i$ and $\Par_i^c$, the following lemma can be considered a refinement of the above result that holds for any binary matrix $B\in\M(L,n)$. This will serve as a useful tool to understand the structure of the recording tableau arising from collapsing. 
The proof of this lemma is also found in \cref{sec:mlq/d_RSK}. 

\begin{lemma}\label{lem:pairing-position-word-Q}
    Let $B\in\M(L,n)$ and $Q = \rho_Q(B)$. Then $\Pair_i(B_{i+1},B_{i}) = \Pair_i(\crw(Q)).$
\end{lemma}

\begin{example}\label{ex:collapse}
Applying $\rho$ to the multiline queue $(\{1,3,4\},\{1,4,5\},\{2,5\},\{1,3\},\{4\})$ below 
yields the following pair:

\begin{figure}[h!]
    \centering
    \def \scale{0.35}
    \begin{subfigure}[t]{0.25\textwidth}
    \centering
    \includegraphics[scale = \scale]{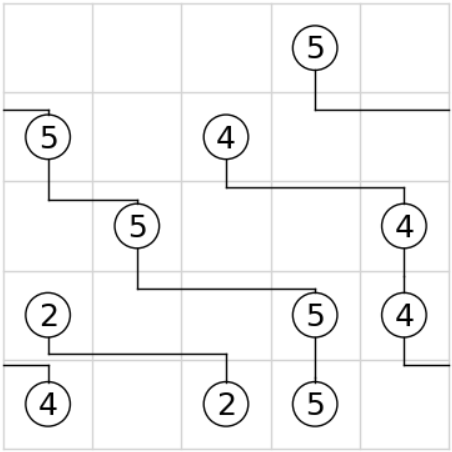}
    \caption*{{\Large $\substack{\lambda \;=\; (5,4, 2) \\ \\ \maj(M) = 4}$}}
    \end{subfigure}
    \raisebox{45pt}{\Large {$\longrightarrow\;\;$}}
    \begin{subfigure}[t]{0.035\textwidth}
        \begin{tikzpicture}
            \draw (7.5,0) arc (-30:30:3.5);
            \draw (0,0) arc (210:150:3.5);
        \end{tikzpicture}        
    \end{subfigure}
    \begin{subfigure}[t]{0.2\textwidth}
    \centering
    \raisebox{8pt}{\includegraphics[scale = 0.3]{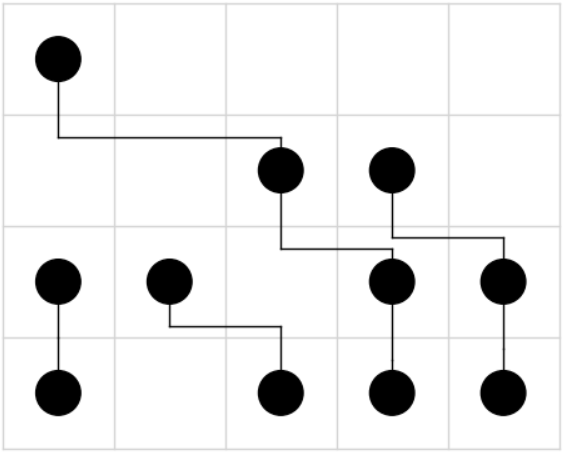}}
    \caption*{{\Large $\substack{\mu \,=\,  (4,3,2,2)\\ \\ \maj(\rho_N(M)) = 0}$}}
    \end{subfigure}
    $\quad$
    \raisebox{47pt}{\huge $,$}
    \begin{subfigure}[t]{0.18\textwidth}
    \centering
    \raisebox{65pt}{{$\tableau{4\\ 3&4\\ 2&2&3&5\\ 1&1&1&2}$}}
    \caption*{{\Large $\substack{\charge(\rho_Q(M)) = 4 \\ }$}}
    \end{subfigure}
\end{figure}

We show the step-by-step collapsing of the rows from bottom to top, following \cref{rem:collapse in diagram}. For each $r$, the unfilled balls represent those that have not yet been collapsed, $Q_r$ is the partial recording tableau, and $N_r$ is the nonwrapping multiline queue represented by the black balls. 

\medskip

\begin{center}
    \begin{tikzpicture}[scale=0.4]
        \def \w{1};
        \def \h{1};
        \def \r{0.25};
        \node at (-1,6) {$r$};
        \node at (-1,2.5) {$N_r$};
        \node at (-1,-2) {$Q_r$};
        \foreach \i in {0,...,5}
        {
        \node at (2.5+6*\i,6) {\i};
        }
        \begin{scope}[xshift=0cm]
        \foreach \i in {0,...,5}
        {
        \draw[gray!50] (0,\i*\h)--(\w*5,\i*\h);
        }
        \foreach \i in {0,...,5}
        {
        \draw[gray!50] (\w*\i,0)--(\w*\i,5*\h);
        }
        \foreach \xx\yy in {0/0,2/0,3/0,0/1,3/1,4/1,1/2,4/2,0/3,2/3,3/4}
        {
        \draw (\w*.5+\w*\xx,\h*.5+\h*\yy) circle (\r cm);
        }
        \node at (2.5,-2) {\scalebox{0.8}{\tableau{\emptyset\\}}};
        \end{scope}
    
        
        \begin{scope}[xshift=6cm]
        \foreach \i in {0,...,5}
        {
        \draw[gray!50] (0,\i*\h)--(\w*5,\i*\h);
        }
        \foreach \i in {0,...,5}
        {
        \draw[gray!50] (\w*\i,0)--(\w*\i,5*\h);
        }
        \foreach \xx\yy in {0/1,3/1,4/1,1/2,4/2,0/3,2/3,3/4}
        {
        \draw (\w*.5+\w*\xx,\h*.5+\h*\yy) circle (\r cm);
        }
        \foreach \xx\yy in {0/0,2/0,3/0}
        {
        \fill[black] (\w*.5+\w*\xx,\h*.5+\h*\yy) circle (\r cm);
        }
        \node at (2.5,-2) {\scalebox{0.8}{\tableau{1&1&1\\}}};
        \end{scope}
    
        
        \begin{scope}[xshift=12cm]
        \foreach \i in {0,...,5}
        {
        \draw[gray!50] (0,\i*\h)--(\w*5,\i*\h);
        }
        \foreach \i in {0,...,5}
        {
        \draw[gray!50] (\w*\i,0)--(\w*\i,5*\h);
        }
        \foreach \xx\yy in {1/2,4/2,0/3,2/3,3/4}
        {
        \draw (\w*.5+\w*\xx,\h*.5+\h*\yy) circle (\r cm);
        }
        \foreach \xx\yy in {0/0,2/0,3/0,4/0,0/1,3/1}
        {
        \fill[black] (\w*.5+\w*\xx,\h*.5+\h*\yy) circle (\r cm);
        }
        \node at (2.5,-2) {\scalebox{0.8}{\tableau{2&2\\1&1&1&2}}};
        \end{scope}
    
        
        \begin{scope}[xshift=18cm]
        \foreach \i in {0,...,5}
        {
        \draw[gray!50] (0,\i*\h)--(\w*5,\i*\h);
        }
        \foreach \i in {0,...,5}
        {
        \draw[gray!50] (\w*\i,0)--(\w*\i,5*\h);
        }
        \foreach \xx\yy in {0/3,2/3,3/4}
        {
        \draw (\w*.5+\w*\xx,\h*.5+\h*\yy) circle (\r cm);
        }
        \foreach \xx\yy in {0/0,2/0,3/0,4/0,0/1,3/1,4/1,1/2}
        {
        \fill[black] (\w*.5+\w*\xx,\h*.5+\h*\yy) circle (\r cm);
        }
        
        \node at (2.5,-2.5) {\scalebox{0.8}{\tableau{3\\2&2&3\\1&1&1&2}}};
        \end{scope}
    
        
        \begin{scope}[xshift=24cm]
        \foreach \i in {0,...,5}
        {
        \draw[gray!50] (0,\i*\h)--(\w*5,\i*\h);
        }
        \foreach \i in {0,...,5}
        {
        \draw[gray!50] (\w*\i,0)--(\w*\i,5*\h);
        }
        \foreach \xx\yy in {3/4}
        {
        \draw (\w*.5+\w*\xx,\h*.5+\h*\yy) circle (\r cm);
        }
        \foreach \xx\yy in {0/0,2/0,3/0,4/0,0/1,3/1,4/1,1/2,2/2,0/3}
        {
        \fill[black] (\w*.5+\w*\xx,\h*.5+\h*\yy) circle (\r cm);
        }
        \node at (2.5,-2.5) {\scalebox{0.8}{\tableau{4\\3&4\\2&2&3\\1&1&1&2}}};
        \end{scope}

        \begin{scope}[xshift=30cm]
        \foreach \i in {0,...,5}
        {
        \draw[gray!50] (0,\i*\h)--(\w*5,\i*\h);
        }
        \foreach \i in {0,...,5}
        {
        \draw[gray!50] (\w*\i,0)--(\w*\i,5*\h);
        }
        \foreach \xx\yy in {0/0,2/0,3/0,4/0,0/1,3/1,4/1,1/1,2/2,0/3,3/2}
        {
        \fill[black] (\w*.5+\w*\xx,\h*.5+\h*\yy) circle (\r cm);
        }
        \node at (2.5,-2.5) {\scalebox{0.8}{\tableau{4\\3&4\\2&2&3&5\\1&1&1&2}}};
        \end{scope}
    \end{tikzpicture}
\end{center}

\end{example}

The collapsing procedure acts trivially on $M\in\MLQ_0(\lambda,n) \subseteq \M(\lambda_1,n)$. Thus we have the simple, but useful lemma below.

\begin{lemma}\label{lem:Q-nonwrapping-case}
Let $M\in \MLQ_0(\lambda)$. Then the tableau $\rho_Q(M)$ is the semistandard tableau with shape and content equal to $\lambda'$, and $\charge(M)=0$. 
\end{lemma}

By sequentially applying braid relations to the operators $\{e^\star_i\}$, we can equivalently write $\rho_N$ 
as follows. 

\begin{lemma}\label{prop:top-to-bot-collapse}
    Let $L$ and $n$ be positive integers and let $B=(B_1,\ldots,B_L)\in\M(L,n)$. Then 
    \[
    \rho_N(B) = e_{[L-1,L-1]}^\star e_{[L-2,L-1]}^\star \cdots e_{[2,L-1]}^\star e_{[1,L-1]}^\star(B).
    \]
\end{lemma}

\begin{proof}
    Recall from the proof of \cref{prop:rho} that $\rho_N^{(k)}=e_{[1,k-1]}^\star\cdots e_{[1,1]}^\star$ is the collapsing operator applied to the bottom $k$ rows. We will show by induction on $k$ that 
    \[\rho_N^{(k)} = e_{[k-1,k-1]}^\star e_{[k-2,k-1]}^\star \cdots  e_{[1,k-1]}^\star.\] 
    The identity is trivially true for $k=1$. Supposing the result holds for $\rho_N^{(r)}$ with $r\geq 2$, we will show it also holds for $\rho_N^{(r+1)}$. Using a sequence of commutation and braid relations from \cref{thm:dropping-operators-algebra}, we have that for $2\leq k\leq r$, 
\[
 e_{[1,r]}^\star e_{[k-1,r-1]}^\star=e_{[k,r]}^\star e_{[1,r]}^\star
\]
We apply this relation sequentially for $k=r,r-1,\ldots,2$ to obtain
\begin{align*}
        \rho_N^{(r+1)} = e_{[1,r]}^\star \rho_N^{(r)} 
        &= e_{[1,r]}^\star e_{[r-1,r-1]}^\star e_{[r-2,r-1]}^\star\cdots  e_{[1,r-1]}^\star
        = e_{[r,r]}^\star e_{[1,r]}^\star e_{[r-2,r-1]}^\star \cdots e_{[1,r-1]}^\star \\
        &= e_{[r,r]}^\star e_{[r-1,r]}^\star e_{[1,r]}^\star \cdots e_{[1,r-1]}^\star = \ldots = e_{[r,r]}^\star e_{[r-1,r]}^\star \cdots e_{[1,r]}^\star.
\end{align*}
Since $\rho_N(B)=\rho_N^{(L)}(B)$, we have the desired identity.
\end{proof}

\begin{remark}
    In \cref{def:collapsing_mlqs} we scan the rows of the multiline queue from bottom to top in each step of the collapsing. This is reflected in the order of the operators. In contrast, in \cref{prop:top-to-bot-collapse}, we start from the top and sequentially create the bottom rows of the nonwrapping multiline queue. Thus, we will refer to this version of collapsing as \emph{top-to-bottom collapsing}, while to the first definition we will refer as \emph{bottom-to-top collapsing} or simply \emph{collapsing.}
\end{remark}

\begin{remark}\label{rem:top-to-bot-adapated}  
    Following the discussion in \cref{rem:adapted-strings}, top-to-bottom collapsing corresponds to picking the reduced word for the longest permutation to be $$w_0 = (s_{L-1}s_{L-2}\cdots s_{1})(s_{L-1}s_{L-2}\cdots s_{2})\cdots(s_{L-1}s_{L-2})s_{L-1}.$$
\end{remark}

We describe an alternative way to obtain the recording tableau $\rho_Q(B)$, by keeping track of where particles from each row of $B$ end up in $\rho_N(B)$, which lets us store both the insertion and recording data in the same object; we call this \emph{labeled collapsing}.

\begin{defn}[Labeled collapsing]\label{def:labeled_collapsing}

Let $B = (B_1,\ldots,B_L)\in\M(L,n)$. For each $r$, assign the \emph{tracking label} $r$ to each particle in $B_r$. Next, perform the collapsing procedure (either top-to-bottom or bottom-to-top) on the labeled configuration, and reconfigure the tracking labels according to the following \emph{local} rule: if a particle $b$ with tracking label $\ell$ bumps a particle $b'$ with tracking label $k$ with $k<\ell$, swap the labels so that label $k$ stays in the current row; once the tracking labels are swapped, continue with the collapsing. After all particles are collapsed, the final configuration of balls is $\rho_N(B)$, and we construct the semistandard tableau $Q'(B)$ by reordering the tracking labels of the final configuration. 
\end{defn}

\begin{remark}
    The previous local rule can be restated in a more \emph{global} fashion to account for the collapsing of a full row on top of another, i.e. the application of $e_i^\star$ instead of a single $e_i$.

    For a given row $i$, say $B_i = \{b_1 < b_2 < \ldots < b_k\}$. For a particle $b\in B_i$, let $S_b$ be the set of particles that are in row $i+1$, weakly to the left of $b$. Then in the collapsing of row $i+1$, the tracking label of the particle paired to $b_1$ is $c_1 \coloneqq \min(S_{b_1})$ and, inductively, for $2\leq j \leq k$ the tracking label of the particle paired to $b_j$ is $c_j \coloneqq \min(S_{b_j}\setminus\{c_1,c_2,\ldots,c_{j-1}\}).$ To drop the remaining particles, we move the tracking label $c_j$ to the particle paired to $b_j$ and shift the tracking labels of the remaining particles to the left, maintaining the previous order. Then, the particles that are unmatched collapse to row $i$.
\end{remark}

See \cref{fig:labeled collapse} for an explicit example of the previous definition. We now show that this version of collapsing coincides with the one in \cref{def:collapsing_mlqs}. Since we use collapsing that is either bottom-to-top or top-to-bottom, \cref{prop:top-to-bot-collapse} already proves that the final ball arrangement in both cases is the same. Therefore we just have to check the construction of the recording tableau.

\begin{prop}\label{prop:labeled-collapsing-Q}
    For any ball arrangement $B\in\M(L,n)$, the tableau $Q'(B)$ from \cref{def:labeled_collapsing} obtained by collapsing coincides with $\rho_Q(B).$
\end{prop}

\begin{proof}

    We proceed by induction on $L$, and we analyze the cases of bottom-to-top and top-to-bottom collapsing separately. In both cases, the base case $L=2$ is trivial. Assume that labeled collapsing in produces the correct recording tableau for a ball arrangement of $L-1$ rows where $L \geq 3$. Suppose $B = (B_1,B_2,\ldots,B_L)$ has $L$ rows. By assumption, for $B^{(L-1)} = (B_1,B_2,\ldots,B_{L-1})$ we have that $Q'(B^{(L-1)}) = \rho_Q(B^{(L-1)})$.
    
    First we analyze bottom-to-top collapsing. In the last part of the labeled collapsing of $B_2 = (\rho_N(B^{(L-1)}),B_L)$, namely when $e_{[1,L-1]}^\star$ acts on the ball arrangement, the possible dropping of labels occur only when a particle with label "$L$" bumps another particle. Indeed, since all the labels appearing in $\rho_N(B^{(L-1)})$ are less than $L$, the dropped labels are always "$L$" unless they are the only label pairing to a given particle, in which case they stay in the given row. This shows that the rule in \cref{def:labeled_collapsing} is a reformulation of the recording tableau from \cref{def:collapsing_mlqs}. 

    Now we analyze top-to-bottom collapsing. Consider the partial collapsing $e_{[1,L-1]}^\star(B)$. Note that the set of labels that stay in each row during the labeled collapsing is the same as in $e_{[1,L-1]}^\star(B^{(L-1)})$ since the labels from $B_L$ are all "$L$". Thus, any ball in row $1$ that is in $e_{[1,L-1]}^\star(B)$ but not in $e_{[1,L-1]}^\star(B^{(L-1)})$ has label "$L$". The same argument shows that any ball that is in row $j$ of $e_{[j,L-1]}^\star\cdots e_{[1,L-1]}^\star(B)$ but not in $e_{[j,L-2]}^\star\cdots e_{[1,L-2]}^\star(B^{(L-1)})$ has label "$L$". Thus, the tableau $Q'(B)$ is also recording the new particles that end up in each row of the previously collapsed ball arrangement, which means that it coincides with the recording tableau from \cref{def:collapsing_mlqs}.
\end{proof}

\begin{example}\label{fig:labeled collapse}
We show the step-by-step labeled collapsing of the multiline queue from \cref{ex:collapse}. First, we consider the case of bottom-to-top collapsing.

\medskip

\begin{center}
    \begin{tikzpicture}[scale=0.4]
        \def \w{1};
        \def \h{1};
        \def \r{0.3};
        \begin{scope}[xshift=0cm]
        \foreach \i in {0,...,5}
        {
        \draw[gray!50] (0,\i*\h)--(\w*5,\i*\h);
        }
        \foreach \i in {0,...,5}
        {
        \draw[gray!50] (\w*\i,0)--(\w*\i,5*\h);
        }
        \foreach \xx\yy\c in {0/0/1,2/0/1,3/0/1,0/1/2,3/1/2,4/1/2,1/2/3,4/2/3,0/3/4,2/3/4,3/4/5}
        {
       \filldraw[black,fill=red!10] (\w*.5+\w*\xx,\h*.5+\h*\yy) circle (\r cm);
        \node at (\w*.5+\w*\xx,\h*.5+\h*\yy) {\tiny \c};
        }
        \node at (2.5,-2) {\scalebox{0.8}{\tableau{\emptyset\\}}};
        \end{scope}
    
        
        \begin{scope}[xshift=6cm]
        \foreach \i in {0,...,5}
        {
        \draw[gray!50] (0,\i*\h)--(\w*5,\i*\h);
        }
        \foreach \i in {0,...,5}
        {
        \draw[gray!50] (\w*\i,0)--(\w*\i,5*\h);
        }
        \foreach \xx\yy\c in {0/1/2,3/1/2,4/1/2,1/2/3,4/2/3,0/3/4,2/3/4,3/4/5}
        {
        \filldraw[black,fill=red!10] (\w*.5+\w*\xx,\h*.5+\h*\yy) circle (\r cm);
        \node at (\w*.5+\w*\xx,\h*.5+\h*\yy) {\tiny \c};
        }
        \foreach \xx\yy\c in {0/0/1,2/0/1,3/0/1}
        {
        \draw (\w*.5+\w*\xx,\h*.5+\h*\yy) circle (\r cm);
        \node at (\w*.5+\w*\xx,\h*.5+\h*\yy) {\tiny \c};
        }
        \node at (2.5,-2) {\scalebox{0.8}{\tableau{1&1&1\\}}};
        \end{scope}
    
        
        \begin{scope}[xshift=12cm]
        \foreach \i in {0,...,5}
        {
        \draw[gray!50] (0,\i*\h)--(\w*5,\i*\h);
        }
        \foreach \i in {0,...,5}
        {
        \draw[gray!50] (\w*\i,0)--(\w*\i,5*\h);
        }
        \foreach \xx\yy\c in {1/2/3,4/2/3,0/3/4,2/3/4,3/4/5}
        {
       \filldraw[black,fill=red!10] (\w*.5+\w*\xx,\h*.5+\h*\yy) circle (\r cm);
        \node at (\w*.5+\w*\xx,\h*.5+\h*\yy) {\tiny \c};
        }
        \foreach \xx\yy\c in {0/0/1,2/0/1,3/0/1,4/0/2,0/1/2,3/1/2}
        {
        \draw (\w*.5+\w*\xx,\h*.5+\h*\yy) circle (\r cm);
        \node at (\w*.5+\w*\xx,\h*.5+\h*\yy) {\tiny \c};
        }
        \node at (2.5,-2) {\scalebox{0.8}{\tableau{2&2\\1&1&1&2}}};
        \end{scope}
    
        
        \begin{scope}[xshift=18cm]
        \foreach \i in {0,...,5}
        {
        \draw[gray!50] (0,\i*\h)--(\w*5,\i*\h);
        }
        \foreach \i in {0,...,5}
        {
        \draw[gray!50] (\w*\i,0)--(\w*\i,5*\h);
        }
       \foreach \xx\yy\c in {0/3/4,2/3/4,3/4/5}
        {
       \filldraw[black,fill=red!10] (\w*.5+\w*\xx,\h*.5+\h*\yy) circle (\r cm);
        \node at (\w*.5+\w*\xx,\h*.5+\h*\yy) {\tiny \c};
        }
        \foreach \xx\yy\c in {0/0/1,2/0/1,3/0/1,4/0/2,0/1/2,3/1/2,4/1/3,1/2/3}
        {
        \draw (\w*.5+\w*\xx,\h*.5+\h*\yy) circle (\r cm);
        \node at (\w*.5+\w*\xx,\h*.5+\h*\yy) {\tiny \c};
        }
        
        \node at (2.5,-2.5) {\scalebox{0.8}{\tableau{3\\2&2&3\\1&1&1&2}}};
        \end{scope}
    
        
        \begin{scope}[xshift=24cm]
        \foreach \i in {0,...,5}
        {
        \draw[gray!50] (0,\i*\h)--(\w*5,\i*\h);
        }
        \foreach \i in {0,...,5}
        {
        \draw[gray!50] (\w*\i,0)--(\w*\i,5*\h);
        }
        \foreach \xx\yy\c in {3/4/5}
        {
       \filldraw[black,fill=red!10] (\w*.5+\w*\xx,\h*.5+\h*\yy) circle (\r cm);
        \node at (\w*.5+\w*\xx,\h*.5+\h*\yy) {\tiny \c};
        }
        \foreach \xx\yy\c in {0/0/1,2/0/1,3/0/1,4/0/2,0/1/2,3/1/2,4/1/3,1/2/4,2/2/3,0/3/4}
        {
        \draw (\w*.5+\w*\xx,\h*.5+\h*\yy) circle (\r cm);
        \node at (\w*.5+\w*\xx,\h*.5+\h*\yy) {\tiny \c};
        }
        \node at (2.5,-2.5) {\scalebox{0.8}{\tableau{4\\3&4\\2&2&3\\1&1&1&2}}};
        \end{scope}

        \begin{scope}[xshift=30cm]
        \foreach \i in {0,...,5}
        {
        \draw[gray!50] (0,\i*\h)--(\w*5,\i*\h);
        }
        \foreach \i in {0,...,5}
        {
        \draw[gray!50] (\w*\i,0)--(\w*\i,5*\h);
        }
        \foreach \xx\yy\c in {0/0/1,2/0/1,3/0/1,4/0/2,0/1/2,3/1/2,4/1/3,1/1/5,2/2/4,0/3/4,3/2/3}
        {
        \draw (\w*.5+\w*\xx,\h*.5+\h*\yy) circle (\r cm);
        \node at (\w*.5+\w*\xx,\h*.5+\h*\yy) {\tiny \c};
        }
        \node at (2.5,-2.5) {\scalebox{0.8}{\tableau{4\\3&4\\2&2&3&5\\1&1&1&2}}};
        \end{scope}
    \end{tikzpicture}
\end{center} 
\noindent On the other hand, with top-to-bottom collapsing we have the following. 

\medskip

\begin{center}
    \begin{tikzpicture}[scale=0.4]
        \def \w{1};
        \def \h{1};
        \def \r{0.3};
        \begin{scope}[xshift=0cm]
        \foreach \i in {0,...,5}
        {
        \draw[gray!50] (0,\i*\h)--(\w*5,\i*\h);
        }
        \foreach \i in {0,...,5}
        {
        \draw[gray!50] (\w*\i,0)--(\w*\i,5*\h);
        }
        \foreach \xx\yy\c in {0/0/1,2/0/1,3/0/1,0/1/2,3/1/2,4/1/2,1/2/3,4/2/3,0/3/4,2/3/4,3/4/5}
        {
       \filldraw[black,fill=red!10] (\w*.5+\w*\xx,\h*.5+\h*\yy) circle (\r cm);
        \node at (\w*.5+\w*\xx,\h*.5+\h*\yy) {\tiny \c};
        }
        \node at (2.5,-2) {\scalebox{0.8}{\tableau{\emptyset\\}}};
        \end{scope}
    
        
        \begin{scope}[xshift=6cm]
        \foreach \i in {0,...,5}
        {
        \draw[gray!50] (0,\i*\h)--(\w*5,\i*\h);
        }
        \foreach \i in {0,...,5}
        {
        \draw[gray!50] (\w*\i,0)--(\w*\i,5*\h);
        }
        \foreach \xx\yy\c in {0/1/2,3/1/2,2/2/3,4/2/3,0/3/4,3/3/4,1/1/5}
        {
        \filldraw[black,fill=red!10] (\w*.5+\w*\xx,\h*.5+\h*\yy) circle (\r cm);
        \node at (\w*.5+\w*\xx,\h*.5+\h*\yy) {\tiny \c};
        }
        \foreach \xx\yy\c in {0/0/1,2/0/1,3/0/1,4/0/2}
        {
        \draw (\w*.5+\w*\xx,\h*.5+\h*\yy) circle (\r cm);
        \node at (\w*.5+\w*\xx,\h*.5+\h*\yy) {\tiny \c};
        }
        \node at (2.5,-2) {\scalebox{0.8}{\tableau{1&1&1&2\\}}};
        \end{scope}
    
        
        \begin{scope}[xshift=12cm]
        \foreach \i in {0,...,5}
        {
        \draw[gray!50] (0,\i*\h)--(\w*5,\i*\h);
        }
        \foreach \i in {0,...,5}
        {
        \draw[gray!50] (\w*\i,0)--(\w*\i,5*\h);
        }
        \foreach \xx\yy\c in {2/2/3,0/3/4,3/3/4}
        {
        \filldraw[black,fill=red!10] (\w*.5+\w*\xx,\h*.5+\h*\yy) circle (\r cm);
        \node at (\w*.5+\w*\xx,\h*.5+\h*\yy) {\tiny \c};
        }
        \foreach \xx\yy\c in {0/0/1,2/0/1,3/0/1,4/0/2,3/1/2,0/1/2,1/1/5,4/1/3}
        {
        \draw (\w*.5+\w*\xx,\h*.5+\h*\yy) circle (\r cm);
        \node at (\w*.5+\w*\xx,\h*.5+\h*\yy) {\tiny \c};
        }
        \node at (2.5,-2) {\scalebox{0.8}{\tableau{2&2&3&5\\1&1&1&2}}};
        \end{scope}
    
        
        \begin{scope}[xshift=18cm]
        \foreach \i in {0,...,5}
        {
        \draw[gray!50] (0,\i*\h)--(\w*5,\i*\h);
        }
        \foreach \i in {0,...,5}
        {
        \draw[gray!50] (\w*\i,0)--(\w*\i,5*\h);
        }
        \foreach \xx\yy\c in {0/3/4}
        {
        \filldraw[black,fill=red!10] (\w*.5+\w*\xx,\h*.5+\h*\yy) circle (\r cm);
        \node at (\w*.5+\w*\xx,\h*.5+\h*\yy) {\tiny \c};
        }
        \foreach \xx\yy\c in {0/0/1,2/0/1,3/0/1,4/0/2,3/1/2,0/1/2,1/1/5,4/1/3,3/2/3,2/2/4}
        {
        \draw (\w*.5+\w*\xx,\h*.5+\h*\yy) circle (\r cm);
        \node at (\w*.5+\w*\xx,\h*.5+\h*\yy) {\tiny \c};
        }
        
        \node at (2.5,-2.5) {\scalebox{0.8}{\tableau{3&4\\2&2&3&5\\1&1&1&2}}};
        \end{scope}
    
        
        \begin{scope}[xshift=24cm]
        \foreach \i in {0,...,5}
        {
        \draw[gray!50] (0,\i*\h)--(\w*5,\i*\h);
        }
        \foreach \i in {0,...,5}
        {
        \draw[gray!50] (\w*\i,0)--(\w*\i,5*\h);
        }
        \foreach \xx\yy\c in {}
        {
        \filldraw[black,fill=red!10] (\w*.5+\w*\xx,\h*.5+\h*\yy) circle (\r cm);
        \node at (\w*.5+\w*\xx,\h*.5+\h*\yy) {\tiny \c};
        }
        \foreach \xx\yy\c in {0/0/1,2/0/1,3/0/1,4/0/2,3/1/2,0/1/2,1/1/5,4/1/3,3/2/3,2/2/4,0/3/4}
        {
        \draw (\w*.5+\w*\xx,\h*.5+\h*\yy) circle (\r cm);
        \node at (\w*.5+\w*\xx,\h*.5+\h*\yy) {\tiny \c};
        }
        
        \node at (2.5,-2.5) {\scalebox{0.8}{\tableau{4\\3&4\\2&2&3&5\\1&1&1&2}}};
        \end{scope}
    \end{tikzpicture}
\end{center}
\end{example}

\begin{theorem}\label{thm:wrapping-to-nonwrapping}
    The collapsing procedure from \cref{def:collapsing_mlqs} gives a 
    bijection such that $x^B=x^{\rho_N(B)}$ for $B\in\M(L,n)$. 
\end{theorem} 

\begin{proof}
Let $B\in\M(L,n)$ and let $(N,Q)=\rho(B)$. The equality $x^B=x^N$ 
is due to the fact that the collapsing procedure leaves the column content (and hence the $x$-weight) invariant. We will show $\rho$ is a bijection by constructing an inverse using the invertibility of the dropping operators in \eqref{eq:dropping}.

For $1\leq r< L-1$ with $L=\lambda_1$, let $B^{(r-1)}=e^{\star}_{[1,r-2]}\circ\cdots\circ e^{\star}_{[1,1]}(B)$ be a partially collapsed multiline queue. Then $e^{\star}_{[1,r-1]}(B^{(r)})$ generates the collapsing of row $r$. To construct the inverse, one needs to keep track of the number of times each operator $e_i$ is applied within $e^{\star}_{[1,r-1]}(B^{(r-1)})$. For $1\leq j\leq r-1$, let $\phi(r,j)$ be minimal such that $e^{\star}_{[j,r-1]}(B^{(r-1)})=e^{\phi(r,j)}_j(e^{\star}_{[j+1,r-1]}(B^{(r-1)}))$. Then
\begin{equation}\label{eq:Mr}
B^{(r)}=e^{\star}_{[1,r-1]}(B^{(r-1)})=e^{\phi(r,1)}_1\circ\cdots\circ e^{\phi(r,r-1)}_{r-1}(B^{(r-1)}).
\end{equation}
In particular, the sequence $(\phi(r,1),\ldots,\phi(r,r-1))$ dictates the multiset of rows containing the entries ``$r$'' in the recording tableau as follows: if $(m^{(r)}_r,\ldots,m^{(r)}_1)$ is the multiplicity vector with $m^{(r)}_j$ equal to the number of entries ``$r$'' in row $j$ of $Q$, then $m^{(r)}_j=\phi(r,j)-\phi(r,j-1)$, where $\phi(r,r)=\lambda_r'$ and $\phi(r,0)=0$.

By minimality of the $\phi(r,j)$'s, each $e_j^{\phi(r,j)}$ is invertible with inverse $f_j^{\phi(r,j)}$, and thus
\[f^{\phi(r,r-1)}_{r-1}\circ\cdots\circ f^{\phi(r,1)}_{1}(M_{r})=M_{r-1}.\]
Therefore, to construct $\rho^{-1}(N,Q)$ from $(N,Q)\in \MLQ_0(\mu,n) \times \SSYT(\mu',\lambda')$ for some partition $\mu$, we define the tuples $(\phi(r,j)\colon1\leq j<r)$ from the multiplicity vectors  $(m_j^{(r)}\colon1\leq j\leq r)$ in $Q$ for $2\leq r\leq L$ , and set
\begin{equation}\label{eq:inverse-collapsing-mlq-ssyt}
    \rho^{-1}(N,Q) = \bigcirc_{r=2}^{L}\left( f^{\phi(r,r-1)}_{r-1}
     \circ\cdots\circ f^{\phi(r,2)}_2\circ f^{\phi(r,1)}_1 \right)(N).
\end{equation}
where the notation $\bigcirc_{i=a}^b Y_i (N)$ represents the composition $Y_b\circ\cdots\circ Y_{a+1}\circ Y_a(N)$. 

By \eqref{eq:Mr}, the composition of \eqref{eq:inverse-collapsing-mlq-ssyt} with \eqref{eq:dropping} is the identity, 
so $\rho^{-1}$ is a left inverse of $\rho$. A similar argument shows this is also a right inverse of $\rho$.
\end{proof}

In fact, the bijection $\rho$ can be directly considered an analogue of RSK on multiline queues, since the lifting operators required to recover the corresponding multiline queue are read from the tableau $Q$, in the same way as one would construct the inverse RSK map. In \cref{sec:mlq/d_RSK}, we will interpret collapsing as a map from multiline queues to pairs of nonwrapping multiline queues, to strengthen the comparison with RSK.

We restrict the collapsing procedure to the set of multiline queues $\MLQ(\lambda,n)$, identified as the set of binary matrices on $n$ columns with row content given by $\lambda'$. Endowing this set with the weight $\wt(M) = q^{\maj(M)}x^M$ gives the following map. In particular, with this map, we recover the expansion of $q$-Whittaker polynomials in the Schur basis. 

\begin{theorem}\label{thm:collapsing_mlqs}

    Let $\lambda$ be a partition and $n \geq \ell(\lambda)$ a positive integer. Then collapsing restricted to $\MLQ(\lambda,n)$ is a weight-preserving bijection with $x^{|M|}=x^{\rho_N(M)}$ and $\maj(M)=\charge(\rho_Q(M))$:
    \begin{equation}
        \rho
        \; 
        \colon
        \; \MLQ(\lambda,n) \;\to\; \bigcup_{\mu\leq \lambda} \MLQ_0(\mu,n)\times \SSYT(\mu',\lambda')\\
    \end{equation}
\end{theorem}

We end this section by examining some properties of collapsing on twisted multiline queues. 

\begin{prop}\label{prop:collapsing and sigma}
For a composition $\alpha$ and $B\in\GMLQ(\alpha)$,  $\rho(B)=\rho(\sigma_i(B))$ for any $i\geq 1$.
\end{prop}

\begin{proof}
Write $B'=\sigma_i(B)$. Without loss of generality, assume $\alpha_{i+1}>\alpha_i$ (if $\alpha_i=\alpha_{i+1}$, the claim is trivial since $B=B'$). 

Let $A\subseteq B_{i+1}$ and $C\subseteq B_{i}$ be the sets of particles matched above and below, respectively, in $\Par^c(B_{i+1},B_{i})$. Note that these are also the sets matched above and below in $\Par^c(B'_{i+1},B'_{i})$ by the definition of $\sigma_i$. Let $D$ be the set of particles that is moved between rows $i$ and $i+1$ by the involution $\sigma_i$ (i.e. the set of unmatched particles above in $\Par^c(B_{i+1},B_{i})$ and below in $\Par^c(B'_{i+1},B'_{i})$, respectively). Then $B_i=C$, $B_{i+1}=A\cup D$, $B'_i=C\cup D$, and $B'_{i+1}=A$. 

Balls matched above (resp. below) in $\Par(X,Y)$ are necessarily also matched above (resp. below) in $\Par^c(X,Y)$. Thus a ball is matched above in $\Par(A,C\cup D)$ only if it is matched above in $\Par^c(A,C\cup D)$, which is true only if it is matched above in $\Par^c(A\cup D,C)$. Since no ball in $D$ is matched, every ball matched above in $\Par(A\cup D,C)$ must also be matched above in $\Par(A,C\cup D)$. Thus the set of balls matched above in $\Par(B'_{i+1},B'_i)$ is equal to the set of balls matched above in $\Par(B_{i+1},B_i)$, which implies that $e_i^\star(\sigma_i(B))=e_i^\star(B)$. 

Now we write $\rho=e_{[1,L]}^\star\cdots e_{[1,i+1]}^\star \rho^{(i+1)}_N$. Since the first operator applied in $\rho^{(i+1)}_N$ is $e_i^\star$, from the above we have that $\rho^{(i+1)}_N\circ\sigma_i(B)=\rho^{(i+1)}_N(B)$, from which we get $\rho(\sigma_i(B))=\rho(B)$.
\end{proof}

\subsubsection{Nonwrapping twisted multiline queues}
It is natural to define the analogue of $\MLQ_0(\lambda)$ in the twisted case.
\begin{defn}\label{def:gmlq 0}
For a composition $\alpha$, the set of \emph{nonwrapping twisted multiline queues} is
\[
\GMLQ_0(\alpha)=\{B\in\GMLQ(\alpha):\maj_G(B)=0\}.
\]
\end{defn}
We claim that $B\in \GMLQ_0(\alpha,n)$ if and only if it has no wrapping pairings under the left-to-right pairing order convention of \cref{def:gmlq_particlewise}, hence justifying the choice of the name for these objects. We give a sketch of the argument: first, for any twisted multiline queue, within each pair of rows, the number of particle pairings wrapping to the right must equal the number of anti-particle pairings wrapping to the left. Second, if $\maj_G(B)=0$, the contribution to $\maj_G$ from each pair of rows is 0. Since the label of any anti-particle is less than or equal to the label of any particle, in order to get a sum of zero from the particle and anti-particle contributions within each pair of rows, the only possibility is that all wrapping pairings come from particles or anti-particles of the same label (i.e. the smallest label pairing in the particle phase), which necessarily cancel each other out. However, this means the left to right pairing order during the particle phase will necessarily prevent any wrapping particle pairings. In particular, this last fact implies that to check that a twisted multiline queue has $\maj_G$ equal to 0, it is sufficient to only check for particle pairings wrapping to the right.

From \cref{prop:collapsing and sigma}, we obtain the  following commuting diagram. 
\begin{center}
    \begin{tikzpicture}
        \node[] (GMLQ) at (0,0) {$\GMLQ(\sigma\cdot\lambda)$};
        \node[] (GMLQ0) at (0,-2) {$\GMLQ_0(\sigma\cdot\mu)$};
        \node[] (MLQ) at (4,0) {$\MLQ(\lambda)$};
        \node[] (MLQ0) at (4,-2) {$\MLQ_0(\mu)$};

        \draw[->] (GMLQ) -- (GMLQ0) node[left, midway] {$\rho_\sigma^\downarrow$};
        \draw[->] (MLQ) -- (GMLQ) node[above, midway] {$\sigma$};
        \draw[->] (MLQ) -- (MLQ0) node[right, midway] {$\rho_N^\downarrow$};
        \draw[->] (MLQ0) -- (GMLQ0) node[above, midway] {$\sigma$};
        \draw[->] (GMLQ) -- (MLQ0) node[above right, midway] {$\rho_N^\downarrow$};
        
    \end{tikzpicture}
\end{center}

The composition $\sigma\circ \rho_N^{\downarrow}\circ \sigma^{-1}$ defines a map $\rho_\sigma^\downarrow \colon\GMLQ(\alpha) \to \bigcup_{\beta}\GMLQ_0(\beta)$ where the union runs over all compositions $\beta$ such that $\sigma^{-1}\cdot \beta \leq \sigma^{-1}\cdot\alpha$ with respect to dominance order on partitions. For a concrete example, see \cref{ex:commuting}. 

\begin{question}
    Give a combinatorial description of the map $\rho_\sigma^\downarrow$ on twisted multiline queues that is analogous to \cref{def:collapsing_mlqs}. 
\end{question}

\begin{example}\label{ex:commuting}

For $\lambda=(4,4,3,2,1)$ and $\sigma = 4\,2\,1\,5\,3 = s_3\,s_2\,s_1\,s_2\,s_4$, the twisting of the partition $\lambda$ is the composition $\alpha = \sigma\cdot\lambda=(3,4,1,4,2)$. For $M \in\GMLQ(\lambda,6)$, in \cref{fig:gmlq_square} we show the twisted multiline queue $B = \sigma(M) \in\GMLQ(\alpha,6)$, and the corresponding collapsed twisted multiline queues $\rho^\downarrow_\sigma(B)$ and $\rho^\downarrow(M)$. Note that $\rho_N^{\downarrow}(B) = \rho_N^{\downarrow}(M)$ in accordance with \cref{prop:collapsing and sigma}.

\begin{figure}
    \centering
    \begin{tikzpicture}[scale=1]
            \def \scale{0.4} 
            \node at (-3.25,0) {$B = $};
            \node at (9,0) {$=M$};
            \node at (-3.25,-5) {$\rho_\sigma^{\downarrow}(B)=$};
            \node at (9.25,-5) {$=\rho^{\downarrow}(M)$};
            \node (left) at (6,0) {\rotatebox{0}{\includegraphics[scale = \scale]{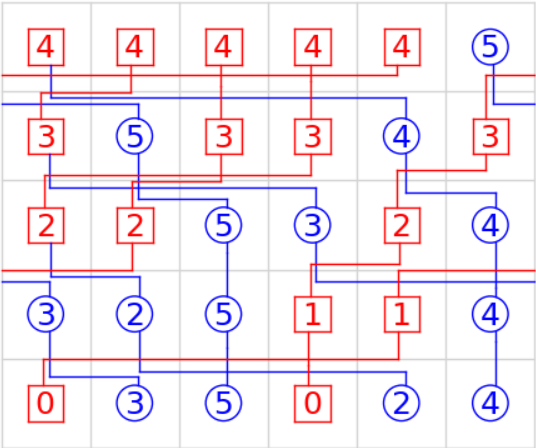}}};
            \node (original) at (0,0) {\includegraphics[scale = \scale]{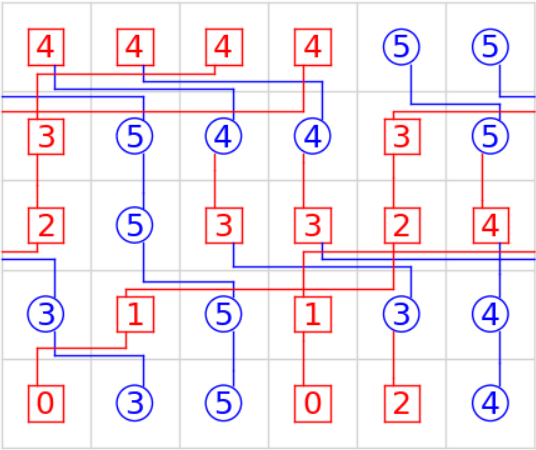}};
            \node (down) at (0,-5) {\rotatebox{0}{\includegraphics[scale = \scale]{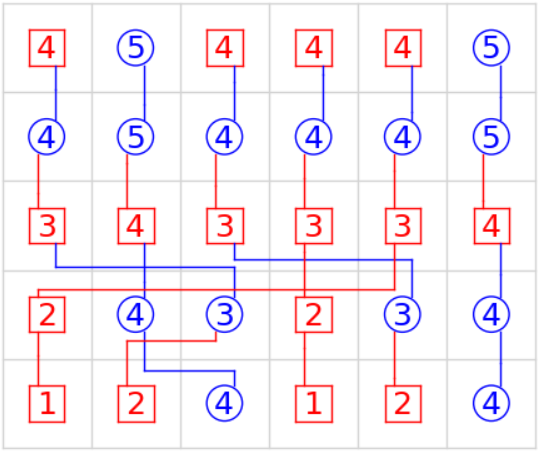}}};
            \node (double) at (6,-5) {\includegraphics[scale = \scale]{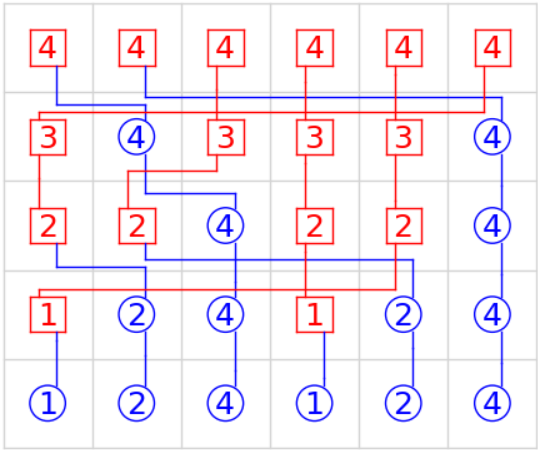}};
            \draw[->] (left) -- (original) node[above, midway] {\footnotesize{$\sigma$}};
            \draw[->] (original) -- (down) node[right, midway] {\footnotesize{$\rho_\sigma^\downarrow$}};
            \draw[->] (double) -- (down) node[above, midway] {\footnotesize{$\sigma$}};
            \draw[->]  (left) -- (double) node[right, midway] {\footnotesize{$\rho_N^\downarrow$}};
            \draw[->]  (original) -- (double) node[above right, midway] {\footnotesize{$\rho_N^\downarrow$}};
        \end{tikzpicture}  
    \caption{The twisting of a multiline queue by the permutation $\sigma = 4\,2\,1\,5\,3$ (written in one-line notation) together with the corresponding collapsing results.}
    \label{fig:gmlq_square}
\end{figure}

\end{example}

\subsection{Bijection from $\MLQ_0(\lambda,n)$ to $\SSYT(\lambda)$}\label{sec:bijection_mlq0_ssyt}
In this section, we establish a weight-preserving bijection between semistandard tableaux and nonwrapping multiline queues, which gives an alternate proof of \eqref{eq:schur-mlqs}. 
Denote by $\mathcal{W}$ the set of words in the alphabet $\mathbb{N}$, and the reverse of $w_1\ldots w_k\in \mathcal{W}$ by $\rev(w_1\ldots w_k)=w_k\ldots w_1$. We write $\revcrw(T)\coloneqq \rev(\crw(T))$ to mean the \emph{reverse column reading word} of $T \in \SSYT$ (see \cref{def:readingorder}). 

\begin{defn}
    For $w=w_1\ldots w_k\in\mathcal{W}$, denote by $\Icol(w)\in\SSYT$ the \emph{column} insertion of $w$ in the empty tableau:
    \[
\Icol(w)=w_k\rightarrow(\cdots\rightarrow(w_2\rightarrow(w_1\rightarrow \emptyset))\cdots).
    \]
\end{defn}

See \cite[Section A.2]{fulton-young} and the references therein for a complete description of column insertion. Using collapsing, we define an insertion procedure on nonwrapping multiline queues.

\begin{definition}
    Let $N=(N_1,\ldots,N_L)$ be a nonwrapping multiline queue on $n$ columns and let $k\in\{1,2\ldots,n\}$. The \emph{insertion of $k$ into $N$} is given by 
    \[
    k\rightarrow N \coloneqq \rho_N(N'),\qquad N'=(N_1,\ldots,N_L,\{k\}).\]
\end{definition}

For $w\in\mathcal{W}$, denote by $\rho(w)$ the nonwrapping multiline queue obtained by sequentially inserting the entries $w_1,\ldots,w_k$ into an empty multiline queue. We justify this slight abuse of notation by defining the multiline queue $M_w = (\{w_1\},\{w_2\},\ldots\{w_k\}) \in \MLQ((k),n)$, so that $\rho(w)=\rho_N(M_w)$.

In particular, checking the properties of column insertion and collapsing of nonwrapping multiline queues, respectively, we have the following.

\begin{lemma}\label{lemma:property-rho-and-Icol}
    For $T \in \SSYT(\lambda)$ and $N \in \MLQ_0(\lambda)$, $\Icol(\revcrw(T)) = T$ and $\rho(\rw(N)) = N.$ 
\end{lemma}

Our main result of this section is the bijection below, illustrated in \cref{ex:rho bijection}.

\begin{theorem}\label{thm:bijectionSSYT-MLQ0}
    Let $\lambda$ be a partition. The following maps are inverses and are content-preserving:
    \begin{align*}
    \mlq \;& \coloneqq \; \rho\circ\revcrw \, \colon \, \SSYT(\lambda) \to \MLQ_0(\lambda)\\
    \tab \;& \coloneqq \; \Icol\circ\rw \, \colon \, \MLQ_0(\lambda) \to \SSYT(\lambda)
    \end{align*}
\end{theorem}

\begin{example}\label{ex:rho bijection}
The tableau $T\in\SSYT(7,6,3)$ shown below has reversed column word $$\revcrw(T)=5\,|\,3\,4\,|\,3\,4\,|\,2\,3\,|\,1\,2\,5\,|\,1\,2\,4\,|\,1\,2\,4.$$
The corresponding multiline queue is $M=\rho(\revcrw(T))\in\MLQ(7,6,3)$, with row word $$\rw(M) = 3\,4\,5\,|\,2\,3\,5\,|\,1\,3\,4\,|\,2\,4\,|\,1\,4\,|\,1\,2\,|\,2.$$
One may check that $\Icol(\rw(M))=T$. 

\begin{center}
$T=\quad$\raisebox{6pt}{\tableau{4&4&5\\2&2&2&3&4&4\\1&1&1&2&3&3&5}}
$\qquad \xrightarrow{\mlq} \qquad $ 
\raisebox{-52pt}{\includegraphics[scale = 0.4]{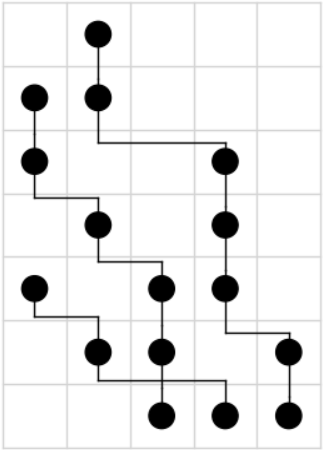}}$\quad=M$
\end{center}

\end{example}

The rest of the section is devoted to the proof of \cref{thm:bijectionSSYT-MLQ0}. We start with lemmas on the relation of Knuth equivalence to insertion on multiline queues. Recall that two words $w,w'\in\mathcal{W}$ are \emph{Knuth-equivalent} if $w$ can be transformed to $w'$ by means of the following \emph{elementary Knuth transformations} acting on triples of letters:
\begin{center}
    $bac \longmapsto bca$ for $a<b\leq c$  $\qquad$ and $\qquad$ $acb \longmapsto cab$ for $a\leq b < c$.
\end{center}
The following is a well-known close relation between column insertion and Knuth equivalence. 
\begin{prop}[{\cite[Prop. 2.3.14]{butler1994subgroup}}]
    Two words $w,w'\in\mathcal{W}$ have the same column-insertion tableau, i.e. $\Icol(w) = \Icol(w')$, if and only if $\rev(w)$ and $\rev(w')$ are Knuth-equivalent.
\end{prop}

\begin{lemma}\label{lemma:knuth-equivalence-collapsing}
    For $w\in\mathcal{W}$, $\rev(\rw(\rho(w)))$ is Knuth-equivalent to $\rev(w)$.
\end{lemma}
\begin{proof}
    Since $\rho$ can be built from the operators $e_i$, it is enough to show that $\rev(\rw(e_i(B)))$ and $\rev(\rw(B))$ are Knuth equivalent for any binary matrix $B$. Moreover, let us assume that $B$ only has two rows and $i=1$. 

    If $e_i(B)=B$, the statement is trivial, so suppose $x$ is the leftmost unmatched above ball in row $i+1$ of $B$. Let $\{a_1,\ldots,a_k\}$ be the set of balls to the left of $x$ in row $i+1$ with $a_1<\cdots<a_k<x$, and let $\{b_1,\ldots,b_\ell\}$ be the set of balls in row $i$ that are between $a_1$ (inclusive) and $x$: $a_1\leq b_1<\cdots<b_\ell<x$; the balls $\{a_1,\ldots,a_k\}$ are thus matched above to some set of balls $\{b_{i_1},\ldots,b_{i_k}\}\subseteq\{b_1,\ldots,b_\ell\}$. Similarly, let $\{d_1,\ldots,d_m\}$ be the set of balls to the right of $x$ in row $i$ with $x<d_1<\cdots<d_m$, and let $\{c_1,\ldots,c_n\}$ be the set of balls in row $i+1$ that are between $x$ and $d_m$ (inclusive): $x<c_1<\cdots<c_n\leq d_m$; the balls $\{d_1,\ldots,d_m\}$ are thus matched below to some set of balls $\{c_{j_1},\ldots,c_{j_m}\}\subseteq\{c_1,\ldots,c_n\}$. Any other balls (left of $b_1$ in row $i$ or right of $c_n$ in row $i+1$) are irrelevant to the calculations that follow, so we will not record them. Schematically, $B$ has the following structure:

    \begin{center}
        \begin{tikzpicture}[scale=0.7]
            \def \w{1};
            \def \h{1};
            \def \r{0.25};
            
            \foreach \i in {0,...,2}
            {
            \draw[gray!50] (0,\i*\h)--(\w*19,\i*\h);
            }
            \foreach \i in {0,...,19}
            {
            \draw[gray!50] (\w*\i,0)--(\w*\i,2*\h);
            }
            \foreach \xx\yy in {0/0,2/0,4/0,5/0,7/0,12/0,13/0,15/0}
            {
                \fill (\w*.5+\w*\xx,\h*.5+\h*\yy) circle (\r cm);
            }
            
            \foreach \xx\yy in {3/1,4/1,6/1,11/1,12/1,14/1,16/1,18/1}
            {
                \fill (\w*.5+\w*\xx,\h*.5+\h*\yy) circle (\r cm);
        
            }
        
            \foreach \xx\yy in {9/1}
            {
            \draw[thick] (\w*.5+\w*\xx,\h*.5+\h*\yy) circle (\r cm);
            }
        
            \foreach \xx\yy in {1/0,6/0,14/0,5/1,13/1,17/1}
            {
            \node at (\w*.5+\w*\xx,\h*.5+\h*\yy) {\textbf{$\ldots$}};
            }

            \node at (\w*.5+\w*4,\h*.5-\h) {$b_1$};
            \node at (\w*.5+\w*5,\h*.5-\h) {$b_2$};
            \node at (\w*.5+\w*7,\h*.5-\h) {$b_\ell$};
        
            \node at (\w*.5+\w*3,\h*.5+\h*2) {$a_1$};
            \node at (\w*.5+\w*4,\h*.5+\h*2) {$a_2$};
            \node at (\w*.5+\w*6,\h*.5+\h*2) {$a_k$};
        
            \node at (\w*.5+\w*9,\h*.5+\h*2) {$x$};
        
            \node at (\w*.5+\w*12,\h*.5-\h) {$d_1$};
            \node at (\w*.5+\w*13,\h*.5-\h) {$d_2$};
            \node at (\w*.5+\w*15,\h*.5-\h) {$d_m$};
        
            \node at (\w*.5+\w*11,\h*.5+\h*2) {$c_1$};
            \node at (\w*.5+\w*12,\h*.5+\h*2) {$c_2$};
            \node at (\w*.5+\w*14,\h*.5+\h*2) {$c_n$};
        \end{tikzpicture}
    \end{center}

     \noindent Observe that in $e_i(B)$, $x$ drops from row $i+1$ to row $i$. Then we have (ignoring irrelevant balls):
    \begin{align*}
    \rev(\rw(B))& = c_n \, c_{n-1} \, \ldots \, c_1 \, \mathbf{x} \, a_k \, a_{k-1} \, \ldots \, a_1 \, d_m \, d_{m-1} \, \ldots \, d_1 \, b_\ell \, b_{\ell-1} \, \ldots \, b_1
    &\coloneqq \textbf{c}\,x\,\textbf{a}\,\textbf{d}\,\textbf{b}\\ 
    \rev(\rw(e_1(B)))& = c_n \, c_{n-1} \, \ldots \, c_1 \, a_k \, a_{k-1} \, \ldots \, a_1 \, d_m \, d_{m-1} \, \ldots \, d_1  \, \mathbf{x} \, b_\ell \, b_{\ell-1} \, \ldots \, b_1 
    &\coloneqq \textbf{c}\,\textbf{a}\,\textbf{d}\,x\,\textbf{b}
    \end{align*}
    In fact, from the picture above, we have the following Knuth relations:\\
    \begin{itemize}
        \item $a_{t-1} < a_t < d_s$ for $t=1,2\ldots,k-1$ and $s=1,2,\ldots,m$ so $a_{t} \, a_{t-1} \, d_s \longmapsto a_{t} \, d_s \, a_{t-1}.$
        \item $a_t < x < d_s$ for $t=1,\ldots,k$ and $s=1,\ldots,m$ so $x\, a_t \, d_s\longmapsto x\, d_s \, a_t.$
        \item $x < c_1 \leq d_s$ for $s=1,\ldots,m$ so 
        $c_1 \, x \, d_ s \longmapsto c_1 \, d_s \, x.$
        \item $x<d_{t}<d_{t+1}$ for $1\leq t<m$ so
        $x d_t d_{t-1}\longmapsto d_t x d_{t-1}.$
    \end{itemize}
    By repeatedly applying the relations above, we obtain that $\rev(\rw(B))$ is equivalent to the word $u_1=\textbf{c}\textbf{d}x\textbf{a}\textbf{b}.$
    Now observe that \\
    \begin{itemize}
        \item $a_{q-1} < a_q \leq b_p$ for all $q=1,2,\ldots,k$ and $p=i_q,\ldots,\ell$, so $a_q \, a_{q-1} \, b_p \longmapsto a_q \, b_p \, a_{q-1}.$
        \item $a_{q} < b_{p-1} \leq b_p$ for all $p=2,\ldots,m$ and $q=1,\ldots,r$ where $i_r$ is the closest index to $p-1$; therefore $a_q \, b_{p} \, b_{p-1} \longmapsto b_p \, a_q \, b_{p-1}.$
        \item $a_q \leq b_\ell < x$ so $x \, a_q \, b_\ell \longmapsto a_q \, x \, b_\ell.$
    \end{itemize}
    By repeatedly applying the previous relations and $d_1\,x\,a_k \longmapsto d_1\, a_k\, x$ we obtain that $u_1$ is equivalent to $u_2 = \textbf{c}\textbf{d}\textbf{a}x\textbf{b}$. Finally, we have
    \begin{itemize}
        \item $c_t < d_{s-1} < d_s$ for $s=2,3,\ldots,m$ and $t=1,2,\ldots,j_{s-1}$; so in this case $c_t \, d_s \, d_{s-1} \longmapsto d_s \, c_t \, d_{s-1}.$
    \end{itemize}
    and using this relation, together with the second and fourth relations listed in the first set of transformations we obtain that $u_2$ is equivalent to $\textbf{c}\textbf{a}\textbf{d}x\textbf{b}$ as desired. 
\end{proof}

\begin{lemma}\label{lemma:collapse-equivalent-words}
    Let $w, w'\in\mathcal{W}$. If $\rev(w)$ and $\rev(w')$ are Knuth-equivalent, then $\rho(w) = \rho(w')$.
\end{lemma}

\begin{proof}
   Suppose
   \[
   w = v\,  c\,  a\,  b\qquad \text{and}\qquad w' = v \, a \, c\,  b
   \]
   for some $v\in\mathcal{W}$ with $a < b\leq c$ so that the reversed words are Knuth equivalent. The multiline queues $M_w$ and $M_{w'}$ only differ in the top 3 rows. Using top-to-bottom collapsing from \cref{prop:top-to-bot-collapse}, both outputs coincide after the $e_{[2,k-1]}^\star e_{[1,k-1]}^\star$ is applied. Thus, $M_w$ and $M_{w'}$ collapse to the same nonwrapping multiline queue $\rho_N(M_w)=\rho_N(M_{w'})$, and hence $\rho(w) = \rho(w')$. 
   If $w=v\,c\,a\,b\,v'$ and $w'=v\,a\,c\,b\,v'$ for some $v,v'\in\mathcal{W}$, then $\rho_N(M_w)=\rho_N(M_{w'})$, given by inserting $v'$ into $\rho_N(M_{v\,c\,a\,b})=\rho_N(M_{v\,a\,c\,b})$. The other case of Knuth-equivalency is analogous.
\end{proof}

\begin{proof}[Proof of \cref{thm:bijectionSSYT-MLQ0}]

Let $T\in\SSYT(\lambda)$. Recall that $\revcrw(T)$ is a concatenation of increasing subwords obtained from reading each column of $T$ from bottom to top. Our proof is by induction on the columns of $T$. Let $M_k$ be the output of the MLQ insertion of the last $k$ columns of $T$, and assume it has type $\lambda^{(k)}\coloneqq (\lambda'_{n-k+1},\ldots,\lambda'_{n-1},\lambda'_n)'$ where $n\coloneqq \lambda_1$  (the partition corresponding to the last $k$ columns of $\lambda$). In other words, row $j$ of $M_k$ has $\lambda'_{n-k+j}$ balls for $1\leq j\leq k$.
Then $M_{k+1}$ is obtained by inserting the $(n-k)$'th
column of $T$. 
This column has length $\lambda'_{n-k}$. By the row-semi-strict condition on $T$, no bumping of balls in $M_k$ occurs when the entries of this column are inserted into $M_k$. Thus the top $k$ rows of $M_{k+1}$ have the same shape as $M_k$, and its first row has $\lambda'_{n-k}$ balls, making the shape of $M_{k+1}$ equal to $\lambda^{(k+1)}$ as desired. Completing the inductive argument, we conclude that $\rho(\revcrw(T)) \in \MLQ_0(\lambda').$

By \cref{lemma:knuth-equivalence-collapsing}, $\rev(\sigma(\colw(T)))$ and $\rev(\colw(T))$ are Knuth equivalent, so \cref{lemma:collapse-equivalent-words} and \cref{lemma:property-rho-and-Icol} together imply
\[
\Icol(\rw(\rho(\colw(T)))) = \Icol(\colw(T)) = T\] 
Thus, $(\Icol \circ \rw) \circ (\rho \circ \colw)$ is the identity map in $\SSYT(\lambda)$.

A similar argument for $M\in\MLQ_0(\lambda)$ shows that $\Icol(\rw(M))\in\SSYT(\lambda')$, since $\rw(M)$ is a concatenation of increasing subwords coming from the rows of $M$. Indeed, the rows are increasing entry-wise due to the nonwrapping condition on $M$. We have that $\Icol(\colw(\Icol(\rw(M)))) = \Icol(\rw(M))$ by \cref{lemma:property-rho-and-Icol}. 
Then $\rev(w')$ and $\rev(\rw(M))$ are Knuth equivalent words, so \cref{lemma:collapse-equivalent-words} implies 
\[
\rho(\colw(\Icol(\rw(M)))) = \rho(\rw(M)) = M.
\]
This shows $(
\rho \circ \colw) \circ (\Icol \circ \rw)$ is the identity map in $\MLQ_0(\lambda)$.
\end{proof}

\subsection{Collapsing on bosonic multiline queues}

We define collapsing on bosonic multiline queues analogously to collapsing on multiline queues by modifying the definition of parentheses matching to account for the "strictly left" pairing in bosonic multiline queues.

\begin{defn}\label{def:MLD operators}
    Let $D = (D_1,\ldots,D_L)$ be a bosonic multiline queue. The \emph{dropping operator} $\widetilde e_i$ acts on $D$ by moving the largest entry unmatched above in $\Par_i(\widetilde{\cw}(D))$ from $D_{i+1}$ to $D_i$. The \emph{lifting operator} $\widetilde f_i$ acts by moving the smallest entry unmatched below in $\Par_i(\widetilde{\cw}(D))$ from $D_{i}$ to $D_{i+1}$. 
\end{defn}

We get an analog to \cref{prop:rais/low-and-drop/lift}:
\begin{align*}
    \widetilde{cw}(\widetilde e_i(D)) = E_i(\widetilde{\cw}(D))\qquad \text{and}\qquad \widetilde{\cw}(\widetilde f_i(M)) = F_i(\widetilde{\cw}(M)).
\end{align*}

%
The above implies the $\widetilde e_i$'s and the $\widetilde f_i$'s are mutual inverses when they act non-trivially.
\begin{lemma}\label{lem:tilde ei inverse}
For a bosonic multiline queue $D$, if $\widetilde e_i(D)\neq D$, then $\widetilde f_i(\widetilde e_i(D))=D$. If $\widetilde f_i(D)\neq D$, then $\widetilde e_i(\widetilde f_i(D))=D$.
\end{lemma}

\begin{definition}\label{def:collapsing_mlds}
    Let $D$ be a bosonic multiline queue. Define $\tilde{\rho}(D)=(\tilde{\rho}_N(D),\tilde{\rho}_Q(D))$ 
    by applying the collapsing procedure $\rho$ in \cref{def:collapsing_mlqs} to $D$, with the notion of ``(un)matched above'' coming from \cref{def:matched above below MLD word} and the operators from \cref{def:MLD operators}. See \cref{ex:collapse MLD}.
\end{definition}

With the previous definitions we get an analog of \cref{thm:wrapping-to-nonwrapping} for bosonic multiline queues.

\begin{theorem}\label{thm:wrapping-to-nonwrapping_mlds}
    Let $n$ be an integer and $\lambda$ a partition. The collapsing procedure from \cref{def:collapsing_mlds} gives a weight-preserving bijection  
    \[
    \tilde \rho \, : \, \MLD(\lambda,n) \longrightarrow \bigcup_{\mu\leq \lambda} \MLD_0(\mu,n) \times \SSYT(\mu',\lambda'),
    \]
    with the following weights:
    \begin{itemize}
        \item for $D\in \MLD(\lambda,n)$, $\wt(D) = q^{\widetilde\maj(D)}x^D$ and
        \item for $(\widetilde N,\widetilde Q) \in \MLD_0(\mu,n)\times \SSYT(\mu',\lambda')$, $\wt(\widetilde N,\widetilde Q) = q^{\cocharge(\widetilde Q)}x^{\widetilde N}.$
    \end{itemize}
\end{theorem}


The following result follows from the same proof as \cref{thm:charge} with $\widetilde{\cw}$ replacing $\cw$. 

\begin{theorem}\label{thm:mld_charge}
    Let $D$ be a bosonic multiline queue with $(\widetilde N, \widetilde Q)=\widetilde{\rho}(D)$. 
    Then 
    \[ \widetilde\maj(D)=\cocharge(\widetilde Q)\;.
    \]
\end{theorem}

\begin{example}\label{ex:collapse MLD}
For the bosonic multiline queue $D\in\MLD((4,4,4))$ below, we show $\rho(D)=(\widetilde{N},\widetilde{Q})$ with $\widetilde{N}\in\MLD_0(\mu)$ and $\widetilde{Q}\in\SSYT(\mu)$ with $\mu=(4,3,2,1,1)$. 
\begin{figure}[h!]
    \centering
    \def \scale{0.3}
    \begin{subfigure}[t]{0.25\textwidth}
    \centering
    \includegraphics[scale = \scale]{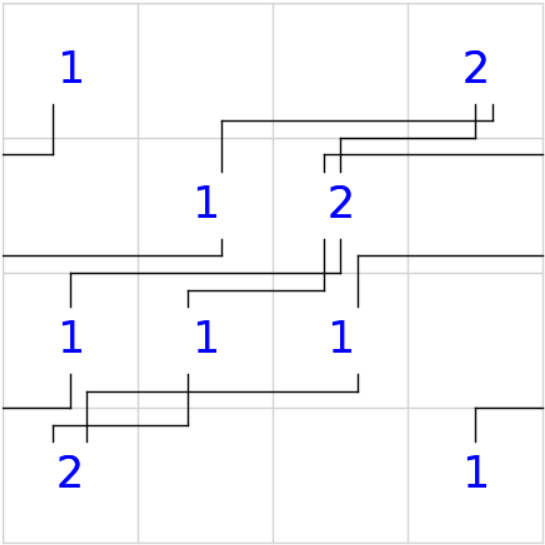}
    \caption*{{\Large $\substack{\lambda \;=\; (4,4,4) \\ \\ \widetilde\maj(D) = 12}$}}
    \end{subfigure}
    \raisebox{45pt}{\Large {$\longrightarrow$}}
    \raisebox{-7pt}{\begin{subfigure}[t]{0.05\textwidth}
        \begin{tikzpicture}
            \draw (8.75,0) arc (-30:30:4);
            \draw (0,0) arc (210:150:4);
        \end{tikzpicture}        
    \end{subfigure}}
    \begin{subfigure}[t]{0.22\textwidth}
    \centering
    \includegraphics[scale = \scale]{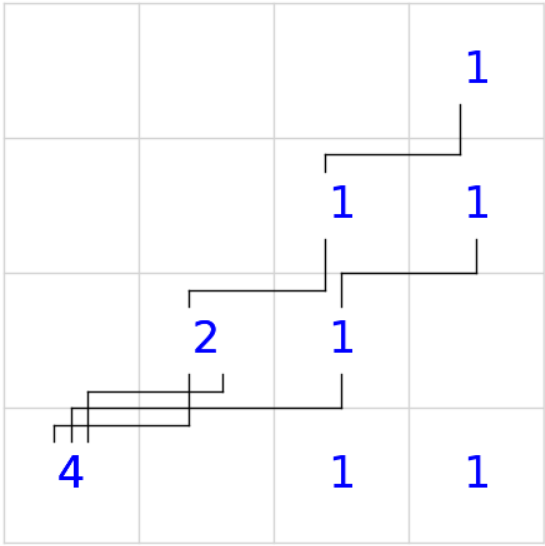}
    \caption*{{\Large $\substack{\mu \,=\,  (4,3,2,1,1)\\ \\ \widetilde\maj(\widetilde N) = 0}$}}
    \end{subfigure}
    $\quad\quad$
    \raisebox{45pt}{\huge $,$}
    $\;$
    \begin{subfigure}[t]{0.2\textwidth}
    \centering
    \raisebox{65pt}{$\tableau{4\\ 3&4\\2&2&3\\1&1&1&2&3&4}$}
    \caption*{{\Large $\substack{\cocharge(\widetilde Q) = 12 \\ }$}}
    \end{subfigure}
\end{figure}
Below, we show the collapsing procedure of $M$ row by row, as well as the construction of $\widetilde{\rho}_Q(D)$.

\begin{center}
    \begin{tikzpicture}[scale=0.5]
        \def \w{1};
        \def \h{1};
        \def \r{0.25};
        \node at (-1,5) {$r$};
        \node at (-1,2) {$\widetilde N_r$};
        \node at (-1,-1.75) {$Q_r$};
        \foreach \i in {0,...,4}
        {
        \node at (2+6*\i,5) {\i};
        }
        \begin{scope}[xshift=0cm]
        \foreach \i in {0,...,4}
        {
        \draw[gray!50] (0,\i*\h)--(\w*4,\i*\h);
        }
        \foreach \i in {0,...,4}
        {
        \draw[gray!50] (\w*\i,0)--(\w*\i,4*\h);
        }
        \foreach \xx\yy\rr in {0/0/2,3/0/1,0/1/1,1/1/1,2/1/1,1/2/1,2/2/2,0/3/1,3/3/2}
        {
            \node at (\w*.5+\w*\xx,\h*.5+\h*\yy) {\tcr{\rr}};
        }
        \node at (2,-1.75) {\scalebox{0.8}{\tableau{\emptyset\\}}};
        \end{scope}
    
        
        \begin{scope}[xshift=6cm]
        \foreach \i in {0,...,4}
        {
        \draw[gray!50] (0,\i*\h)--(\w*4,\i*\h);
        }
        \foreach \i in {0,...,4}
        {
        \draw[gray!50] (\w*\i,0)--(\w*\i,4*\h);
        }
        \foreach \xx\yy\rr in {0/1/1,1/1/1,2/1/1,1/2/1,2/2/2,0/3/1,3/3/2}
        {
            \node at (\w*.5+\w*\xx,\h*.5+\h*\yy) {\tcr{\rr}};
        }
        \foreach \xx\yy\rr in {0/0/2,3/0/1}
        {
            \node at (\w*.5+\w*\xx,\h*.5+\h*\yy) {\tcb{\rr}};
        }
        
        \node at (2,-1.75) {\scalebox{0.8}{\tableau{1&1&1\\}}};
        \end{scope}
    
        
        \begin{scope}[xshift=12cm]
        \foreach \i in {0,...,4}
        {
        \draw[gray!50] (0,\i*\h)--(\w*4,\i*\h);
        }
        \foreach \i in {0,...,4}
        {
        \draw[gray!50] (\w*\i,0)--(\w*\i,4*\h);
        }
        \foreach \xx\yy\rr in {1/2/1,2/2/2,0/3/1,3/3/2}
        {
            \node at (\w*.5+\w*\xx,\h*.5+\h*\yy) {\tcr{\rr}};
        }
        \foreach \xx\yy\rr in {0/0/3,3/0/1,1/1/1,2/1/1}
        {
            \node at (\w*.5+\w*\xx,\h*.5+\h*\yy) {\tcb{\rr}};
        }
        \node at (2,-1.75) {\scalebox{0.8}{\tableau{2&2\\1&1&1&2}}};
        \end{scope}
    
        
        \begin{scope}[xshift=18cm]
        \foreach \i in {0,...,4}
        {
        \draw[gray!50] (0,\i*\h)--(\w*4,\i*\h);
        }
        \foreach \i in {0,...,4}
        {
        \draw[gray!50] (\w*\i,0)--(\w*\i,4*\h);
        }
        \foreach \xx\yy\rr in {0/3/1,3/3/2}
        {
            \node at (\w*.5+\w*\xx,\h*.5+\h*\yy) {\tcr{\rr}};
        }
        \foreach \xx\yy\rr in {0/0/3,3/0/1,1/1/2,2/1/1,2/0/1,2/2/1}
        {
            \node at (\w*.5+\w*\xx,\h*.5+\h*\yy) {\tcb{\rr}};
        }
        
        \node at (2,-1.75) {\scalebox{0.8}{\tableau{3\\2&2&3\\1&1&1&2&3}}};
        \end{scope}
    
        
        \begin{scope}[xshift=24cm]
        \foreach \i in {0,...,4}
        {
        \draw[gray!50] (0,\i*\h)--(\w*4,\i*\h);
        }
        \foreach \i in {0,...,4}
        {
        \draw[gray!50] (\w*\i,0)--(\w*\i,4*\h);
        }
        \foreach \xx\yy\rr in {0/0/4,3/0/1,1/1/2,2/1/1,2/0/1,2/2/1,3/2/1,3/3/1}
        {
            \node at (\w*.5+\w*\xx,\h*.5+\h*\yy) {\tcb{\rr}};
        }
        
        \node at (2,-1.75) {\scalebox{0.8}{\tableau{4\\ 3&4\\2&2&3\\1&1&1&2&3&4}}};
        \end{scope}
    \end{tikzpicture}
\end{center}

\end{example} 

\subsection{Formula for $P_{\lambda}(X;q,0)$ 
via twisted multiline queues
}\label{sec:GMLQ maj}

In this section, we define the statistic $\maj_G$ as the extension of $\maj$ to twisted multiline queues, by thinking of particles as pairing to the right and anti-particles as pairing to the left; every pairing wrapping to the right contributes a positive term to the $\maj$, and every pairing wrapping to the left contributes a negative term. We show in \cref{prop:local_majG_energy} that this definition coincides with the \emph{energy function} of \cite{NY95} We then generalize results of \cite{AGS20} to obtain a family of formulas for the $q$-Whittaker polynomials as a sum over twisted multiline queues with the $\maj_G$ statistic. 

\begin{defn}\label{def:majg}
.Let $M
\in\GMLQ(\alpha,n)$ with an associated labeling $L(M)$. For $1\leq r,\ell \leq L$, let $m_{r,\ell}$ (\emph{resp.}~$a_{r,\ell}$) be the number of particles (\emph{resp.}~anti-particles) of type $\ell$ that wrap when pairing to the right (\emph{resp.}~left) from row $r$ to row $r-1$, as shown in \cref{ex:GMLQ}. Define
\begin{equation}\label{eq:majG}
\maj_G(M)=\sum_{1\leq r,\ell\leq L} m_{r,\ell}(\ell-r+1)-a_{r,\ell}(\ell-r+1).
\end{equation}
\end{defn}

\begin{lemma}\label{lem:basecase}
    When $M\in\MLQ(\lambda,n)$, $\maj_G(M)=\maj(M)$.
\end{lemma}  

\begin{proof}
    The particle phase of \cref{def:gmlq_particlewise} is identical to the FM algorithm. Thus when $M$ is a straight multiline queue, the labeling of the particles is identical to that in \cref{def:FM2}, and so $L_G(M)^+=L(M)$. Furthermore, for $1\leq r\leq \lambda_1$, one may check that all anti-particles in row $r$ are labeled $r-1$ in $L_G(M)$. Thus the contribution to $\maj$ from anti-particles wrapping during either phase of the generalized pairing procedure is $(r-1)-r+1=0$. Therefore, only wrapping particles contribute to $\maj_G(M)$, and the contribution is identical to that for $\maj(M)$.
\end{proof}

\begin{example}\label{ex:GMLQ} For $\alpha=(2,2,3)$, the labeled particles (circles) and anti-particles (squares) are shown for $M=(\{2,3\},\{1,4\},\{2,3,4\})\in\GMLQ(\alpha,4)$, \linebreak[0]$\sigma_2(M)=(\{2,3\},\{1,2,4\},\{3,4\})\in\GMLQ(s_2\cdot\alpha,4)$ and $\sigma_1(\sigma_2(M))=(\{2,3,4\},\{1,2\},\{3,4\})\in\GMLQ(s_1s_2\cdot\alpha,4)$. The positive (blue) and negative (red) contributions to $\maj_G$ are shown, totalling $\maj_G=2$ in each case.
\end{example}
\vspace{-0.1in}
\begin{center}
   
    \resizebox{\linewidth}{!}{
    \begin{tikzpicture}[scale=0.7]
    \def \w{1};
    \def \h{1};
    \def \r{0.2};

    \def \r{0.3};
    \def \off{0.09};
    \begin{scope}[xshift=9cm]
    \foreach \i in {0,...,3}
    {
    \draw[gray!50] (0,\i*\h)--(\w*4,\i*\h);
    }
    \foreach \i in {0,...,4}
    {
    \draw[gray!50] (\w*\i,0)--(\w*\i,3*\h);
    }
    \foreach \xx\yy\i in {1/0,2/0,0/1,3/1,1/2,2/2,3/2}
    {
    \draw[blue] (\w*.5+\w*\xx,\h*.5+\h*\yy) circle (\r cm);
    }
    \foreach \xx\yy\i in {0/0,3/0,1/1,2/1,0/2}
    {
    \draw[red] (\w*.5+\w*\xx-\w*\r,\h*.5+\h*\yy+\h*\r) -- (\w*.5+\w*\xx+\w*\r,\h*.5+\h*\yy+\h*\r) -- (\w*.5+\w*\xx+\w*\r,\h*.5+\h*\yy-\h*\r)--(\w*.5+\w*\xx-\w*\r,\h*.5+\h*\yy-\h*\r)--(\w*.5+\w*\xx-\w*\r,\h*.5+\h*\yy+\h*\r);
    }
    \node at (-1,2*\h) {\scriptsize\tcr{2-3+1}};
    \node at (-1,\h) {\scriptsize\tcr{2-2+1}};
    \node at (1+4*\w,2*\h) {\scriptsize\tcb{3-3+1}};
    \node at (1+4*\w,\h) {\scriptsize\tcb{3-2+1}};

    \draw[blue] (\w*2.5-\off,\h*2.5-\r)--(\w*2.5-\off,\h*1.9)--(\w*3.5-\off,\h*1.9)--(\w*3.5-\off,\h*1.5+\r);

    \draw[blue,-stealth] (\w*3.5-\off,\h*2.5-\r)--(\w*3.5-\off,\h*2.1)--(\w*4.2,\h*2.1);
    \draw[blue] (-.2,\h*2.1)--(\w*0.5-\off,\h*2.1)--(\w*0.5-\off,\h*1.5+\r);%

    \draw[blue] (\w*0.5-\off,\h*1.5-\r)--(\w*0.5-\off,\h*0.9)--(\w*1.5-\off,\h*0.9)--(\w*1.5-\off,\h*0.5+\r);

    \draw[blue,-stealth] (\w*3.5-\off,\h*1.5-\r)--(\w*3.5-\off,\h*1.1)--(\w*4.2,\h*1.1);
    \draw[blue] (-.2,\h*1.1)--(\w*2.5-\off,\h*1.1)--(\w*2.5-\off,\h*0.5+\r);%


    \draw[red] (\w*1.5-\off,\h*2.5-\r)--(\w*1.5-\off,\h*1.5+\r);%
    \draw[red,-stealth] (\w*0.5+\off,\h*2.5-\r)--(\w*0.5+\off,\h*2)--(-0.2,\h*2.0);
    \draw[red] (4.2,\h*2.0)--(\w*2.5+\off,\h*2.0)--(\w*2.5+\off,\h*1.5+\r);

    \draw[red] (\w*1.5+\off,\h*1.5-\r)--(\w*1.5+\off,\h*0.95)--(\w*0.5+\off,\h*0.95)--(\w*0.5+\off,\h*0.5+\r);

    \draw[red,-stealth] (\w*2.5+\off,\h*1.5-\r)--(\w*2.5+\off,\h*1.025)--(-0.2,\h*1.025);
    \draw[red] (4.2,\h*1.025)--(\w*3.5+\off,\h*1.025)--(\w*3.5+\off,\h*0.5+\r);


    \node[red] at (0.5,0.5) {0};
    \node[blue] at (1.5,0.5) {3};
    \node[blue] at (2.5,0.5) {3};
    \node[red] at (3.5,0.5) {1};

    \node[blue] at (0.5,1.5) {3};
    \node[red] at (1.5,1.5) {2};
    \node[red] at (2.5,1.5) {1};
    \node[blue] at (3.5,1.5) {3};

    \node[red] at (0.5,2.5) {2};
    \node[blue] at (1.5,2.5) {3};
    \node[blue] at (2.5,2.5) {3};
    \node[blue] at (3.5,2.5) {3};

    \draw[black,<->] (6,1.5) -- (7,1.5);
    \node[black] at (6.5,1.85) {$\sigma_2$};
    \end{scope}

    \begin{scope}[xshift=18cm]
    \foreach \i in {0,...,3}
    {
    \draw[gray!50] (0,\i*\h)--(\w*4,\i*\h);
    }
    \foreach \i in {0,...,4}
    {
    \draw[gray!50] (\w*\i,0)--(\w*\i,3*\h);
    }
    \foreach \xx\yy\i in {1/0,2/0,0/1,1/1,3/1,2/2,3/2}
    {
    \draw[blue] (\w*.5+\w*\xx,\h*.5+\h*\yy) circle (\r cm);
    }
    \foreach \xx\yy\i in {0/0,3/0,2/1,1/2,0/2}
    {
    \draw[red] (\w*.5+\w*\xx-\w*\r,\h*.5+\h*\yy+\h*\r) -- (\w*.5+\w*\xx+\w*\r,\h*.5+\h*\yy+\h*\r) -- (\w*.5+\w*\xx+\w*\r,\h*.5+\h*\yy-\h*\r)--(\w*.5+\w*\xx-\w*\r,\h*.5+\h*\yy-\h*\r)--(\w*.5+\w*\xx-\w*\r,\h*.5+\h*\yy+\h*\r);
    }

    \draw[blue] (\w*3.5-\off,\h*2.5-\r)--(\w*3.5-\off,\h*1.5+\r);
    \draw[blue] (\w*1.5-\off,\h*2.5-\r)--(\w*1.5-\off,\h*1.5+\r);%
    \draw[blue,-stealth] (\w*2.5-\off,\h*2.5-\r)--(\w*2.5-\off,\h*2.1)--(\w*4.2,\h*2.1);
    \draw[blue] (-.2,\h*2.1)--(\w*0.5-\off,\h*2.1)--(\w*0.5-\off,\h*1.5+\r);%

    \draw[blue] (\w*0.5-\off,\h*1.5-\r)--(\w*0.5-\off,\h*0.9)--(\w*1.5-\off,\h*0.9)--(\w*1.5-\off,\h*0.5+\r);

    \draw[blue,-stealth] (\w*3.5-\off,\h*1.5-\r)--(\w*3.5-\off,\h*1.1)--(\w*4.2,\h*1.1);
    \draw[blue] (-.2,\h*1.1)--(\w*2.5-\off,\h*1.1)--(\w*2.5-\off,\h*0.5+\r);%


    
    \draw[red,-stealth] (\w*0.5+\off,\h*2.5-\r)--(\w*0.5+\off,\h*2)--(-0.2,\h*2.0);
    \draw[red] (4.2,\h*2.0)--(\w*2.5+\off,\h*2.0)--(\w*2.5+\off,\h*1.5+\r);

    \draw[red] (\w*1.5+\off,\h*1.5-\r)--(\w*1.5+\off,\h*0.95)--(\w*0.5+\off,\h*0.95)--(\w*0.5+\off,\h*0.5+\r);

    \draw[red,-stealth] (\w*2.5+\off,\h*1.5-\r)--(\w*2.5+\off,\h*1.025)--(-0.2,\h*1.025);
    \draw[red] (4.2,\h*1.025)--(\w*3.5+\off,\h*1.025)--(\w*3.5+\off,\h*0.5+\r);


    \node[red] at (0.5,0.5) {0};
    \node[blue] at (1.5,0.5) {3};
    \node[blue] at (2.5,0.5) {3};
    \node[red] at (3.5,0.5) {1};

    \node[blue] at (0.5,1.5) {3};
    \node[blue] at (1.5,1.5) {2};
    \node[red] at (2.5,1.5) {1};
    \node[blue] at (3.5,1.5) {3};

    \node[red] at (0.5,2.5) {2};
    \node[red] at (1.5,2.5) {2};
    \node[blue] at (2.5,2.5) {3};
    \node[blue] at (3.5,2.5) {3};

    \node at (-1,2*\h) {\scriptsize\tcr{2-3+1}};
    \node at (-1,\h) {\scriptsize\tcr{2-2+1}};
    \node at (1+4*\w,2*\h) {\scriptsize\tcb{3-3+1}};
    \node at (1+4*\w,\h) {\scriptsize\tcb{3-2+1}};

    \draw[black,<->] (6,1.5) -- (7,1.5);
    \node[black] at (6.5,1.85) {$\sigma_1$};
    
    \end{scope}

    \begin{scope}[xshift=27cm]
    \foreach \i in {0,...,3}
    {
    \draw[gray!50] (0,\i*\h)--(\w*4,\i*\h);
    }
    \foreach \i in {0,...,4}
    {
    \draw[gray!50] (\w*\i,0)--(\w*\i,3*\h);
    }
    \foreach \xx\yy\i in {1/0,2/0,3/0,0/1,1/1,2/2,3/2}
    {
    \draw[blue] (\w*.5+\w*\xx,\h*.5+\h*\yy) circle (\r cm);
    }
    \foreach \xx\yy\i in {0/0,2/1,3/1,1/2,0/2}
    {
    \draw[red] (\w*.5+\w*\xx-\w*\r,\h*.5+\h*\yy+\h*\r) -- (\w*.5+\w*\xx+\w*\r,\h*.5+\h*\yy+\h*\r) -- (\w*.5+\w*\xx+\w*\r,\h*.5+\h*\yy-\h*\r)--(\w*.5+\w*\xx-\w*\r,\h*.5+\h*\yy-\h*\r)--(\w*.5+\w*\xx-\w*\r,\h*.5+\h*\yy+\h*\r);
    }

    \draw[blue,-stealth] (\w*2.5-\off,\h*2.5-\r)--(\w*2.5-\off,\h*2.125)--(\w*4.2,\h*2.125);
    \draw[blue] (-.2,\h*2.125)--(\w*0.5-\off,\h*2.125)--(\w*0.5-\off,\h*1.5+\r);%

    \draw[blue,-stealth] (\w*3.5-\off,\h*2.5-\r)--(\w*3.5-\off,\h*1.875)--(\w*4.2,\h*1.875);
    \draw[blue] (-.2,\h*1.875)--(\w*1.5-\off,\h*1.875)--(\w*1.5-\off,\h*1.5+\r);%

    \draw[blue] (\w*0.5-\off,\h*1.5-\r)--(\w*0.5-\off,\h*0.9)--(\w*1.5-\off,\h*0.9)--(\w*1.5-\off,\h*0.5+\r);

    \draw[blue] (\w*1.5-\off,\h*1.5-\r)--(\w*1.5-\off,\h*1.1)--(\w*2.5-\off,\h*1.1)--(\w*2.5-\off,\h*0.5+\r);

    \draw[blue,] (\w*3.5-\off,\h*1.5-\r)--(\w*3.5-\off,\h*0.5+\r);%


    \draw[red,-stealth] (\w*0.5+\off,\h*2.5-\r)--(\w*0.5+\off,\h*2.05)--(-0.2,\h*2.05);
    \draw[red] (4.2,\h*2.05)--(\w*2.5+\off,\h*2.05)--(\w*2.5+\off,\h*1.5+\r);

    \draw[red,-stealth] (\w*1.5+\off,\h*2.5-\r)--(\w*1.5+\off,\h*1.95)--(-0.2,\h*1.95);
    \draw[red] (4.2,\h*1.95)--(\w*3.5+\off,\h*1.95)--(\w*3.5+\off,\h*1.5+\r);

    \draw[red] (\w*2.5+\off,\h*1.5-\r)--(\w*2.5+\off,\h*1.0)--(\w*0.5+\off,\h*1.0)--(\w*0.5+\off,\h*0.5+\r);


    \node[red] at (0.5,0.5) {0};
    \node[blue] at (1.5,0.5) {3};
    \node[blue] at (2.5,0.5) {3};
    \node[blue] at (3.5,0.5) {1};

    \node[blue] at (0.5,1.5) {3};
    \node[blue] at (1.5,1.5) {3};
    \node[red] at (2.5,1.5) {1};
    \node[red] at (3.5,1.5) {1};

    \node[red] at (0.5,2.5) {2};
    \node[red] at (1.5,2.5) {2};
    \node[blue] at (2.5,2.5) {3};
    \node[blue] at (3.5,2.5) {3};

    \node at (-1,2*\h+3*\off) {\scriptsize\tcr{2-3+1}};
    \node at (-1,2*\h-3*\off) {\scriptsize\tcr{2-3+1}};
    \node at (1+4*\w,2*\h+3*\off) {\scriptsize\tcb{3-3+1}};
    \node at (1+4*\w,2*\h-3*\off) {\scriptsize\tcb{3-3+1}};
    
    \end{scope}
    
    \end{tikzpicture}
}
\end{center}

The main result of \cite{AGS20} states that the involution $\sigma_i$ preserves the labeling of the bottom row of a twisted multiline queue, thus preserving the distribution of states of the TASEP onto which the twisted multiline queues project, and thereby giving a twisted analogue of \eqref{eq:FM result}. We claim that $\maj_G$ is also preserved under $\sigma_i$. We will show that our definition of $\maj_G$ can be reformulated in terms of an \emph{energy function} on tensors of KR crystals that was introduced in   
\cite{NY95}.

\begin{defn}\label{def:local_maj}
    Let $w$ be a word in $\{0,1,\ldots,L\}$ in $n$ letters and let $D \subseteq [n]$ represent a row of a twisted multiline queue. Define 
    \begin{equation}\label{eq:localmaj}
    \maj_G(w;D) = \sum_{a} \left( \ell(a)-i\right ) - \sum_{b} \left (\ell(b)-i\right)
    \end{equation}
    where the first sum is over all pairings $a$ wrapping (to the right) during the particle phase and the second is over all pairings $b$ wrapping during the anti-particle phase at that row.
\end{defn}
In a twisted multiline queue $M=(B_1,\ldots,B_L)$, if we let $w^{(i)}$ be the labeling word of row $i$ in $L_G(M)$, then $\maj_G(w^{(i+1)};B_i)$ is the contribution to $\maj_G(M)$ of the pairings from row $i+1$ to row $i$, and thus we may write 
\[
\maj_G(M)=\sum_{i=1}^{L-1}\maj_G(w^{(i+1)};B_i).
\]

\begin{defn}\label{def:indicator_words}
    Let $w$ be a word in $\{0,1,\ldots,L\}$. The \emph{nested indicator decomposition} of $w$ is the set of words $\{w_1,w_2,\ldots,w_L\}$ given by $w_j = w_{j1}\,w_{j2}\,\ldots\,w_{jn}\in\{0,1\}^n$ where $w_{jk} = 1$ if and only if the $k$-th entry of $w$ is greater than or equal to $j$, so that $w=w_1+\cdots+w_L$ if addition is performed componentwise. For an indicator word $u\in\{0,1\}^n$, define $\iota(u)$ to be the subset of $[n]$ associated to $u$: $\iota(u)=\{i\in[n]:u_i=1\}$. In particular, $w_1,\ldots,w_L$ are nested indicator words if and only if $\iota(w_L)\subseteq\cdots\subseteq \iota(w_2)\subseteq\iota(w_1)$.
\end{defn}

\begin{example}
    For $w = 2\,5\,2\,3\,4\,2$ the nested indicator decomposition is:
    \[w_1 =w_2= 1\,1\,1\,1\,1\,1,\ w_3 = 0\,1\,0\,1\,1\,0,\ w_4 = 0\,1\,0\,0\,1\,0,\ \mbox{and}\ w_5 = 0\,1\,0\,0\,0\,0.\]
The corresponding subsets are
\[\iota(w_1)=\iota(w_2)=\{1,2,3,4,5,6\},\ \iota(w_3)=\{2,4,5\},\ \iota(w_4)=\{2,5\},\ \mbox{and}\ \iota(w_5)=\{2\}.\]
\end{example}

\begin{defn}\label{def:energy_function}
Given a queue $D\subseteq[n]$ and a word $w\in\{0,1,\ldots,L\}^n$ with nested indicator decomposition $w=w_1+\cdots+w_L$, the \emph{energy function} $H(w_j;D)$ is defined as the number of wrapping pairings in $\Par^c(D,\iota(w_j))$. 
Then, define $H(w;D)\coloneqq \sum_{j=1}^L H(w_j;D)$. Finally, for a multiline queue $M=(B_1,\ldots,B_L)$ with $L_G(M)=(w^{(1)},\ldots,w^{(L)})$, define $H(M)\coloneqq \sum_{i=1}^{L-1} H(w^{(i+1)};B_i)$.
\end{defn}

\begin{remark}\label{rem:H}
\cref{def:energy_function} is a translation of the energy function as defined in \cite{NY95}. The authors only consider objects corresponding to straight multiline queues in their paper; however, the energy function is shown to be invariant under the combinatorial $R$-matrix (in our setup, these are the operators~$\sigma_i$).  Since the $\sigma_i$'s give crystal isomorphisms between twisted and straight configurations, the same energy function extends immediately to the twisted multiline queue setting. In particular, $H(M)=H(\sigma_i(M))$ for any (twisted) multiline queue $M$. 
\end{remark}

\begin{example}\label{ex:energy function}
    Consider the twisted multiline queue $M = (\{2,3\},\{1,4\},\{2,3,4\})$ from \cref{ex:GMLQ}. Then $$L(G) = (w^{(1)} = 0\,3\,3\,1\,,\,w^{(2)} = 3\,2\,1\,3\,,\, w^{(3)} = 2\,3\,3\,3)$$ and the computation of the energy function $H(M)=2$ can be decomposed as follows:
    \begin{table}[h]
        \centering
        \begin{tabular}{|c|c||c|c|}
        \hline
             $w^{(3)} = 2\,3\,3\,3$ & $B_2 = 1\,0\,0\,1$ & $w^{(2)} = 3\,2\,1\,3$ & $B_1 = 0\,1\,1\,0$ \\
             \hline
             \hline
             $w^{(3)}_1 = 1\,1\,1\,1$ & $H(w^{(3)}_1,B_2) = 0$ & $w^{(2)}_1 = 1\,1\,1\,1$ & $H(w^{(2)}_1,B_1) = 0$ \\
             $w^{(3)}_2 = 1\,1\,1\,1$ & $H(w^{(3)}_2,B_2) = 0$ & $w^{(2)}_2 = 1\,1\,0\,1$ & $H(w^{(2)}_2,B_1) = 0$ \\
             $w^{(3)}_3 = 0\,1\,1\,1$ & $H(w^{(3)}_3,B_2) = 1$ & $w^{(2)}_3 = 1\,0\,0\,1$ & $H(w^{(2)}_3,B_1) = 1$ \\
             \hline
        \end{tabular}
    \end{table}
\end{example}

\begin{defn}\label{def:reversed-words-and-set}
    For a word $w$, let $\rev(w)$ be the reversed word with letters $\rev(w)_i = w_{n-i+1}$. For a subset $S\subseteq [n]$ let $\rev(S) = \{n-i+1 \,\colon\, i\in S\}$ be the reversed version of $S$ in $[n]$. 
\end{defn}

Our main observation, which yields \cref{thm:GMLQ}, is that $\maj_G$ relates the particle and anti-particle energies, coming from the right-wrapping of particles and left-wrapping of anti-particles, respectively.

\begin{prop}\label{prop:local_majG_energy}
    In the setup of \cref{def:local_maj,def:indicator_words}, let $w=w_1+\cdots+w_L$ be the decomposition of a word $w$ into nested indicator vectors, and let $D\subseteq[n]$ be a queue. Then 
    \begin{equation}\label{eq:maj energy}
    \maj_G(w\,;D) = \sum_{j=1}^L H(w_j\,;D).
    \end{equation}
    In particular, this implies that for a twisted multiline queue $M$, we have $\maj_G(M)=H(M)$.
\end{prop}

\begin{proof}

For an indicator word 
$v\in \{0,1\}^n$ and a set of balls $D\subseteq [n]$, denote by $\overline{D}$ the complement of $D$ in $[n]$ and by $\overline{v}$ the complement of the word $v$ with letters $\overline{v}_i = 1-v_i$. Define an \emph{anti-particle energy function}, $H^{\leftarrow}(v\,;D)$ as the number of wrapping pairings in $\Pair^c(\rev(D),\rev(v))$. We interpret this energy function as counting the number of wrapping pairings \emph{to the left} in the two-row arrangement $(D,\iota(v))$. We claim that $H^{\leftarrow}(\overline{v}\,;\overline{D}) = H(v\,;D)$. This can be seen from the identity $H^{\leftarrow}(\overline{v}\,;\overline{D}) = H(\iota^{-1}(\overline{D})\,;\iota(\overline{w}))$ 
and the fact that a wrapping particle in $\Pair^c(D,\iota(v))$ unbalances the number of anti-particles to its right creating a wrapping in $\Pair^c(\iota(\overline{w}),\overline{D}).$

Let $w=w_1+w_2+\cdots+w_L$ be the word labeling row $i$ of $L_G(M)$, let $D$ be the queue at row $i-1$, and consider $\maj_G(w;D)$ (note that we have decremented the index for nicer computations). We will write 
$h^\rightarrow_j\coloneqq H(w_j;D)$ and $h^\leftarrow_j\coloneqq H^\leftarrow(\overline{w_j};\overline{D})$. Define $\ell=\min_j\{|\iota(w_{j+1})|<|D|\leq |\iota(w_{j})|\}$ to be the smallest label of a ball  that pairs during the particle phase in row $i$ of $L_G(M)$. Let us decorate the letters $\ell$ in $w$: label the $\ell$'s that pair during the particle phase by $\ell^+$ and those that pair during the anti-particle phase by $\ell^-$, and let those sites be indexed by the indicator vectors $u_{\ell^+}, u_{\ell^-}$ so that $w_\ell=w_{\ell+1}+u_{\ell^+}+u_{\ell^-}$. 
Now, by comparing definitions, we observe that 
\begin{align*}
m_{i,j}&=h^\rightarrow_j-h^\rightarrow_{j+1},&\mbox{if}\ \ell<j\leq L\qquad\qquad \mbox{and}\ m_{i,\ell}&=h^\rightarrow_{\ell^+}-h^\rightarrow_{\ell+1} \\
a_{i,j}&=h^\leftarrow_j-h^\leftarrow_{j+1},&\mbox{if}\ i-1\leq j< \ell\qquad\qquad \mbox{and}\ a_{i,\ell}&=h^\leftarrow_{\ell^+}-h^\leftarrow_{\ell}.
\end{align*}
with $h^\rightarrow_{L+1}\coloneqq 0$, where $m_{i,j}$ and $a_{i,j}$ are the quantities in \cref{def:majg}. Plugging these into \eqref{eq:localmaj}, we obtain
{\allowdisplaybreaks
\begin{align*}
    \maj_G(w;D)&=\sum_{j=\ell}^L m_{i,j}(j-i+1)-\sum_{j=i-1}^{\ell} a_{i,j}(j-i+1)\\
    &\begin{aligned}=m_{i,\ell}(\ell-i+1)&+\sum_{j=\ell+1}^L (h^\rightarrow_j-h^\rightarrow_{j+1})(j-i+1)\\&-a_{i,\ell}(\ell-i+1)-\sum_{j=i-1}^{\ell-1} (h^\leftarrow_j-h^\leftarrow_{j+1})(j-i+1)\end{aligned}\\
    &\begin{aligned}=(-h^\rightarrow_{\ell+1}+h^\leftarrow_{\ell})(\ell-i+1) & + h^\rightarrow_{\ell+1}(\ell+1-i+1) - h^\rightarrow_{L+1}(L-i+1)\\ &+\sum_{j=\ell+2}^L(h^\rightarrow_j(j-i+1)-h^\rightarrow_j((j-1)-i+1))\\&-\sum_{j=i}^{\ell-1}(h^\leftarrow_j((j-1)-i+1)-h^\leftarrow_j(j-i+1))\\&+h^\leftarrow_{i-1}((i-1)-i+1)-h^\leftarrow_{\ell}((\ell-1)-i+1)\end{aligned}\\
    &=h^\rightarrow_{\ell+1}+h^\leftarrow_{\ell}+\sum_{j=\ell+2}^L(h^\rightarrow_j)-\sum_{j=i}^{\ell-1}(-h^\leftarrow_j) = \sum_{j=i}^L h^\leftarrow_j
\end{align*}}

The last equality follows from $h^\leftarrow_j=h^\rightarrow_j$. It should be noted that if $w$ is the labeling word of row $i$ of $L_G(M)$, then it only contains the letters $\{i-1,\ldots,L\}$, which determines the indices of the sums. Since $H(w_j;D)=0$ for all $j<i$, we obtain the right hand side of \eqref{eq:maj energy}.
\end{proof}

\begin{corollary}\label{prop:charge commutes with sigma}
Let $\alpha$ be a composition with $\alpha^+=\lambda'$, $L\coloneqq \ell(\alpha)$,  $M\in\GMLQ(\alpha)$, and let $1\leq i\leq L-1$. Then $\maj_G(M)=\maj_G(\sigma_i(M))$. 
\end{corollary}

\begin{proof}
From \cref{prop:local_majG_energy}, we have that $\maj_G(M)=H(M)=H(\sigma_i (M))=\maj_G(\sigma_i(M))$, where the $\sigma_i$-invariance of $H$ is explained in \cref{rem:H}.
\end{proof}

\begin{prop}\label{lem: majg is charge cw}
    Let $\alpha$ be a composition and let $M\in\GMLQ(\alpha)$. Then 
    \[\maj_G(M)=\charge_G(\cw(M)).
    \]
\end{prop}

\begin{proof}
From \cref{lem:basecase} and \cref{thm:charge}, this is true for a straight multiline queue. Since the left hand side is invariant under the action of $\sigma_i$ by \cref{prop:charge commutes with sigma}, it is enough to show that is true for the right hand side, as well. By definition, $\charge_G$ is invariant under the action of $\LS_i$ on the column reading word, so we have $\charge_G(\cw(M))=\charge_G(\LS_i(\cw(M)))=\charge_G(\cw(\sigma_i(M)))$, proving our claim.
\end{proof}

By \cref{lem:coxeter,prop:charge commutes with sigma,thm:charge}, we obtain \cref{thm:GMLQ}. 
\begin{theorem}\label{thm:GMLQ}
    Let $\lambda$ be a partition, $n$ an integer, and let $\alpha$ be a composition with $\alpha^+=\lambda'$. Then
    \[
    P_{\lambda}(X;q,0)=\sum_{M\in\GMLQ(\alpha,n)} q^{\maj_G(M)}x^M.
    \]
\end{theorem}

\begin{remark} There is an analogous NY rule for bosonic multiline queues that naturally arises from the rank-level duality property of KR crystals of affine type (see, e.g. \cite{KMO-multiline-2016}). We remark that the same arguments yield analogous formulas for $\widetilde{H}_{\lambda}(X;q,0)$ in terms of a $\cocharge$ statistic defined on \emph{bosonic twisted multiline queues}, recently introduced in \cite{ManScrim24}; we omit these details in this article.
\end{remark}

\section{Multiline queue analogues of the RSK correspondence}\label{sec:mlq/d_RSK}

Recall that \cref{thm:wrapping-to-nonwrapping} is a bijection involving multiline queues and tableaux. In this section, we describe an analogue of the RSK correspondence in which all objects are multiline queues. These descriptions are equivalent to the double crystals considered in \cite{vanleeuwen2006double} due to the relation between collapsing and raising/lowering operators described in \cref{sec:doublecollapse}. However, in the context of multiline queues this correspondence becomes a useful tool to give an elementary proof of the $\charge$ formula in \cref{thm:charge} and to derive multiline queue expressions for $K_{\lambda\mu}(q,0)$. 

\subsection{Multiline queue RSK via commuting crystal operators}\label{sec:doublecollapse}
We define two operators $\rho^{\downarrow}$ and $\rho^{\leftarrow}$ acting on $\M$ and a $90^{\circ}$ rotation of it by treating the elements of this set as (twisted) multiline queues, and collapsing them both. This is equivalent to collapsing the multiline queue in two directions: downwards and to the left, respectively.

\begin{remark}\label{rem:directions}
    The interpretation of elements of $\M$ as multiline queues involves several choices, each leading to variations in the results.
    \begin{itemize}
        \item The pair of collapsing directions can be any orthogonal pair in the set $\{\uparrow\;,\;\downarrow\;,\;\leftarrow\;,\,\rightarrow\}$.
        \item Particles can be paired either weakly to the right or weakly to the left. 
    \end{itemize}
    We analyze the case of downward and leftward collapsing with particles pairing to their left, as downward collapsing aligns with previous sections, and leftward collapsing provides an elementary proof of  \cref{thm:charge}. Moreover, the pair (down, left) yields a symmetry as shown in \cref{prop:double-collapse}.
\end{remark}

\begin{defn}\label{def:rot}
For $B\in\M$, define $\rot(B)$ to be the rotation of $B$ by $90^{\circ}$ counterclockwise.  
We define two collapsing operations on multiline queues: \emph{downward collapsing} ($\rho^\downarrow(B)$) and \emph{leftward collapsing} ($\rho^\leftarrow(B)$):
\[
    \rho^{\downarrow}(B)\coloneqq  \rho_N(B)\qquad\text{and}\qquad
\rho^{\leftarrow}(B)\coloneqq \rot^{-1}(\rho_N(\rot(B))).
\]
\end{defn}

To illustrate the above definition, see \cref{ex:mRSK} below.

\begin{example}\label{ex:mRSK}

For the matrix $B\in\M(5,6)$ in the upper left, we show $\rho^{\leftarrow}(B)$, 
$\rho^{\downarrow}(B)\in\MLQ_0(\lambda,6)$, 
and $\rho^{\downarrow}(\rho^{\leftarrow}(B))=M(\lambda)\in\MLQ_0(\lambda,5)$ for $\lambda=(5,3,2,2,1)$. As $\rho^{\leftarrow}$ is equivalent to downward collapsing of a rotated multiline queue, we show the rotated pairing lines in $\rho^{\leftarrow}(B)$, viewing it as a rotated element of $\MLQ_0(\lambda',5)$. 

\begin{center}
        \begin{tikzpicture}
        \def \scale{0.3} 
        \node at (-2.25,0) {$B=$};
        \node at (8.5,0) {$=\rho^{\leftarrow}(B)$};
        \node at (-2.6,-4) {$\rho^{\downarrow}(B)=$};
        \node at (8.75,-4) {$=\rho^{\downarrow}(\rho^{\leftarrow}(B))$};
            \node (left) at (6,0) {\rotatebox{-90}{\includegraphics[scale = \scale]{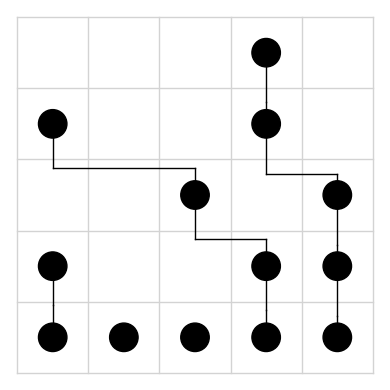}}};
            \node (original) at (0,0) {\includegraphics[scale = \scale]{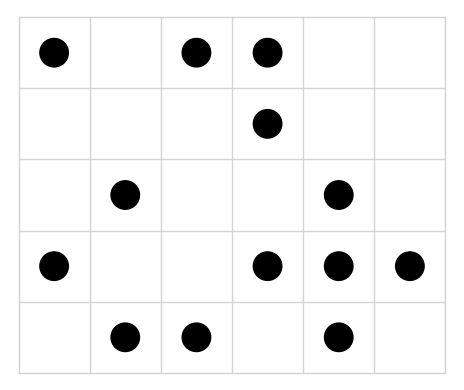}};
            \node (down) at (0,-4) {\rotatebox{0}{\includegraphics[scale = \scale]{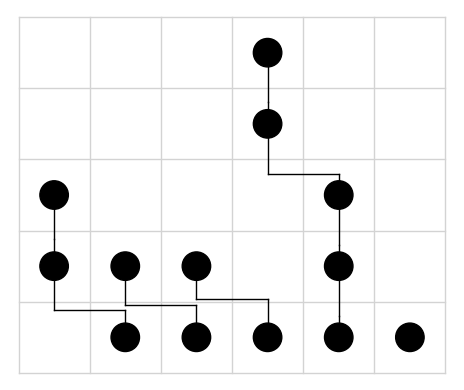}}};
            \node (double) at (6,-4) {\includegraphics[scale = \scale]{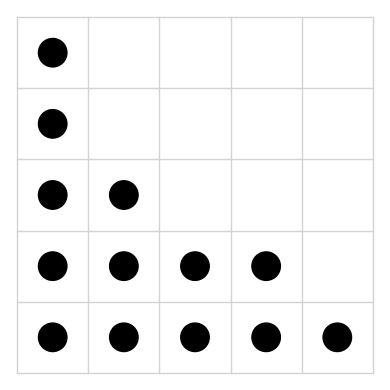}};
            \draw[->] (original) -- (left) node[above, midway] {\footnotesize{$\rho^\leftarrow$}};
            \draw[->] (original) -- (down) node[right, midway] {\footnotesize{$\rho^\downarrow$}};
            \draw[->] (down) -- (double) node[above, midway] {\footnotesize{$\rho^\leftarrow$}};
            \draw[->]  (left) -- (double) node[right, midway] {\footnotesize{$\rho^\downarrow$}};
            
        \end{tikzpicture}
    \end{center}

\end{example}

Our main observation is that the pair $(\rho^\leftarrow,\rho^\downarrow)$ defines a multiline queue analogue of the RSK correspondence, as established in the following theorem.

\begin{theorem}\label{thm:mRSK}
    Let $L,n$ be positive integers. Define the map 
    \begin{align*}
    \mRSK \;:\; \M(L,n)& \to \bigcup_{\lambda} \MLQ_0(\lambda,n)\times \MLQ_0(\lambda',L),\\
    B& \longmapsto (\rho^{\leftarrow}(B),\rho^{\downarrow}(B)),
    \end{align*}
    where the union is over partitions $\lambda$ with $\ell(\lambda)\leq n$ and $\ell(\lambda')\leq L$. Then $\mRSK$ is a bijection. 
\end{theorem}

To prove \cref{thm:mRSK} we need some preliminary results. For now, we mention that the theorem immediately gives an elementary proof of the dual Cauchy identity. Note that this result can also be viewed as a special case of the general integer-matrix case \cite{Knuth1970}.
\begin{corollary}[Dual Cauchy identity]\label{cor:dualCauchy}
    \[
    \sum_{\lambda} s_\lambda(x)s_{\lambda'}(y) = \prod_{i,j}(1+x_i y_j).
    \]
\end{corollary}

\begin{proof}
Assign to $B\in\M$ the weight function $\wt(B)=\prod_{i,j\geq 1}(x_iy_j)^{B_{i,j}}$, so that the right hand side is the weight generating function over $\M$. With the weight $\wt(N_1,N_2) = x^{N_1}y^{N_2}$ for $(N_1,N_2)\in\MLQ_0(\lambda)\times\MLQ_0(\lambda')$, the generating function of the codomain of $\mRSK$ is 
$$\displaystyle \sum_{\lambda}\sum_{\substack{N_1\in\MLQ_0(\lambda)\\N_2\in\MLQ_0(\lambda')}}x^{N_1}y^{N_2}.$$ From \eqref{eq:schur-mlqs}, this generating function equals the expression on the left hand side of the Cauchy identity. The identity follows since $\mRSK$ is weight-preserving.
\end{proof}

We now focus on proving \cref{thm:mRSK}. We show that orthogonal dropping operators commute, and then conclude that we can regard one of the collapsed multiline queues as a recording object in the usual RSK sense.

\begin{defn}
Define the operator $e_i^{d}$ for $d\in\{\uparrow\;,\;\downarrow\;,\;\leftarrow\;,\,\rightarrow\}$ as the dropping operator in the direction $d$ acting on $\M$ by moving all unmatched above balls from row $i$ to row $i-1$, where rows are numbered with respect to the orientation $d$. Given that we can rotate and reflect any matrix to apply the operators, we restrict to studying the downwards and leftwards directions. Explicitly, $e_i^{\downarrow}\coloneqq e_i^\star$ from \cref{sec:raising-lowering} and $e_i^{\leftarrow}\coloneqq \rot^{-1}\circ\, e_i^\star \circ \rot$.
\end{defn} 

For any direction $d$, we write $e_{[a,b]}^d \coloneqq e_b^d e_{b-1}^d\cdots e_a^d$ for the composition of a sequence of operators, and define
\[\rho_N^d (B) \coloneqq e_{[1,L-1]}^d e_{[1,L-2]}^d \cdots e_{[1,2]}^de_{[1,1]}^d (B).
\]

While commutation of the $e_i^\downarrow$ and $e_j^\leftarrow$–operators is already known (see, for instance, \cite[Lemma 1.3.7]{vanleeuwen2006double}), we include a direct proof to align with our particular setting.
\begin{lemma}\label{lem:localcrystalscommute_mlq}
    Let $B\in\M$. Then $e_i^\downarrow(e_j^\leftarrow(B)) = e_j^\leftarrow(e_i^\downarrow(B))$ for all $i$ and $j$. 
\end{lemma}

\begin{proof}
    Let $X$ be a (possibly empty) set of balls in $B_{i+1}$ that are unmatched above. Let $Y$ be the (possibly empty) set of balls in row $j+1$ of $\rot(B)$ (i.e.~column $j+1$ of $B$) that are unmatched above. Then $e_i^\downarrow(B)$ is $B$ with the balls in columns indexed by $X$ moved from row $i+1$ to row $i$, and $e_i^\leftarrow(B)$ is $B$ with the balls in rows indexed by $Y$ moved from column $j+1$ to column $j$. We consider the following four cases: (i) $j+1\in X$ and $i+1\in Y$, (ii) $j+1\in X$ and $i+1\not\in Y$, (iii) $j+1\not\in X$ and $i+1\in Y$, and (iv) $j+1\not\in X$ and $i+1\not\in Y$. 
    
    \medskip
    
    \noindent\textit{Case i.} $B$ corresponds to the following configuration on columns $j,j+1$.
    \begin{center}
    \begin{tikzpicture}[scale=0.5]
    \node at (0,0) {\tableau{\ &\bullet\\
               \ &\ }};

               \node at (-1.5,-.5) {\tiny $i$};
               \node at (-1.5,.5) {\tiny $i+1$};
               
    \end{tikzpicture}
    \end{center}
     Since the ball depicted is unmatched above in both $B$ and $\rot(B)$, $X\cup\{j\}\backslash\{j+1\}$ is the set of balls unmatched above in row $i+1$ of $e_i^\leftarrow(B)$ and $Y\cup\{i\}\backslash\{i+1\}$ is the set of balls unmatched above in row $j+1$ of $\rot(e_i^\downarrow(B))$, so the operators commute in this case.
    
    \medskip
    
     \noindent\textit{Case ii.} $B$ corresponds to one of the following configurations on columns $j,j+1$:
     \begin{center}
    \begin{tikzpicture}[scale=0.5]
    \def\h{7};
    \def\hh{14};
    \node at (-\h-3,0) {(a)};
    \node at (-\h,0) {\tableau{\bullet&\bullet\\
               \ &\ }};
               \node at (-\h-1.5,-.5) {\tiny $i$};
               \node at (-\h-1.5,.5) {\tiny $i+1$};
               \draw[black,thick] (-\h+0.5,0.425) -- (-\h-0.5,0.425);
    \node at (-3,0) {(b)};
    \node at (0,0) {\tableau{\bullet&\bullet\\
               \bullet&\ }};
               \node at (-1.5,-.5) {\tiny $i$};
               \node at (-1.5,.5) {\tiny $i+1$};
               \draw[black,thick] (0.5,0.425) -- (-0.5,0.425);
    \node at (\h-3,0) {(c)};
    \node at (\h+0,0) {\tableau{\ &\bullet\\
               \bullet&\ }};
               \node at (\h-1.5,-.5) {\tiny $i$};
               \node at (\h-1.5,.5) {\tiny $i+1$};
               \draw[black,thick] (\h+0.5,0.425) -- (\h + 0.2,0.425) -- (\h+0.2,-0.425) -- (\h-0.5,-0.425);
    \node at (\hh-3,0) {(d)};
    \node at (\hh+0,0) {\tableau{\ &\bullet\\
               \ &\ }};
               \node at (\hh-1.5,-.5) {\tiny $i$};
               \node at (\hh-1.5,.5) {\tiny $i+1$}; 
                \draw[black,thick] (\hh+0.5,0.425) -- (\hh + 0.2,0.425) -- (\hh+0.2,-1.1);
    \end{tikzpicture}
    \end{center}
    
    In configuration (a), if $(i+1,j+1)$ is unmatched above in $B$ then so is $(i+1,j)$, and so $j\in X$ as well. Thus the two balls will both move to row $i$ and still be paired to each other in column $i$ of $\rot(e_i^\downarrow(B))$. In configuration (b), the ball at site $(j+1,i)$ is matched above in $\rot(e_i^\downarrow(B))$, and the ball at site $(j,i+1)$ is matched below in $\rot(e_i^\downarrow(B))$ if and only if the ball at site $(j,i)$ is matched below in $\rot(B)$. In configurations (c) and (d), the ball at site $(j+1,i)$ is matched above in $\rot(e_i^\downarrow(B))$, and all other pairings are unaffected. Thus in all four cases, $Y$ is the set of balls unmatched above in row $j+1$ of $\rot(e_i^\downarrow(B))$ and there's no change to rows $i,i+1$ in $e_i^\leftarrow(B)$, and so the operators commute in this case as well.

    \medskip
    
    \noindent\textit{Case iii.} $B$ corresponds to one of the following configurations on columns $j,j+1$:
     \begin{center}
    \begin{tikzpicture}[scale=0.5]
    \def\h{7};
    \def\hh{14};
    \node at (-\h-3,0) {(e)};
    \node at (-\h,0) {\tableau{\ &\bullet\\
               \bullet &\bullet }};
               \node at (-\h-1.5,-.5) {\tiny $i$};
               \node at (-\h-1.5,.5) {\tiny $i+1$};
               \draw[black,thick] (-\h+0.425,0.5) -- (-\h+0.425,-0.5);
    \node at (-3,0) {(f)};
    \node at (0,0) {\tableau{\ &\bullet\\
               \ &\ }};
               \node at (-1.5,-.5) {\tiny $i$};
               \node at (-1.5,.5) {\tiny $i+1$};  
               \draw[black,thick] (0.425,0.5) -- (0.425,0.2) -- (1.1, 0.2);
               
    \node at (\h-3,0) {(g)};
    \node at (\h,0) {\tableau{\ &\bullet\\
               \ &\bullet }};
               \node at (\h-1.5,-.5) {\tiny $i$};
               \node at (\h-1.5,.5) {\tiny $i+1$};
               \draw[black,thick] (\h+0.425,0.5) -- (\h+0.425,-0.5);
    \end{tikzpicture}
    \end{center}
    The analysis of configurations (e) and (f) are similar to that of configurations (b) and (d) respectively from \textit{Case (ii)}. In configuration (g), whether or not $i\in Y$, the two balls will still be paired to each other in $e_i^{\leftarrow}(B)$, and since $e_i^{\downarrow}(B)$ doesn't affect columns $j$ and $j+1$, $X$ (resp. $Y$) is the set of balls in row $i+1$ (row $j+1$) of $e_i^{\leftarrow}(B)$ (resp. $e_i^{\downarrow}(B)$) that are unmatched above, and so the operators commute in this case as well.

    \medskip
    
    \noindent\textit{Case iv.} The only configurations in which $e_i^{\downarrow}$ may affect columns $j$ and $j+1$ or $e_i^{\leftarrow}$ may affect rows $i$ and $i+1$ are the following:
     \begin{center}
    \begin{tikzpicture}[scale=0.5]
    \def\aa{2.7};
    \def\h{5};
    \def\ha{10};
    \def\hh{15};
    \def\hha{20};
    \node at (-\h-\aa,0) {(h)};
    \node at (-\h,0) {\tableau{\ &\ \\
               \ &\bullet }};
               \node at (-\h-1.5,-.5) {\tiny $i$};
               \node at (-\h-1.5,.5) {\tiny $i+1$};
    \node at (-\aa,0) {(i)};
    \node at (0,0) {\tableau{\bullet&\ \\
               \ &\bullet}};
               \node at (-1.5,-.5) {\tiny $i$};
               \node at (-1.5,.5) {\tiny $i+1$};  
    \node at (\h-\aa,0) {(j)};
    \node at (\h,0) {\tableau{\ &\bullet\\
               \ &\bullet }};
               \node at (\h-1.5,-.5) {\tiny $i$};
               \node at (\h-1.5,.5) {\tiny $i+1$};
               \draw[black,thick] (\h+0.425,0.5) -- (\h+0.425,-0.5);

    \node at (\ha-\aa,0) {(j')};
    \node at (\ha,0) {\tableau{\ &\bullet\\
               \ &\bullet }};
               \node at (\ha-1.5,-.5) {\tiny $i$};
               \node at (\ha-1.5,.5) {\tiny $i+1$};
               \draw[black,thick] (\ha+0.425,0.425) -- (\ha+0.175,0.425) -- (\ha+0.175,-1.1);

     \node at (\hha-\aa,0) {(k')};
    \node at (\hha+0,0) {\tableau{\bullet&\bullet\\
               \ &\ }};
               \node at (\hha-1.5,-.5) {\tiny $i$};
               \node at (\hha-1.5,.5) {\tiny $i+1$}; 
               \draw[black,thick] (\hha + 0.5, 0.425) -- (\hha - 0.5, 0.425);

    \node at (\hh-\aa,0) {(k)};
    \node at (\hh+0,0) {\tableau{\bullet&\bullet\\
               \ &\ }};
               \node at (\hh-1.5,-.5) {\tiny $i$};
               \node at (\hh-1.5,.5) {\tiny $i+1$}; 
               \draw[black,thick] (\hh + 0.425, 0.5) -- (\hh + 0.425, 0.15) -- (\hh+1.1,0.15);
    \end{tikzpicture}
    \end{center}
    In the previous cases, $(j)$ and $(j')$ (as well as $(k)$ and $(k')$) correspond to the same scenario, but we separate them to show that the ball in position $(i+1,j+1)$ is matched above in $B$ and $\rot(B)$. It is a straightforward check that in all four cases, a ball is unmatched above in row $j+1$ (resp. row $i+1$) of $e_i^{\leftarrow}(B)$ (resp. $\rot(e_i^{\downarrow}(B))$) if the corresponding ball in $B$ (resp. $\rot(B)$ is unmatched above, which implies the two operators commute.

    Thus we can conclude that $e_j^\leftarrow\circ e_i^\downarrow$ and $e_i^\downarrow \circ e_j^\leftarrow$ coincide for any $i$ and $j$.
    \end{proof}

    As a consequence of the local commutativity of the operators we obtain the following result.

\begin{prop}\label{prop:double-collapse}
    Let $B\in\M$, and let $Q=\rho_Q(B)$ be the recording tableau of the standard collapsing of $B$, with $Q\in\SSYT(\mu')$ for some partition $\mu$. Then the following hold. 
    \begin{itemize}
    \item[i.] $\rho_N^{\downarrow}(B)\in\MLQ_0(\mu)$, $\rho_N^{\leftarrow}(B)\in\MLQ_0(\mu')$, and their double collapsing satisfies $\rho_N^{\downarrow}(\rho_N^{\leftarrow}(B)) = \rho_N^{\leftarrow}(\rho_N^{\downarrow}(B)) = M(\mu)\in\MLQ_0(\mu)$. 
   \item[ii.] $\cw(\rho_N^{\leftarrow}(B))=\crw(Q)$.
   \end{itemize}
\end{prop}

\begin{proof}
    The equality $\rho^{\downarrow}_N(\rho^{\leftarrow}_N(B)) = \rho^{\leftarrow}_N(\rho^{\downarrow}_N(B))$ is immediate from \cref{lem:localcrystalscommute_mlq}, since $\rho^\leftarrow_N$ and $\rho^\downarrow_N$ are built from sequences of $e_i^\leftarrow$'s and $e_i^\downarrow$'s. 
 The fact that $\rho^\downarrow_N(B)$ has the same shape $\mu$ as the conjugate shape of the recording tableau follows from the construction of the tableau from \cref{def:collapsing_mlqs}. Since $\rho^\leftarrow_N$ preserves the number of balls in each row, $\rho^\leftarrow_N(\rho^\downarrow_N(B))\in\MLQ_0(\mu)$, as well. This in turn means that $\rot(\rho^\downarrow_N(\rho^\leftarrow_N(B)))\in\MLQ_0(\mu')$, which must have the same shape as $\rot(\rho^\leftarrow_N(B))$ since $\rho^\downarrow_N$ preserves column content. Finally, the only possible configuration for $\rho^\downarrow_N(\rho^\leftarrow_N(B))$ having shape $\mu$ and its rotation having shape $\mu'$ is the configuration $M(\mu)$.  

    To prove (ii.), we use labeled collapsing from \cref{def:labeled_collapsing}. Let $Q'$ be the unique semistandard tableau whose row content is identified with the labels in the labeled (downward) collapsing of $\rho_N^{\leftarrow}(B)$. 
    By \cref{prop:labeled-collapsing-Q}, $Q'=\rho_Q(\rho_N^{\leftarrow}(B))$. Denote $B'=\rho^{\leftarrow}_N(B)$ and $B''=\rho^{\downarrow}_N(B')$. Again due to the commutativity of the orthogonal operators, we have 
    \[
        \rho^{-1}(B'',Q')=B'=\rho^{\leftarrow}_N(B)=\rho^{\leftarrow}_N(\rho^{-1}(\rho^{\downarrow}_N(B),Q))=\rho^{-1}(\rho^{\leftarrow}_N(\rho^{\downarrow}_N(B)),Q)
        =\rho^{-1}(B'',Q),
    \]
    which implies that $Q'=Q$ since $\rho^{-1}$ is a bijection.
    
     We claim that during the labeled collapsing of $\rho^\leftarrow(B)$, every ball preserves its original tracking label (i.e. no swapping of tracking labels occurs). In the proof of \cref{lem:collapsing_condition}, we showed that at each step of the collapsing procedure, each particle is sent to its corresponding terminal row without bumping any other particle; in terms of labeled collapsing, this means no tracking labels are exchanged in the process. Since $B'=\rho^\leftarrow(B)$ collapses to $M(\mu)$, whose row lengths are given by $\mu'$ (the row lengths of $Q$), the number of particles in row $r$ with tracking label $\ell$ in $\rho^\leftarrow(B')$ precisely corresponds to the number of entries $\ell$ in row $r$ of $Q$.  Finally, since the configuration of $B'$ is bottom and left-justified and all particles have tracking label equal to their original row number in $B$, we immediately obtain that the column reading word of $\rho^\leftarrow(B)$ is identically the column reading word of $Q'=Q$, that is,  $\cw(\rho_{N}^{\leftarrow}(B))=\crw(Q')=\crw(Q)$, proving item (ii). 
\end{proof}

In particular, \cref{prop:double-collapse} implies that the map $\mRSK$ is well-defined. We show the map is a bijection by constructing an inverse.

\begin{proof}[Proof of \cref{thm:mRSK}]

    We will describe an inverse to $\mRSK$ by regarding $\rho^\downarrow(B)$ and $\rho^\leftarrow(B)$ as the insertion and recording components of the correspondence, respectively.

    Let $(M_1,M_2)\in \MLQ_0(\lambda,n)\times \MLQ_0(\lambda',L)$. 
    Considering $\rot(M_2)$ as a binary matrix with row sums equal to the parts of $\lambda$, let $Q_2=\rho_Q(\rot(M_2))$ be the recording tableau of its collapsing. By \cref{prop:double-collapse}, we have that $\rho^{\leftarrow}_N(M_1) = \rho^{\downarrow}_N(\rot(M_2))$,  
    and so $Q_2\in\MLQ_0(\lambda',L)$. Then, we can apply the inverse of map $\varphi$ from \cref{thm:wrapping-to-nonwrapping} to $(M_1,Q_2)\in \MLQ_0(\lambda,n)\times \SSYT(\lambda',L)$ to obtain a binary matrix $M=\varphi^{-1}(M_1,Q_2)\in\M(L,n)$. 
    The fact that this is indeed the inverse map to $\mRSK$ follows from the commutativity of $e^\leftarrow_i$ and $e^\downarrow_i$ which define $\rho^{\leftarrow}_N$ and $\rho^{\downarrow}_N$, respectively: indeed, if $M'\in\M(L,n)$ such that $\mRSK(M')=(M_1,M_2)$, then $\rho_Q(M')=\rho_Q(\rot(M_2))=Q_2$, and so $M'=M$ as desired. 
\end{proof} 

\begin{example}\label{ex:proof of 5.8}

We illustrate \cref{prop:labeled-collapsing-Q} with the twisted multiline queue from \cref{ex:mRSK}. The tracking labels on $B$ and $\rho^{\leftarrow}(B)$ are shown below (not to be confused with the labeled multiline queue $L_G(B)$), together with the results of the labeled collapsing of each. Observe that the column content of these are the same, and we obtain $\rho_Q(B)$ by identifying the tracking labels in $\rho^\downarrow_N(\rho^\leftarrow_N(B))$ with a filling of a Young diagram with the same row-shape. 

\begin{center}
    \begin{tikzpicture}[scale=0.5]
        \def \w{1};
        \def \h{1};
        \def \r{0.3};
        \begin{scope}[xshift=0cm]
        \node at (-1.1,2.5) {$B = $};
        \foreach \i in {0,...,5}
        {
        \draw[gray!50] (0,\i*\h)--(\w*6,\i*\h);
        }
        \foreach \i in {0,...,6}
        {
        \draw[gray!50] (\w*\i,0)--(\w*\i,5*\h);
        }
        \foreach \xx\yy\c in {1/0/1,2/0/1,4/0/1,0/1/2,3/1/2,4/1/2,5/1/2,1/2/3,4/2/3,3/3/4,0/4/5,2/4/5,3/4/5}
        {
       \filldraw[black,fill=red!0] (\w*.5+\w*\xx,\h*.5+\h*\yy) circle (\r cm);
        \node at (\w*.5+\w*\xx,\h*.5+\h*\yy) {\tiny \c};
        }
        \end{scope}

        \draw[->] (7,2.5*\h) -- (9,2.5*\h) node[above, midway] {\footnotesize{$\rho^\leftarrow$}};

        \begin{scope}[xshift=10cm]
        \foreach \i in {0,...,5}
        {
        \draw[gray!50] (0,\i*\h)--(\w*6,\i*\h);
        }
        \foreach \i in {0,...,6}
        {
        \draw[gray!50] (\w*\i,0)--(\w*\i,5*\h);
        }
        \foreach \xx\yy\c in {0/0/1,1/0/1,2/0/1,0/1/2,1/1/2,3/1/2,4/1/2,0/2/3,2/2/3,0/3/4,0/4/5,1/4/5,3/4/5}
        {
       \filldraw[black,fill=red!0] (\w*.5+\w*\xx,\h*.5+\h*\yy) circle (\r cm);
        \node at (\w*.5+\w*\xx,\h*.5+\h*\yy) {\tiny \c};
        }
        \node at (7.9,2.5) {$=\rho^{\leftarrow}(B)$};
        \end{scope};

        \begin{scope}[yshift=-7cm]
        \foreach \i in {0,...,5}
        {
        \draw[gray!50] (0,\i*\h)--(\w*6,\i*\h);
        }
        \foreach \i in {0,...,6}
        {
        \draw[gray!50] (\w*\i,0)--(\w*\i,5*\h);
        }
        \foreach \xx\yy\c in {1/0/1,2/0/1,3/0/2,4/0/1,5/0/2,0/1/3,1/1/2,2/1/5,4/1/2,0/2/5,4/2/3,3/3/4,3/4/5}
        {
       \filldraw[black,fill=red!0] (\w*.5+\w*\xx,\h*.5+\h*\yy) circle (\r cm);
        \node at (\w*.5+\w*\xx,\h*.5+\h*\yy) {\tiny \c};
        }
        \end{scope}

        \draw[->] (3,-.1) -- (3,-1.9) node[left, midway] {\footnotesize{$\rho^\downarrow$}};

        \begin{scope}[xshift=10cm, yshift = -7cm]
        \foreach \i in {0,...,5}
        {
        \draw[gray!50] (0,\i*\h)--(\w*6,\i*\h);
        }
        \foreach \i in {0,...,6}
        {
        \draw[gray!50] (\w*\i,0)--(\w*\i,5*\h);
        }
        \foreach \xx\yy\c in {0/0/1,1/0/1,2/0/1,0/1/2,1/1/2,3/0/2,4/0/2,0/2/3,2/1/3,0/3/4,0/4/5,1/2/5,3/1/5}
        {
       \filldraw[black,fill=red!0] (\w*.5+\w*\xx,\h*.5+\h*\yy) circle (\r cm);
        \node at (\w*.5+\w*\xx,\h*.5+\h*\yy) {\tiny \c};
        }

        \draw[->] (\w*7,2.5*\h)--(\w*8,2.5*\h);

        \node at (\w*11,2.5*\h) {\tableau{5\\4\\3&5\\2&2&3&5\\1&1&1&2&2}};
        \node at (\w*14,2.5*\h) {$=\rho_Q(B)$};
        \end{scope};

        \draw[->] (13,-.1) -- (13,-1.9) node[right, midway] {\footnotesize{$\rho^\downarrow$}};
    \end{tikzpicture}
\end{center}
\end{example}

\begin{lemma}\label{lem:maj_left_coll} 
  Let $M$ be a (twisted) multiline queue. Then
    \[\maj_G(\rho^{\leftarrow}_N(M))=\maj_G(M).
    \]
\end{lemma}

\begin{proof}
Let $M=(B_1,\ldots,B_k)$ and let $u_j^{(i+1)}$ be the set of balls labeled $\geq j$ in row $i+1$ of $L_G(M)$. By \cref{prop:local_majG_energy}, we write $\maj_G(M)$ in terms of the \emph{energy function} from \cref{def:energy_function} as
\begin{equation}\label{eq:dec}
\maj_G(M)=\sum_i\sum_j H(u_j^{(i+1)};B_i).
\end{equation}
Recall that $H(A;C)$ is equal to the number of balls paired above in $\Par^c(A,C)$ minus the number of balls paired above in $\Par(A,C)$. Since collapsing left leaves the total number of particles (of each label) in each row unchanged, $\Par^c(B_i,u_j^{(i+1)})$ is invariant of $\rho^\leftarrow$ for each component in the decomposition \eqref{eq:dec}. On the other hand, the number of matched/unmatched above particles in $\Par(B_i,u_j^{(i+1)})$ is a function of the set of $e_i^\downarrow$ operators applied to obtain $\rho^\downarrow_N$ (specifically, of the sequence of the $\phi(r,j)$'s from the proof of \cref{thm:wrapping-to-nonwrapping}). Since the $e_j^\leftarrow$'s and $e_i^\downarrow$'s commute, this set of operators is invariant of $\rho^\leftarrow_N$, and hence so is the right hand side of \eqref{eq:dec}.    
\end{proof}

\begin{example}\label{ex:rot} 
    For $\lambda=(5,4,3,2)$, we show the multiline queues $M\in\MLQ(\lambda,7)$ and  $\rho^{\leftarrow}_N(M)=\rot^{-1}(\rho_N(\rot(M)))\in\MLQ(\lambda,5)$. Notice that $\rho_N(\rot(M))\in \MLQ_0(\lambda',5)$ is nonwrapping, which can be seen from the rotated pairing lines.

    \begin{center}
        \begin{tikzpicture}
        \def \scale{0.3} 
        \node at (-2.5,0) {$M=$};
        \node at (9.5,0) {$=\rho^{\leftarrow}(M)$};
            \node (original) at (0,0) {\rotatebox{0}{\includegraphics[scale = \scale]{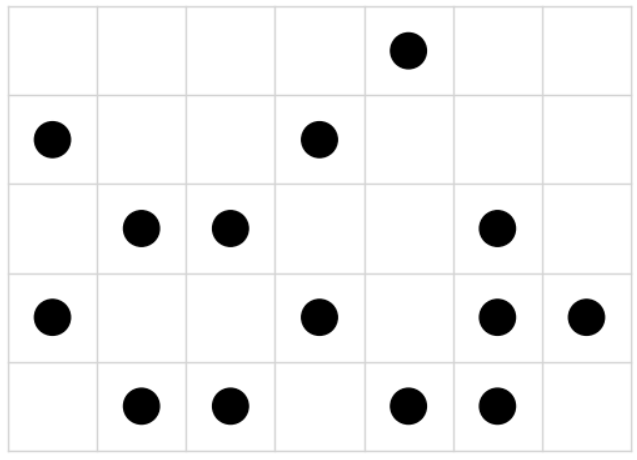}}};

            \node (leftward) at (7.25,0){\rotatebox{270}{\includegraphics[scale = \scale]{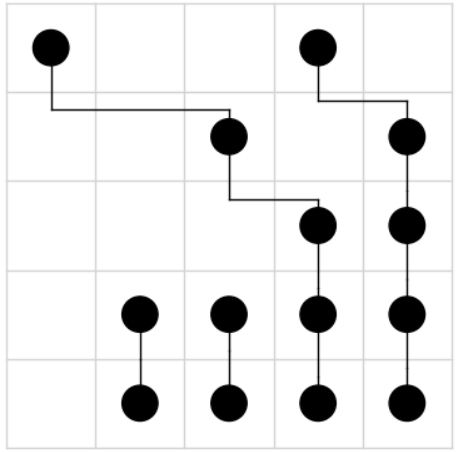}}};
            
            \draw[->]  (3,0)--(5,0) node[above, midway] {$\rho^\leftarrow$};
        \end{tikzpicture}
    \end{center}
    However, the $\maj$ of $\rho^{\leftarrow}_N(M)$, which is the rotation of $\rho_N(\rot(M))$, is equal to $\maj(M)=3$, as we can see from the pairing lines drawn below.

    \begin{center}
        \begin{tikzpicture}
        \def \scale{0.3} 
        \node at (-2.5,0) {$M=$};
        \node at (9.5,0) {$=\rho^{\leftarrow}(M)$};
        \node (downward) at (0,0) {\includegraphics[scale = \scale]{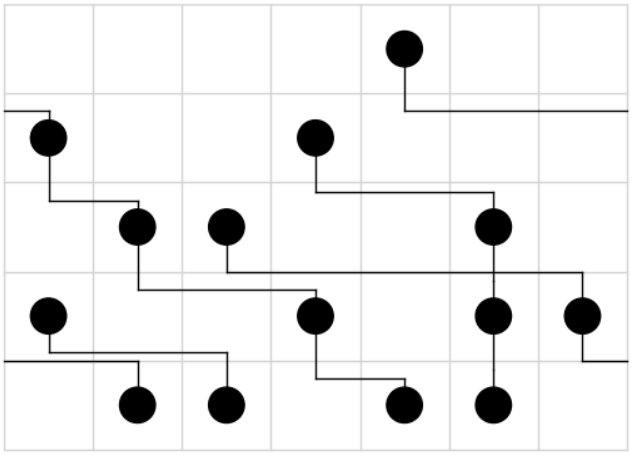}};
        \node (double) at (7.25,0) {\includegraphics[scale = \scale]{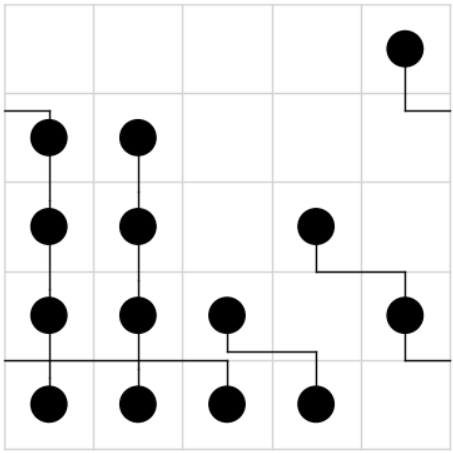}};
 
        \end{tikzpicture}
    \end{center}
    
\end{example}

To conclude this section, we will use double collapsing to give simple proofs of \cref{thm:charge,lem:pairing-position-word-Q}. We start by proving \cref{thm:charge}, which states: for a (straight) multiline queue $M$ with $Q=\rho_Q(M)$, $\maj(M)=\charge(Q)$.
\begin{proof}[Proof of \cref{thm:charge}]
    Let $M' = \rho^{\leftarrow}_N(M)$, and let $Q'=\rho_Q(M')$. By \cref{lem:charge cw}, $\maj(M) = \charge(\cw(Q)).$ Moreover, from \cref{prop:double-collapse}, $\cw(M')=\crw(Q')=\crw(Q)$.  Therefore,  
    \[
    \maj(M)=\maj(M') = \charge(\cw(M')) = \charge(\crw(Q)) = \charge(Q)
    \]
    where the first equality is from \cref{lem:maj_left_coll}, the second from \cref{lem:charge cw}, and the fourth from \cref{lem:rw colw}.
\end{proof}

The same argument yields a proof of \cref{lem:pairing-position-word-Q}: for a binary matrix (or twisted multiline queue) $B$ with $Q = \rho_Q(B)$, $\Pair_i(B_{i+1},B_{i}) = \Pair_i(\crw(Q))$.

\begin{proof}[Proof of \cref{lem:pairing-position-word-Q}]
    Let $B' = \rho^{\leftarrow}_N(B)$. From \cref{prop:double-collapse} we know that $\cw(B') = \crw(Q)$. Since $\Pair_i(B'_{i+1},B'_{i}) = \Pair_i(\cw(B'))$, it is enough to show that $\Pair_i(B_{i+1},B_{i})=\Pair_i(B'_{i+1},B'_{i}).$ However, we have that both $e_i^\downarrow$ and $e_i^\uparrow$ commute with $\rho^{\leftarrow}$. Moreover, $\Pair_i$ can be computed from the number of times that the operators $e_i^\downarrow$ and $e_i^\uparrow$ can be applied to $(B_{i+1},B_{i})$ combined with the sizes of the original sets, since this information determines the number of balls that are matched and unmatched above and below. Then the equality holds from \cref{lem:localcrystalscommute_mlq}. 
\end{proof}

\subsection{Connection to crystal operators on semistandard key tabloids}\label{sec:key tabloids}
The leftward dropping operator $e_i^{\leftarrow}$ corresponds to the raising crystal operator introduced in \cite{AssafGonzales21} on \emph{semistandard key tabloids}. We explain the connection in the context of $q$-Whittaker polynomials.

A semistandard key tabloid corresponding to a partition $\lambda$ is a filling $\sigma\colon \dg(\lambda)\rightarrow \mathbb{N}_{>0}$ in which there are no repeated entries within any row, and which satisfies a certain non-attacking and coinversion-free condition on the entries, obtained from restricting the set of non-attacking Haglund--Haiman--Loehr tableaux of \cite{HHL08} to those which are coinversion-free.\footnote{To make the connection with multiline queues, our definition of key tabloids is a $90^{\circ}$ counterclockwise rotation of that in \cite{AssafGonzales21}. Moreover, their objects are defined more generally for composition shapes, but we restrict to the partition case.} 
These objects are similar to those we describe in \cref{sec:tableaux}, except that the coinversion-free property results in a different order of entries within each row. 

The major index statistic on semistandard key tabloids is given by \eqref{eq:maj tab}. 
There is a unique semistandard key tabloid satisfying the coinversion-free condition for every given row content. Omitting the details of the definition, let us identify a multiline queue with the unique semistandard key tabloid that has the same row content. Similar to the proof of the analogous fact in \cref{thm:queue to tableau} for \emph{coquinv-free tableaux}, for key tabloids it is also the case that $\maj(T)=\maj(M)$, where $T$ is a key tabloid associated to the multiline queue $M$ with the same row content as $T$ (see \cite[Lemma 5.14]{CMW18} for details, which proves the statement for tableaux which coincide with Haglund--Haiman--Loehr tableaux in the $t=0$ case, and are in bijection with multiline queues).

\begin{lemma}
    Let $e_i$ be the raising operator on key tabloids defined in \cite{AssafGonzales21}. Let $M=(B_1,\ldots,B_L)$ be a multiline queue and let $T$ be a key tabloid with row content given by the sets $B_1,\ldots,B_L$. Then the multiline queue corresponding to $e_i(T)$ is equal to $e_i^{\leftarrow}(M)$.
\end{lemma}
\begin{proof}
    Let the row content of $\rot(M)$ be $B'=(B'_1,\cdots,B'_n)$. Comparing definitions confirms that the $i$-pairing process defined in \cite{AssafGonzales21} to identify the unpaired $i+1$'s in $T$ is equivalent to $\Par_i(\cw(B'))$. Thus the action of the operator $e_i$ on $T$ is equivalent to the action of the operator $e^\downarrow_i$ on $\rot(M)$, and so the multiline queue corresponding to $e_i(T)$ is $\rot^{-1}(e_i(\rot(M)))=e_i^{\leftarrow}(M)$.
\end{proof}

\subsection{Bosonic multiline queue RSK via commuting crystal operators}\label{sec:bMLQ RSK}

The bijection from \cref{thm:mRSK} naturally extends to the set $\mathcal{M}$ of \emph{positive integer} matrices with finite support. Such matrices with weakly decreasing row sums correspond to bosonic multiline queues, just as binary matrices correspond to multiline queues. Given a positive integer matrix $B=(b_{ij})$ of size $L\times n$ whose row sums from bottom to top are weakly decreasing, let $D_j$ be the multiset where each $i\in[n]$ has multiplicity $b_{ji}$. Then $D=(D_1,\ldots,D_L)$ is the bosonic multiline queue corresponding to $B$. 

Recall from \cref{def:MLDpairing} that bosonic multiline queue pairing is done strictly to the left via a modified parenthesis matching. Modifying the parenthesis word $\widetilde\cw(D)$ from \cref{def:cw_bmlqs} by scanning the columns from \emph{right to left} instead of left to right, we can define a strictly-to-the-right bosonic multiline queue pairing. The interplay between the two pairing directions will be important when defining the RSK analog for bosonic multiline queues in \cref{thm:dRSK}.

\begin{definition}
    Let $P\in\{L,R\}$ (left or right) be a direction of pairing. The set of bosonic multiline queues of shape $\lambda$ on $n$ columns with pairing direction $P$ is denoted by $\MLD_{P}(\lambda,n)$. Similarly, the set of nonwrapping bosonic multiline queues with pairing in direction $P$ is denoted by $\MLD_{0,P}(\lambda,n)$.
\end{definition}

As with multiline queues, the parenthesis matching defines a collapsing procedure that gives rise to a bijection between wrapping objects and pairs of nonwrapping objects and semistandard tableaux of fixed shape and content. 

\begin{theorem}\label{thm:wrapping-to-nonwrapping-general-diagram-version}
    Let $\lambda$ be a partition and $P \in\{ L,R \}$ be the pairing direction. The downward collapsing with pairing direction $P$ defines a bijection $$\tilde\rho^{\downarrow}_{P} \, : \, \MLD_P(\lambda,n) \longrightarrow \bigcup_{\mu}\MLD_{0,P}(\mu,n)\times \SSYT(\mu',\lambda')$$
\end{theorem}

The proof of this theorem is analogous to the proof of \cref{thm:wrapping-to-nonwrapping} for collapsing of multiline queues via raising and lowering operators $e_i,f_i$ after we define the analogous operators acting on bosonic multiline queues. 

\begin{definition}
    For a pairing direction $P \in\{ L,R \}$, let $\tilde{e}_{P,i}^\downarrow$ be the operator acting on matrices $D\in\mathcal{M}$ that drops the ball in row $i+1$ unmatched above according to $P$ that is the furthest in the opposite direction of pairing. Extending the definition of the rotation operators $\rot$ on matrices we obtain a definition of leftward operators $\tilde{e}_{P,i}^\leftarrow = \rot^{-1}\circ\;\tilde{e}_{P,i}^\downarrow \circ \rot$.
\end{definition}

\begin{remark}
    In the previous definition, if $P=L$, the operator $\tilde{e}_{P,i}^\downarrow$ drops the rightmost particle that is unmatched above when pairing strictly to the left, and if $P=R$, the operator drops the leftmost particle that is unmatched above when pairing strictly to the right. 
\end{remark}

The bosonic analogue of \cref{lem:localcrystalscommute_mlq}, which can also be deduced from \cite{vanleeuwen2006double}, can be proved directly, in the same way as the fermionic case.

\begin{lemma}\label{lem:localcrystalscommute_mld}
    Let $D\in\mathcal{M}$. Then $\tilde e_{L,i}^\downarrow(\tilde e_{R,j}^\leftarrow(D)) = \tilde e_{R,j}^\leftarrow(\tilde e_{L,i}^\downarrow(D))$ for all $i$ and $j$. 
\end{lemma}

A delicate point in the bosonic case is that in order to commute, the leftward and downward operators must be defined for opposite pairing directions, see \cref{rem:opposite directions}.

\begin{remark}\label{rem:opposite directions}
Notably, in the bosonic case, keeping the same pairing direction in both downward and leftward collapsing does not yield a common double collapsing of the matrix. For instance, if we let $D = \begin{pmatrix} 0 & 1 \\ 2 & 0 \\ \end{pmatrix}$. Then $\tilde e_{R,1}^\leftarrow(\tilde e_{R,1}^\downarrow(D)) = \begin{pmatrix} 0 & 0 \\ 3 & 0 \\ \end{pmatrix}$ while $\tilde e_{R,1}^\downarrow(\tilde e_{R,1}^\leftarrow(D)) = \begin{pmatrix} 0 & 0 \\ 2 & 1 \\ \end{pmatrix}$.
\end{remark}

We use this property to give the bosonic multiline queue analog of the RSK correspondence, following the same argument as in the proof of \cref{thm:mRSK}.

\begin{theorem}\label{thm:dRSK}
    Let $L,n$ be positive integers, and let $\mathcal{M}(L,n)$ represent the set of $L\times n$ nonnegative integer matrices. The following map, given by given by $\dRSK(B) = (\rho_L^{\downarrow}(B),\rho_R^{\leftarrow}(B))$, is a bijection: 
    \[
    \dRSK \;:\; \mathcal{M}(L,n) \to \bigcup_{\lambda} \MLD_{0,L}(\lambda,n)\times \MLD_{0,R}(\lambda,L),
    \]
    where the union is over partitions $\lambda$ with {$\lambda_1 \leq \min\{L,n\}$}. 
\end{theorem}

\begin{prop}\label{thm:double-collapsing-diagrams}
     For a matrix $B\in\mathcal{M}(n,L)$, the double collapsing $C = \tilde \rho_R^{\leftarrow}(\tilde \rho_L^{\downarrow}(B)) = \tilde \rho_L^{\downarrow}(\tilde \rho_R^{\leftarrow}(B))\in\mathcal{M}$ yields an anti-diagonal matrix. Moreover, $\dRSK$ is well defined, that is, $\dRSK(B)\in\MLD_0(\lambda,n)\times \MLD_0(\lambda,L)$ for some partition $\lambda$, and in this case $C$ is such that $\lambda$ can be read off the anti-diagonal.
\end{prop}

\begin{proof}
    For $B\in\mathcal{M}(n,L)$, $\tilde \rho_L^\downarrow(B)$ has no nonzero entries above the anti-diagonal. Since vertical pairing is not allowed, any particle in the first column necessarily collapses to the lower left corner. Thus, a particle in the second column will necessarily collapse to the anti-diagonal or lower. Repeating the argument for each row, it holds that all particles in $\tilde \rho_L^\downarrow(B)$ collapse to the anti-diagonal or below.

    Now, in $\tilde \rho_L^\downarrow(B)$, the row sums are weakly decreasing from bottom to top. Therefore, when collapsing to the left by $\rho_L^\downarrow$, all particles from the bottom row end up in the lower left corner, and in subsequent rows, all particles pair with the particles in the row below, and are thus pushed to the anti-diagonal. Hence all the particles in $C = \rho_R^\leftarrow(\rho_L^\downarrow(B))$ are in the anti-diagonal. The same argument shows that $C' = \tilde \rho_L^{\downarrow}(\tilde \rho_R^{\leftarrow}(B))$ is an anti-diagonal matrix, and $C=C'$ by \cref{lem:localcrystalscommute_mld}. 
    
    Thus we showed that $\tilde \rho_L^\downarrow(B) \in \MLD_0(\lambda,n)$, where $\lambda'_i$ is the sum of row $i$ of this collapsed bosonic multiline queue. Moreover, the double collapsing yields an anti-diagonal matrix with $\lambda'_i$ in the corresponding entry of the anti-diagonal in row $i$. Since the double collapsing can be performed in the opposite order by \cref{lem:localcrystalscommute_mld}, this also means that $\tilde \rho_R^\leftarrow(B) \in \MLD_0(\lambda,L)$ for the same partition $\lambda$. 
\end{proof}

As a corollary of \cref{thm:mlqs-to-mlds,thm:dRSK}, we obtain the well-known Cauchy identity for integer matrices. 

\begin{corollary}[Cauchy identity]\label{cor:Cauchy}
    $$\sum_{\lambda}s_\lambda(X)s_\lambda(Y) = \prod_{i,j} \frac{1}{1-x_iy_j}.$$
\end{corollary}

We conclude this subsection with a bijection between nonwrapping multiline queues and nonwrapping bosonic multiline queues. For $M\in\MLQ_0(\lambda,n)$, let $w = w_1w_2\ldots w_k = \rw(M)$ and $B_M\in\mathcal{M}(k,n)$ be the binary matrix (considered as a positive integer matrix) with $(B_M)_{ij} = 1$ if and only if $w_{k-i+1} = j$. 

\begin{theorem}\label{thm:mlqs-to-mlds}
    The map $M\longmapsto \rho_L^{\downarrow}(B_M)$ is a bijection from $\MLQ_0(\lambda,n)$ to $\MLD_{0,L}(\lambda',n)$.
\end{theorem}

\begin{proof}
    The inverse map is given by $D \longmapsto \mlq(T_D)$ where $T_D$ is the tableau corresponding to the nonwrapping bosonic multiline queue $D$ as in \cref{lem:schur-mld}.
\end{proof}

\begin{example}
    We show the downward and leftward collapsing, as well as the double collapsing result, for the matrix in $\mathcal{M}(5,6)$ below.  

    \begin{center}
        \begin{tikzpicture}
        \def \scale{0.3} 
            \node (left) at (8,0) {\rotatebox{-90}{\includegraphics[width=4cm]{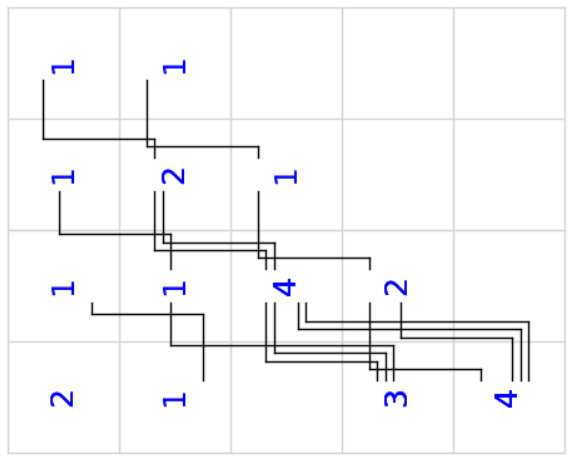}}};
            \node (original) at (0,0) {\includegraphics[width=4cm]{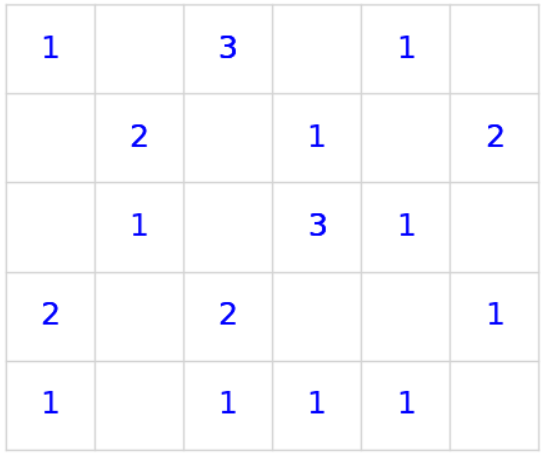}};
            \node (down) at (0,-4.5) {\rotatebox{0}{\includegraphics[width=5cm]{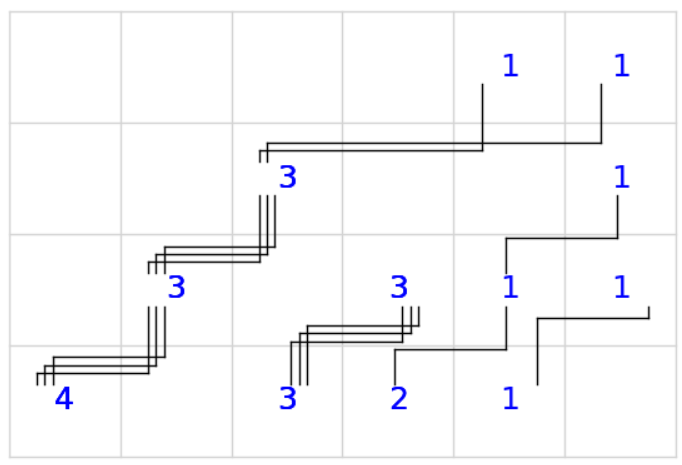}}};
            \node (double) at (8,-4.5) {\includegraphics[width=4cm]{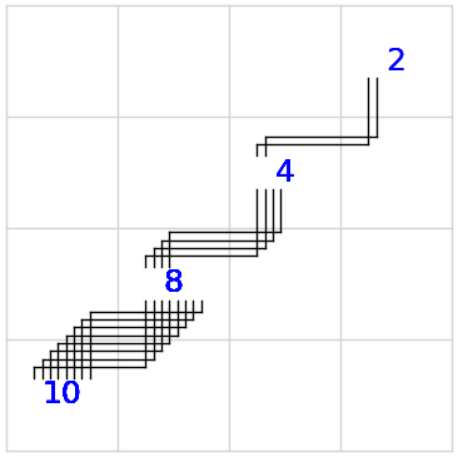}};
            \draw[->] (original) -- (left) node[above, midway] {\footnotesize{$\tilde \rho_R^\leftarrow$}};
            \draw[->] (original) -- (down) node[right, midway] {\footnotesize{$\tilde \rho_L^\downarrow$}};
            \draw[->] (down) -- (double) node[above, midway] {\footnotesize{$\tilde \rho_R^\leftarrow$}};
            \draw[->]  (left) -- (double) node[right, midway] {\footnotesize{$\tilde \rho_L^\downarrow$}};
            
        \end{tikzpicture}
    \end{center}
\end{example}

\subsection{Several formulas for $K_{\lambda\mu}(q,0)$}

In this subsection, we present several variations of formulas for $K_{\lambda\mu}(q;0)$ in terms of nonwrapping multiline queues and bosonic multiline queues. Each of these formulas follows from some version of collapsing on these objects.

\begin{defn}
    For partitions $\nu$ and $\eta$, let $\MLQ(\nu,\eta)=\{M\in\MLQ(\nu)\colon x^M=x^\eta\}$ be the set of multiline queues of shape $\nu$ and content $\eta$. Let $M(\nu)$ be the (unique) left-justified multiline queue of shape $\nu$. That is, $M(\nu) = (Q_1,\ldots,Q_{\nu_1})$ with $Q_j = \{ 1,2,\ldots,\nu'_j\}$ for $1\leq j\leq \nu_1$. In particular, $\MLQ(\nu,\nu)=\{M(\nu)\}$.
\end{defn}

\begin{lemma}\label{lem:collapsing_condition}
    Let $Q$ be a multiline queue with column content $\lambda$. Then $\rho_N(Q) = M(\lambda)$ if and only if $\rw(Q)$ is lattice. 
\end{lemma}

\begin{proof}

    Suppose that $\rw(Q)$ is a lattice word. We show that $\rho_N(Q) = M(\lambda)$ by induction on the length of $\rw(Q)$. The base case for $|\rw(Q)|=1$ is trivially left justified since $\rw(Q) = 1$ is the only lattice word of length one. Now suppose $\rw(Q) = r_1\,r_2\,\ldots\,r_{n-1}\,r_{n} =: w \, r_n.$ Since $w$ is lattice, by induction, collapsing of the multiline queue restricted to $w$ yields $M(\eta)$ for some partition $\eta$. The lattice condition on $\rw(Q)$ implies the last particle, which lies in column $r_n$, collapses to an outer corner of the diagram of $\eta$, since otherwise, the number of particles in column $r_n-1$ would have to be be strictly smaller than the number of particles in column $r_n$ -- a contradiction. Therefore, collapsing the particle in column $r_n$ on $M(\eta)$ yields the left-justified multiline queue $M(\lambda)$. 

    For the converse, suppose that $\rw(M)$ is not lattice, and take the first initial segment $r'$ of $\rw(Q)$ in which the lattice condition does not hold. Then there exists $c\geq 2$ such that $r'$ has $j$ "$c$"s and fewer than $j$ "$c-1$"s. By minimality of $r'$, $r' = w\,c$ where $w$ is a (possibly empty) lattice word. Consider the collapsing of the particles recorded by $r'$, yielding the partial multiline queue $Q'$. By induction, collapsing the particles corresponding to $w$ yields a left-justified multiline queue. Thus, when the particle in column $c$ collapses, it stops at site $(j,c)$, leaving a gap at site $(j,c-1)$. Moreover, once a partially collapsed multiline queue fails to be left-justified, further collapses cannot restore this property. This completes the proof of the claim.
\end{proof}

\begin{corollary}\label{cor:KF_mlq_version_1}
    Let $\lambda$ and $\mu$ be partitions. Then 
    \[
    K_{\lambda\mu}(q,0) = \sum_{\substack{B\in\MLQ(\mu',\lambda') \\ \collapse(B) = M(\lambda')}}q^{\maj(B)} = \sum_{\substack{B\in\MLQ(\mu',\lambda') \\ \rw(B)\text{ is lattice}}}q^{\maj(B)}.
    \]
\end{corollary}

\begin{proof}
For some partitions $\lambda,\mu$, consider the set $\MLQ_0(\mu)\times \SSYT(\mu',\lambda')$, which has generating function $s_{\mu}(X)K_{\mu'\lambda'}(q,0)$. 
For $A\in\MLQ_0(\mu)$, the preimage under $\rho$ of the set $\{A\} \times \SSYT(\mu',\lambda')$ is $\{B\in\MLQ(\lambda):\rho_N(Q)=A\}$, which has generating function $x^A K_{\mu'\lambda'}(q,0)$. Therefore, to extract the coefficient $[s_{\lambda}]P_{\mu}(X;q,t)$, it is sufficient to sum over the preimage of the set $\{A\} \times \SSYT(\lambda',\mu')$ for any choice of $A\in\MLQ_0(\lambda,n)$. Changing the indices $(\mu',\lambda')\mapsto(\lambda,\mu)$ and picking a representative $A=M(\lambda')\in\MLQ(\lambda',\lambda')$, we obtain the desired formula.
\end{proof}

We record a characterization of the preimage of the multiline queue $M(\lambda)$ under $\rho_N$. Recall that a word is \emph{lattice} if each initial segment contains at least as many letters "$i$" as letters "$i+1$" for every $i\geq 1$.

As an application for double collapsing, \cref{cor:KF_mlq_version_2} gives an alternative formula for $K_{\lambda\mu}(q,0)$ in terms of nonwrapping multiline queues. First, we need some preliminary results.


\begin{lemma}\label{lem:KF-mlqs-rot}
    Let $\lambda$ and $\mu$ be partitions. At $t=0$, the Kostka polynomial is given by
    \[
    K_{\lambda\mu}(q,0)=\sum_{N\in\MLQ_0(\lambda,\rev(\mu))} q^{\maj(\rot^{-1}(N))}.
    \]
\end{lemma}
\begin{proof}
    From \cref{cor:KF_mlq_version_1}, we have that 
    \[
    K_{\lambda\mu}(q,0) = \sum_{\substack{B\in\MLQ(\mu',\lambda') \\ \collapse(B) = M(\lambda')}}q^{\maj(B)}.
    \]
    Let $B\in\MLQ(\mu,\lambda')$ with $\collapse(B) = M(\lambda')$, which occurs if and only if $\rw(B)$ is a lattice word. Consider a particle $b$ in row $i+1$ in $\rot(B).$ By the lattice condition on $\rw(B)$, in $\rot(B)$ there are at least as many particles weakly to the northwest of $b$ in row $i$ as there are in row $i+1$. Since this condition holds for any $b$ in row $i+1$, 
    all particles in row $i+1$ are matched above. This holds for every row $i$, so $\rot(B)$ is nonwrapping. Thus $\rho^\leftarrow_N(\rot(B))=\rot(B)$, and so $\mRSK(B) = (M(\lambda'),\rot(B))$, with $\rot(B) \in \MLQ_0(\lambda,\rev(\mu))$ by construction. The result follows from \cref{lem:maj_left_coll}.
\end{proof}

\begin{remark}
In fact, \cref{lem:KF-mlqs-rot} recovers a result of Nakayashiki-Yamada \cite[Corollary 4.2]{NY95}, which is stated as a sum of the NY energy functions over highest weight elements in the crystal $\mathcal{H}_{\mu}'$.  
For a semistandard tableau $T$, the corresponding highest weight element in $\mathcal{H}_{\mu}'$ is $\flip(N)$, where $N=\mlq(T)$ is the corresponding nonwrapping multiline queue from \cref{thm:bijectionSSYT-MLQ0}. Then the NY energy function is computed on $\rot^{-1}(\flip(N))$ by the decomposition in \eqref{eq:dec}, matching our result according to the previous discussion. 
\end{remark}

Next, we change the direction of collapsing. Let $\rot^2 = \rot \circ \rot$ denote the $180^\circ$ rotation. We think of the following map as "collapsing upwards" with particles paired weakly to the right. 

\begin{prop}\label{prop:full_rotation}
    The following map is a bijection: 
    \[\rho_N^{\uparrow} \, \coloneqq \, \collapse^\downarrow \circ \rot^2\, \colon \, \MLQ_0(\lambda,\alpha) \to \MLQ_0(\lambda,\rev(\alpha)).
    \]
\end{prop}

\begin{proof}
    For $B\in\M$, let $\mRSK'(B) = (\rho^{\uparrow}_N(B),\rho^{\leftarrow}_N(B))$. As mentioned in \cref{rem:directions}, this map is also a bijection. The result follows from the identity $\mRSK' \circ \mRSK^{-1}(B,N) = (\rho^{\uparrow}_N(B),N)$.
\end{proof}

\begin{example}

    We show a nonwrapping multiline queue $Q\in\MLQ_0(\lambda,\alpha)$ for $\lambda = (5,5,4,2)$ and $\alpha = (2,3,2,3,3,3)$ on the left and $\flip(Q)=\collapse(\rot^2(Q))$ on the right, with column content $\rev(\alpha) = (3,3,3,2,3,2).$

    \begin{center}
        \begin{tikzpicture}[scale = 0.95]
            \def \scale{0.3}
            \node at (-2.25,0) {$Q=$};
            \node at (0,0) {\includegraphics[scale = \scale]{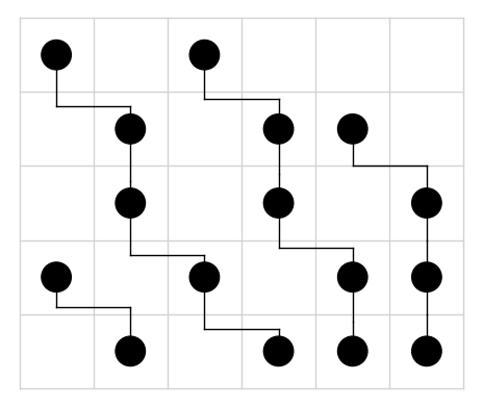}};
            \draw[->]  (2,0)--(3,0) node[above, midway] {$\rot^2$};
            \node at (5,0) {\includegraphics[angle=180,scale = \scale]{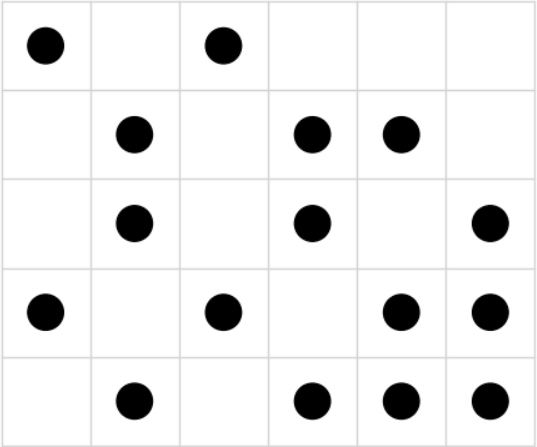}};
            \draw[->]  (7,0)--(8,0) node[above, midway] {$\collapse^\downarrow$};
            \node at (9.75,0) {\includegraphics[angle=180,scale = \scale]{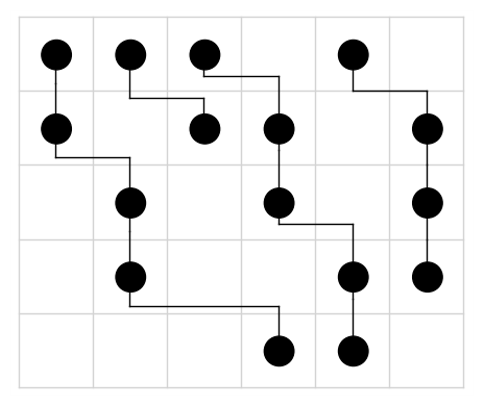}};
            \node at (12.5,0) {$=\flip(Q)$};
        \end{tikzpicture}
    \end{center}
    
\end{example}

Since $\flip$ is the collapsing of a rotation of a multiline queue, the following lemma is a reformulation of \cref{lem:maj_left_coll} in the setting of this new map. 

\begin{lemma}\label{prop:maj_flips}
    Let $\lambda$ and $\mu$ be partition. Then, for any $N\in\MLQ(\lambda,\mu),$ 
    $$\maj(\rot^{-1}(N)) = \maj(\rot(\flip(N))).$$
\end{lemma}

Combining \cref{lem:KF-mlqs-rot} with \cref{prop:maj_flips} we obtain the following result. 

\begin{corollary}\label{cor:KF_mlq_version_2}
    Let $\lambda$ and $\mu$ be partitions. At $t=0$, the Kostka polynomial is given by
    \[
    K_{\lambda\mu}(q,0)=\sum_{N\in\MLQ_0(\lambda,\mu)} q^{\maj(\rot(N))}.
    \]
\end{corollary}

Finally, in view of \cref{thm:mlqs-to-mlds,cor:KF_mlq_version_2} we give a bosonic multiline queue formula for $K_{\lambda\mu}(q,0)$. The proof follows the same arguments as the one for \cref{cor:KF_mlq_version_1}. 

\begin{corollary}\label{cor:KF_mld_version_1}
    Let $\lambda$ and $\mu$ be partitions. Then 
    \[
    \widetilde{K}_{\lambda\mu}(q,0) = \sum_{\substack{D\in\MLD(\mu,\lambda') \\ \widetilde\collapse(D) = D(\lambda')}}q^{\widetilde\maj(D)}.
    \]
    where $D(\lambda')$ is the anti-diagonal bosonic multiline queue with $\lambda'_i$ particles in row and column $i$.  
\end{corollary}

Moreover, from the change of pairing direction in \cref{thm:double-collapsing-diagrams}, using the rotation map $\rot$, we obtain directly the following alternate formula for $K_{\lambda,\mu}(q,0)$. 

\begin{corollary}\label{cor:KF_mld_version_2}
    Let $\lambda$ and $\mu$ be partitions. The Kostka--Foulkes polynomial is given by
    \[
    K_{\lambda\mu}(q,0)=\sum_{D\in\MLD_{0,R}(\lambda',\mu)} q^{\widetilde\maj(\rot(D))}
    \]
    where we think of $\rot(D)$ as a bosonic multiline queue pairing strictly to the left.
\end{corollary}

\section*{Acknowledgements}

Both authors were partially supported by NSERC grant RGPIN-2021-02568. We thank Travis Scrimshaw, Kartik Singh, and Santiago Estupi\~n\'an Salamanca for helpful conversations and insightful observations. We are also grateful to the anonymous referee for their useful comments and for pointing out several references on the existing results presented in this work. 

\printbibliography

\end{document}